**T.C.**

**GEBZE YÜKSEK TEKNOLOJİ ENSTİTÜSÜ**

**FEN FAKÜLTESİ MATEMATİK BÖLÜMÜ**

# FİBONACCİ SAYILARI, ALTIN ORAN VE UYGULAMALARI

**Erdoğan ŞEN**

**LİSANS BİTİRME ÖDEVİ**

**MATEMATİK ANABİLİM DALI**

T.C.

GEBZE YÜKSEK TEKNOLOJİ ENSTİTÜSÜ

FEN FAKÜLTESİ MATEMATİK BÖLÜMÜ

# FİBONACCİ SAYILARI, ALTIN ORAN VE UYGULAMALARI

Erdoğan ŞEN

LİSANS BİTİRME ÖDEVİ

MATEMATİK ANABİLİM DALI

DANIŞMAN

Yrd. Doç. Dr. Mustafa AKKURT

Bu çalışma 17/06/2008 tarihinde aşağıdaki jüri tarafından Matematik Bölümünde Lisans Bitirme Projesi olarak kabul edilmiştir.

Bitirme Projesi Jürisi

| | | |
|---|---|---|
| Danışman Adı | Yard. Doç. Dr. Mustafa AKKURT | |
| Üniversite | Gebze Yüksek Teknoloji Enstitüsü | |
| Fakülte | Fen Fakültesi | |
| Jüri Adı | Yard. Doç. Dr. Coşkun YAKAR | |
| Üniversite | Gebze Yüksek Teknoloji Enstitüsü | |
| Fakülte | Fen Fakültesi | |
| Jüri Adı | Prof. Dr. Tahir AZEROĞLU | |
| Üniversite | Gebze Yüksek Teknoloji Enstitüsü | |
| Fakülte | Fen Fakültesi | |

# ÖZET

Bu tezde Fibonacci sayılarının matematiksel özelliklerini, doğadaki, geometrideki ve ekonomideki uygulamalarını inceledik.

Altın Oran'ı elde ederek, bazı matematiksel özdeşlikleri ispatlamak için kullandık. Asal sayıların sonsuzluğu da Fibonacci sayıları kullanılarak ispatlandı. Doğadaki Fibonacci sayıları ile olan karşılaşmalar örnekleriyle detaylı olarak incelendi. Ayrıca borsa hareketleri ile Fibonacci sayıları arasındaki ilişkiye dair örnekler verildi ve bunlar incelendi.

# SUMMARY

In this thesis we examined mathematical properties of Fibonacci numbers and applications of this numbers in the nature,geometry and economy.

We obtained Golden section and proved some mathematical identities using Golden section. Infinity of the prime numbers proved by using properties of Fibonacci numbers. Encounterings with Fibonacci numbers in the nature are examined with details. Also examples are given for relation about Fibonacci numbers and stock exchange and these are examined.

# TEŞEKKÜR

Bu bitirme çalışmasının ortaya çıkmasındaki eşsiz katkılarından dolayı çok değerli danışmanım Sayın Yrd. Doç. Dr. Mustafa AKKURT’a ve aileme çok teşekkür ederim.

# İÇİNDEKİLER DİZİNİ

Sayfa No:

# 1 GİRİŞ

Fibonacci sayıları günümüzde matematik alanında en çok ilgi çeken konuların başında gelir. Bunun nedeni Fibonacci sayılarının keşfedilmeye uygun yapısıdır. Geometride, doğada, ekonomide ve sanatta bir çok uygulamaları bulunmaktadır. Fibonacci sayılarına özellikle doğada çok sık rastlamaktayız. Keşfedilmeyi bekleyen bir çok özelliği bulunan, keşfedilen fakat ispatlanamayan da birçok özelliği bulunan Fibonacci dizisi ayrıca daha önce ispatlanan bazı matematiksel özelliklere yeni ispat yolları vermektedir(asal sayıların sonsuz olduğunu gösterirken olduğu gibi).

# 2 FİBONACCİ DİZİSİNİN MATEMATİKSEL ÖZELLİKLERİ VE ALTIN ORAN

$u_1 = 1$ , $u_2 = 1$ başlangıç değerleri ve $u_n = u_{n-1} + u_{n-2}$ indirgeme bağıntısı ile tanımlı diziye Fibonacci dizisi denir. İndirgeme bağıntısı, $u_n - u_{n-1} - u_{n-2} = 0$ ve karakteristik denklem $\lambda^2 - \lambda - 1 = 0$ olmak üzere karakteristik denklemin kökleri, $\lambda_1 = \frac{1+\sqrt{5}}{2}$ve $\lambda_2 = \frac{1-\sqrt{5}}{2}$ olmak üzere dizinin genel terimi,

$$u_n = A_1 \left( \frac{1+\sqrt{5}}{2} \right)^n + A_2 \left( \frac{1-\sqrt{5}}{2} \right)^n$$

biçimindedir. Başlangıç değerlerinden

$$A_1 \frac{1+\sqrt{5}}{2} + A_2 \frac{1-\sqrt{5}}{2} = 1$$

$$A_1 \left( \frac{1+\sqrt{5}}{2} \right)^2 + A_2 \left( \frac{1-\sqrt{5}}{2} \right)^2 = 1$$

denklem sistemi yazılır. Denklem sisteminin çözülmesiyle

$$A_1 = \frac{1}{\sqrt{5}}, A_2 = -\frac{1}{\sqrt{5}}$$

elde edilir. Buna göre Fibonacci dizisinin genel terimi,

$$u_n = \frac{1}{\sqrt{5}} \left[ \left( \frac{1+\sqrt{5}}{2} \right)^n - \left( \frac{1-\sqrt{5}}{2} \right)^n \right]$$

olur. Fibonacci dizisinin ilk 20 terimi şu şekildedir:

| $n$ | [illegible] | [illegible] | $n$ | [illegible] | [illegible] |
|---|---|---|---|---|---|
| 1 | 1 | | 11 | 89 | 1.618182 |
| 2 | 1 | 1 | 12 | 144 | 1.617977 |
| 3 | 2 | 2 | 13 | 233 | 1.618056 |
| 4 | 3 | 1.5 | 14 | 377 | 1.618026 |
| 5 | 5 | 1.666667 | 15 | 610 | 1.618037 |
| 6 | 8 | 1.6 | 16 | 987 | 1.618033 |
| 7 | 13 | 1.625 | 17 | 1597 | 1.618034 |
| 8 | 21 | 1.615385 | 18 | 2584 | 1.618034 |
| 9 | 34 | 1.619048 | 19 | 4181 | 1.618034 |
| 10 | 55 | 1.617647 | 20 | 6765 | 1.618034 |

Leonardo Fibonacci 12-13 üncü yüzyıllarda yaşamış bir İtalyan matematikçisidir. Pisa şehrinde doğan Leonardo çocukluğunu babasının çalışmakta olduğu Cezair'de geçirmiştir. İlk matematik bilgilerini müslüman eğiticilerden almış olup küçük yaşlarda onluk Arap sayı sistemini öğrenmiştir. Ülkesi İtalya'da kullanılmakta olan Roma sisteminin hantallığı yanında Arap sisteminin mükemelliğini gören Fibonacci 1201 yılında "Liber abaci" isimli kitabını yazmıştır. Aritmetik ve cebir içeren ticaret ile ilgili bu kitapta Arap sayı sisteminin tanıtımını ve müdafasını yapmıştır. İlk anda kitabın İtalyan tüccarları üzerinde etkisi az olmasına rağmen zamanla bu kitap arab sayı sisteminin batı Avrupa'ya girmesinde büyük rol oynamıştır. Bu kitapta bulunan bir problem ortaçağ matematiğine katkıları olan Fibonacci'yi altı yüzyıl sonra, ondokuzuncu yüzyılın başlarından günümüze meşhur hale gelmesine sebep olmuştur. Bu problem "tavşan problemi" dir. Ergin bir tavşan çiftinin her ay yeni bir yavru çifti verdikleri ve yeni doğan çiftin bir ay zarfında tam erginliğe eriştikleri varsayımıyla, yavru olan bir tavşan çiftinden başlayıp bir yılda (11 ayda) çiftlerin sayısı ne olur?Aylar : 1 2 3 4 5 6 7 8 9 10 11 Çiftlerin sayısı: 1 1 2 3 5 8 13 21 34 55 89.Kitaplara tavşan problemi olarak geçen bu problem halkımız tarafından kuzu-toklu-koyun problemi olarak bilinmektedir. Bilindiği gibi ergin bir kuzu olan toklu genellikle iki yaşına geldiğinde yavrulamaktadır. Buna göre bir yavru dişi kuzu ile başlayıp her yavrulama sonucu bir dişi kuzunun doğduğu ve ölüm olmaması durumunda 11 yıl sonunda kaç baş hayvana sahip olunur? Bu gelişimin ilk beş yılını şekil üzerinde gözleyelim. Dişi kuzuyu k, tokluyu T, koyunu K harfi temsil etmek üzere gelişim aşağıdaki gibidir.

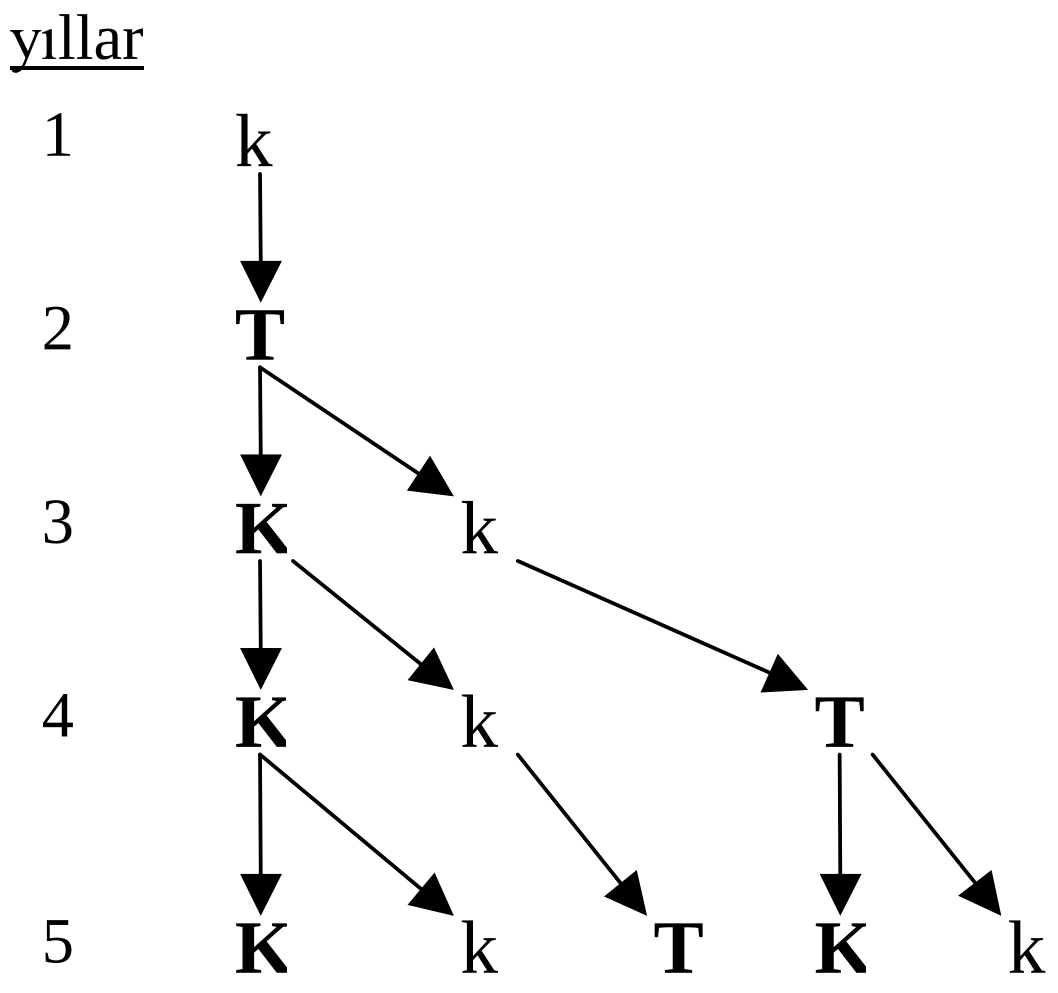

Belli bir yıldaki hayvan sayısı önceki iki yıldakilerin toplamı olmak üzere 11. yılda 89 baş hayvana ulaşılacaktır. Bunların yıllara göre yavru, toklu, koyun sayısı dağılışına gelince şekildeki gibi bir durum ortaya çıkmaktadır.

| yıllar | kuzu sayısı | toklu sayısı | koyun sayısı | toplam |
|---|---|---|---|---|
| 1 | 1 | - | - | 1 |
| 2 | 0 | 1 | - | 1 |
| 3 | 1 | 0 | 1 | 2 |
| 4 | 1 | 1 | 1 | 3 |
| 5 | 2 | 1 | 2 | 5 |
| 6 | 3 | 2 | 3 | 8 |
| 7 | 5 | 3 | 5 | 13 |
| 8 | 8 | 5 | 8 | 21 |
| 9 | 13 | 8 | 13 | 34 |
| 10 | 21 | 13 | 21 | 55 |
| 11 | 34 | 21 | 34 | 89 |

Dikkat edilirse $b_n$ $n.$ yıldaki kuzuların sayısı olmak üzere, kuzular için

$$b_n = b_{n-1} + b_{n-2}, b_1 = 1, b_2 = 0$$

indirgeme bağıntısı sözkonusudur ve çözüm olan 1 , 0 , 1 , 1 , 2 , 3 , 5 , 8 , 13 , . . . dizisinde $n \geq 3$ için Fibonacci dizisi ortaya çıkmaktadır, yani

$$b_n = u_{n-2}, n \geq 3$$

dır. Benzer durum yıllara göre toklu sayısı için de sözkonusudur. $c_n$ $n$.yıldaki tokluların sayısı olmak üzere,

$$c_n = u_{n-3} \quad , n \geq 3$$

dır. Yıllara göre koyunların sayısı ise iki yıl geçikmeli bir Fibonacci dizisidir. Fibonacci'nin kendisi 1 , 1 , 2 , 3 , 5 , 8 , 13 , 21 , ... dizisi üzerinde bir inceleme yapmamıştır. Hatta bu dizi üzerinde ondokuzuncu yüzyılın başlarına kadar ciddi bir araştırma yapılmadığı da belirtilmektedir. Ancak bundan sonra bu dizi üzerine yapılan araştırmaların sayısı Fibonacci'nin tavşanlarının sayısı gibi artmıştır.Fibonacci dizisinde olduğu gibi,

$$u_n = u_{n-1} + u_{n-2}$$

indirgeme bağıntılı ve herhangi $u_1, u_2$ başlangıç değerli genelleştirilmiş Fibonacci dizileri,

$$u_n = u_{n-1} + u_{n-2} + u_{n-3}$$

indirgeme bağıntılı $u_1, u_2, u_3$ başlanğıç değerli tribonacci dizileri gibi genellemeler yapılmıştır. Fibonacci Derneği tarafından 1963 yılından itibaren yayınlanan "The Fibonacci Quarterly" dergisi ilginç özelliklere sahip tamsayıları ve özellikle genelleştirilmiş Fibonacci sayılerını inceleyen araştırmalar yayınlamaktadır.Bu dizilerin çekiciliği bir taraftan matematiği sevenlerin çok az bir temel bilgiyle bu dizileri inceleme imkanı bulmaları diğer taraftan ise hayatta ve araştırmalarda hiç umulmadık yerlerde bu dizilerle karşılaşmalarından kaynaklanmaktadır. Bazısı bilinen, bazısı öne sürülüp ispatlanamayan ve bilinmeyip keşfedilmesi beklenen bir çok özelliğe sahiptir.Genelleştirilmiş Fibonacci dizilerinde de geçerli olmak üzere, $u_1 = 1, u_2 = 1$ ve $u_n = u_{n-1} + u_{n-2}$ olan Fibonacci dizisinin bir terimi öncekine bölündüğünde bölümün $n \to \infty$ için,"altın oran" denen ve irrasyonel bir sayı olan $\frac{1+\sqrt{5}}{2} = 1,61803398...$ sayısına yakınsadığı görülmektedir.

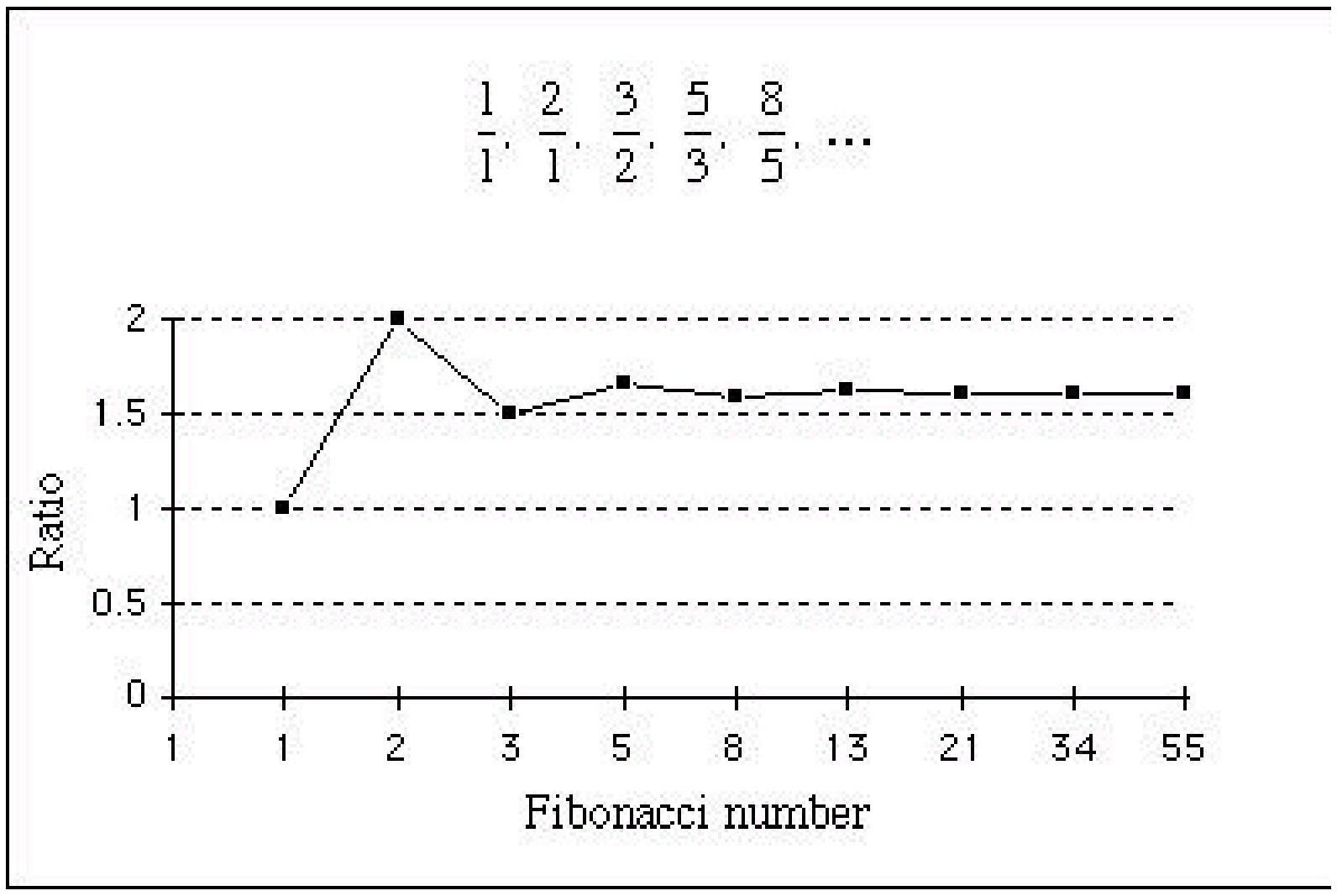

Gerçekten, $n \to \infty$ için limite geçersek

$$\begin{aligned}
\lim_{n\to\infty} \frac{u_{n+1}}{u_n} &= \lim_{n\to\infty} \frac{\frac{1}{\sqrt{5}}\left[\left(\frac{1+\sqrt{5}}{2}\right)^{n+1} - \left(\frac{1-\sqrt{5}}{2}\right)^{n+1}\right]}{\frac{1}{\sqrt{5}}\left[\left(\frac{1+\sqrt{5}}{2}\right)^{n} - \left(\frac{1-\sqrt{5}}{2}\right)^{n}\right]} \\
&= \lim_{n\to\infty} \frac{\left(\frac{1+\sqrt{5}}{2}\right)^{n+1}\left[1 - \left(\frac{1-\sqrt{5}}{1+\sqrt{5}}\right)^{n+1}\right]}{\left(\frac{1+\sqrt{5}}{2}\right)^{n}\left[1 - \left(\frac{1-\sqrt{5}}{1+\sqrt{5}}\right)^{n}\right]} \\
&= \lim_{n\to\infty} \frac{\left(\frac{1+\sqrt{5}}{2}\right)^{n+1}}{\left(\frac{1+\sqrt{5}}{2}\right)^{n}} \lim_{n\to\infty} \frac{\left[1 - \left(\frac{1-\sqrt{5}}{1+\sqrt{5}}\right)^{n+1}\right]}{\left[1 - \left(\frac{1-\sqrt{5}}{1+\sqrt{5}}\right)^{n}\right]} \\
&= \frac{1+\sqrt{5}}{2}
\end{aligned}$$

olur. Buradan Fibonacci dizisinin ıraksak bir dizi olduğunu söyleyebiliriz. Şöyle ki $n \to \infty$ için $L = \lim \frac{u_{n+1}}{u_n} = \frac{1+\sqrt{5}}{2} \approx 1.61803 > 1$ ve D'alambert testine göre $L > 1$ olduğundan ıraksaktır.

Polinomların reel köklerine cebirsel sayı denir. İrasyonel bir sayı olan $\Phi = \frac{1+\sqrt{5}}{2}$ sayısı $x^2-x-1=0$ veya bu denklemin $1/x$ ile çarpılmış hali olan $x = 1/x+1$ denkleminin pozitif değerli $\lambda_1$ kökü olduğundan bir cebirsel sayıdır. $\Phi^2 = \Phi+1$ olduğundan

$$\begin{aligned}\Phi &= \sqrt{1+\Phi} = \sqrt{1+\sqrt{1+\Phi}} = \sqrt{1+\sqrt{1+\sqrt{1+\Phi\ldots}}} \\ &= \sqrt{1+\sqrt{1+\sqrt{1+\ldots}}}\end{aligned}$$

yazılır. Ayrıca

$$\begin{aligned}\Phi &= 1+\frac{1}{\Phi} = 1+\frac{1}{1+\frac{1}{\Phi}} = 1+\frac{1}{1+\frac{1}{1+\frac{1}{\Phi}}} \\ &= 1+\frac{1}{1+\frac{1}{1+\frac{1}{1+\ldots}}}\end{aligned}$$

olmak üzere $\Phi$ sayısı elemanları 1 olan bir sonsuz zincir kesiridir. Bu zincir kesirinin kısmi kesirleri,

$$1+\frac{1}{1} = \frac{2}{1}$$

$$1+\frac{1}{1+\frac{1}{1}} = \frac{3}{2}$$

$$1+\frac{1}{1+\frac{1}{1+\frac{1}{1}}} = \frac{5}{3}$$

$$1+\frac{1}{1+\frac{1}{1+\frac{1}{1+\frac{1}{1}}}} = \frac{8}{5}$$

olmak üzere ardışık iki Fibonacci sayısının oranıdır. $b_0 = 1, b_1 = \Phi$ başlangıç değerleri ve $b_n = b_{n-1} + b_{n-2}$ indirgeme bağıntısı ile tanımlanan genelleştirilmiş Fibonacci dizisi$1, 1+\Phi, 1+2\Phi, 2+3\Phi, 3+5\Phi, \ldots$ olmak üzere, $\Phi^2 = \Phi+1$ bağıntısının gözönüne alınmasıyla, $1, \Phi, \Phi^2, \Phi^3, \Phi^4, \ldots$ biçiminde yazılabilir. Görüldüğü gibi bu genelleştirilmiş Fibonacci dizisi aynı zamanda ortak çarpanı $\Phi$ olan bir geometrik dizidir. Bu diziye Altın Dizi denmektedir. Altın Dizinin terimleri $1, \frac{1+\sqrt{5}}{2}, \frac{3+\sqrt{5}}{2}, \frac{4+2\sqrt{5}}{2}, \ldots$ olarak yazıldığında kesirlerin paylarındaki birinci terimler $1, 3, 4, 7, 11, \ldots$ gibi genelleştirilmiş bir Fibonacci dizisi, ikinci terimlerdeki $\sqrt{5}$ in önündeki katsayılar ise $1, 1, 2, 3, 5, \ldots$ Fibonaci dizisini oluşturmaktadır.Fibonacci dizisinin ilk 20 teriminin toplamını inceleyelim.İlk n terimin toplamı

$$\begin{aligned}
u_1+u_2+\cdots+u_n &= \sum_{k=1}^{n}\frac{1}{\sqrt{5}}\left[\left(\frac{1+\sqrt{5}}{2}\right)^k-\left(\frac{1-\sqrt{5}}{2}\right)^k\right] \\
&= \frac{1}{\sqrt{5}}\left[\sum_{k=1}^{n}\left(\frac{1+\sqrt{5}}{2}\right)^k-\sum_{k=1}^{n}\left(\frac{1-\sqrt{5}}{2}\right)^k\right] \\
&= \frac{1}{\sqrt{5}}\left[\frac{\left(\frac{1+\sqrt{5}}{2}\right)\left(1-\left(\frac{1+\sqrt{5}}{2}\right)^n\right)}{1-\left(\frac{1+\sqrt{5}}{2}\right)}-\frac{\left(\frac{1-\sqrt{5}}{2}\right)\left(1-\left(\frac{1-\sqrt{5}}{2}\right)^n\right)}{1-\left(\frac{1-\sqrt{5}}{2}\right)}\right] \\
&= u_{n+2}-1
\end{aligned}$$

olmak üzere,

$$u_{n+2}=u_1+u_2+u_3+...+u_n+1$$

yazılabilir. Buna göre ilk 20 terimin toplamı,

$$u_1+u_2+u_3+...+u_{20}=u_{22}-1=17711-1=17710$$

olur.$b_1=\alpha, b_2=\beta$ başlangıç değerli genelleştirilmiş Fibonacci dizisinin terimleri için,

$$\begin{aligned}
b_3 &= b_1+b_2=\alpha+\beta \\
b_4 &= b_2+b_3=\alpha+2\beta \\
b_5 &= b_3+b_4=2\alpha+3\beta \\
b_6 &= b_4+b_5=3\alpha+5\beta \\
&\vdots
\end{aligned}$$

yazılabilir. $a_n$ Fibonacci dizisi olmak üzere, $b_n=\alpha.a_{n-2}+\beta.a_{n-1}$ olur. Buna göre,

$$\begin{aligned}
b_1+b_2+...+b_n &= \alpha(1+a_1+a_2+...+a_{n-2})+\beta(1+a_1+a_2+...+a_{n-1}) \\
&= \alpha a_n+\beta\left(a_{n+1}-1\right)
\end{aligned}$$

elde edilir. $b_1=1, b_2=3$ başlangıç değerli 1 , 3 , 4 , 7 , 11 , 18 , 29 , 47

, 76 ,... genelleştirilmiş Fibonacci dizisinin ilk 15 teriminin toplamı, $b_1 + b_2 + ... + b_{15} = 1.a_{15} + 3.(a_{16} - 1)$ olmak üzere bu toplam 3568 dir.Genelleştirilmiş bir Fibonacci dizisinde ilk 10 terimin toplamı 7. terimin 11 katıdır. Örneğin $b_1 = 1, b_2 = 3$ başlangıç değerli genelleştirilmiş Fibonacci dizisinde, $b_1 + b_2 + ... + b_{10} = 1.a_{10} + 3.(a_{11} - 1) = 55 + 3(89 - 1) = 319$ ve $b_7 = 29$ olmak üzere $b_1 + b_2 + ... + b_{10} = 11.b_{17}$ dir. Genel olarak $b_1 = \alpha, b_2 = \beta$ başlangıç değerli bir genelleştirilmiş Fibonacci dizisinde, $b_7 = 5\alpha + 8\beta$ olmak üzere

$$b_1 + b_2 + ... + b_{10} = \alpha.a_{10} + \beta(a_{11} - 1) = \alpha.55 + \beta(89 - 1) = 55.\alpha + 88.\beta = 11.b_7$$

olur. Fibonacci sayıları hızlı bir şekilde büyür. $u_{5n+2} > 10^n, n \geq 1$ eşitsizliği de bunu gösteriyor. Örneğin $n = 1$ için $u_7 > 10$, $n = 2$ için $u_{12} > 100$, $n = 3$ için $u_{17} > 1000$, $n = 4$ için $u_{22} > 10000$ .... Eşitsizlik $n$ için tümevarım kullanılarak gerçeklenebilir. $n = 1$ için aşikardır. Çünkü $u_7 = 13 > 10$ olur. Şimdi varsayalım ki eşitsizlik herhangi bir $n$ sayısı için sağlanıyor. $n+1$ için sağlandığını gösterelim. Rekürsiyon bağıntısı $u_k = u_{k-1} + u_{k+1}$ uygulayalım.

$$\begin{aligned} u_{5n+7} &= 8u_{5n+2} + 5u_{5n+1} \\ &> 8u_{5n+2} + 2(u_{5n+1} + u_{5n}) \\ &= 10u_{5n+2} > 10.10^n = 10^{n+1} \end{aligned}$$

**Teorem 1.** Fibonacci dizisi için $(u_n, u_{n+1}) = 1, n \geq 1$ .

**İspat.** Varsayalım ki $d > 1$ sayısı $u_n$ ve $u_{n+1}$ i böler. O halde onların farkı olan $u_{n+1} - u_n = u_{n-1}$ sayısı $d$ tarafından bölünebilir. Buradan ve $u_n - u_{n-1} = u_{n-2}$ bağıntısından $d \mid u_{n-2}$ ye ulaşılabilir. Benzer yöntemi tekrar uygularsak görürüz ki $d \mid u_{n-3}, d \mid u_{n-4}, ...$ ve sonunda $d \mid u_1$.fakat $u_1 = 1$ herhangi bir $d > 1$ tarafından bölünemez. Böylece çelişki elde edildi ve ispat tamamlandı. ■

$u_3 = 2, u_5 = 5, u_7 = 13$ ve $u_{11} = 89$ asal sayılardır. Buradan $n > 2$ için $u_n$ in bir asal olduğunu tahmin edersek yanılgıya düşeriz. Çünkü $u_{19} = 4181 = 37.113$ ve $u_{19}$ asal değildir. Asal Fibonacci sayılarının sonsuz tane olup olmadığı bilinmemektedir. Fakat bir p asalı için p tarafından bölünen sonsuz tane Fibonacci sayısı vardır. Örnek verilecek olursa 3 her dördüncü terimi, 5 her beşinci terimi ve 7 her sekizinci terimi böler.

$u_1, u_2, u_6, u_{12}$ dışındaki her Fibonacci sayısı yeni bir asal çarpana sahiptir. Şöyle ki; bu yeni asal çarpan daha küçük bir alt indise sahip Fibonacci sayısında bulunmaz. Örneğin 29, $u_{14} = 377 = 13.29$ sayısını böler fakat ondan önceki hiçbir Fibonacci sayısını bölmez.

Bildiğimiz gibi iki pozitif sayının en büyük ortak böleni öklid algoritmasından sonlu sayıda bölme işleminin ardından bulunabilir. Katsayılar uygun seçilerek gerekli olan bölenlerin sayısı yeterince geniş tutulabilir. Kesin durum şudur: Verilen bir $n > 0$ sayısı için $a$ ve $b$ tamsayılarının en büyük ortak bölenini hesaplamak demek; öklid algoritmasında $n$ tane bölene ihtiyaç var demektir. Bunu göstermek için; $a = u_{n+2}$ ve $b = u_{n+1}$ sayılarının en büyük ortak bölenini öklid algoritmasını kullanarak bulmak için aşağıdaki denklem sistemini kuralım.

$$\begin{aligned} u_{n+2} &= 1.u_{n+1} + u_n \\ u_{n+1} &= 1.u_n + u_{n-1} \\ &\vdots \\ u_4 &= 1.u_3 + u_2 \\ u_3 &= 2.u_2 + 0 \end{aligned}$$

aşikar olarak gerekli olan bölenlerin sayısı $n$'dir. Algoritmada son sıfır olmayan kalanın $(u_{n+2}, u_{n+1})$ değerini sağladığını hatırlarsak $(u_{n+2}, u_{n+1}) = u_2 = 1$ olur ki bu da ardışık Fibonacci sayılarının aralarında asal olduğunu gösterir.

**Örnek**. $u_8 = 21$ ve $u_7 = 13$ sayılarının en büyük ortak bölenini bulmak için 6 bölene ihtiyacımız vardır.

$$\begin{aligned} 21 &= 1.13 + 8 \\ 13 &= 1.8 + 5 \\ 8 &= 1.5 + 3 \\ 5 &= 1.3 + 2 \\ 3 &= 1.2 + 1 \\ 2 &= 2.1 + 0 \end{aligned}$$

Gabriel Lame 1844 yılında $(a, b)$ yi hesaplamak için öklid algoritmasında $n$ tane bölme işlemi gerektiğini elde etti($a > 0, b > 0$ ve $a \geq u_{n+2}, b \geq u_{n+1}$ olmak üzere). Sonuç olarak belirli bir zamana kadar $u_n$ dizisine Lame dizisi denildi. Lucas Fibonaccinin bunu altı yüzyıl önce bulduğunu keşfetti ve Amerikan Matematik Bülteninin 1878 yılındaki açılış baskısında Fibonacci dizisi adıyla yayınladı.

Fibonacci dizilerinin dikkat çeken özelliklerinden biri de iki Fibonacci dizisinin en büyük ortak böleninin yine bir Fibonacci sayısı olmasıdır.

$$u_{m+n} = u_{m-1}u_n + u_m u_{n+1} \tag{2.1}$$

eşitliği de bu gerçeği gözler önüne seriyor. Sabit $m \geq 2$ için, bu eşitlik $n$ üzerine kurulu $n = 1$ alırsak 2.1 denklemi

$$u_{m+1} = u_{m-1}u_1 + u_m u_2 = u_{m-1} + u_m$$

şeklini alır. Varsayalım ki $n = k$ için doğru olsun. $n = k+1$ için doğru olduğunu gösterelim.

$$u_{m+k} = u_{m-1}u_k + u_m u_{k+1}$$

$$u_{m+(k-1)} = u_{m-1}u_{k-1} + u_m u_k$$

Bu iki denklemin toplamı bize;

$$u_{m+k} + u_{m+(k-1)} = u_{m-1}(u_k + u_{k-1}) + u_m(u_{k+1} + u_k)$$

denklemini verir. Fibonacci sayılarının tanımından

$$u_{m+(k+1)} = u_{m-1}u_{k+1} + u_m u_{k+2}$$

olur ki bu da 2.1 denkleminde $n$ yerine $k+1$ yazılmış durumdur. Böylece tümevarım basamağı tamamlandı. Yani 2.1 denklemi her $m \geq 2$ ve $n \geq 1$ için sağlanır. Bir örnekle gösterecek olursak;

$$u_9 = u_{6+3} = u_5.u_3 + u_6 u_4 = 5.2 + 8.3 = 34$$

**Teorem 2.** $m \geq 1, n \geq 1$ için $u_{mn}, u_m$ tarafından bölünebilir.

**İspat.** Bir kez daha tümevarımı kullanacağız. $n = 1$ için aşikardır. Şimdi $n = k$ için doğru olduğunu varsayalım. $n = k+1$ için doğru olduğunu gösterelim.

$$u_{m(k+1)} = u_{mk+m}$$

durumuna geçiş 2.1 denklemini kullanarak gerçekleşir. Gerçekten;

$$u_{m(k+1)} = u_{mk-1}u_m + u_{mk}u_{m+1}$$

çünkü altpozisyondan $u_m$ $u_{mk}$ yı böler. Bu açıklamanın sağ kısmı(ve burada sol kısmı) $u_m$ tarafından bölünmeli. Böylece $u_m \mid u_{m(k+1)}$ olur ki bu da ispatı tamamlar. ■

Aşağıda $(u_m, u_n)$ değerini hesaplamak için hazırlayıcı bir teknik lemma verelim.

**Lemma 1.** Eğer $m = nq + r$ ise $(u_m, u_n) = (u_r, u_n)$

**İspat.** 2.1 denklemini kullanarak;

$$(u_m, u_n) = (u_{qn+r}, u_n) = (u_{qn-1}u_r + u_{qn}u_{r+1}, u_n)$$

yazabiliriz. Ayrıca Teorem 2 ve $(a + c, b) = (a, b)$ özelliğine göre $b \mid c$ olduğu zaman

$$(u_{qn-1}.u_r + u_{qn}.u_{r+1}, u_n) = (u_{qn-1}.u_r, u_n)$$

İddiamız $(u_{qn-1}, u_n) = 1$ bunu görmek için $d$ yi $d = (u_{qn-1}, u_n)$ olacak şekilde kuralım. Bağıntılar $d \mid u_n, u_n \mid u_{qn}, d \mid u_{qn}$ olmasını gerektirir. Bu yüzden $d$ ardışık Fibonacci sayıları olan $u_{qn-1}$ ve $u_{qn}$ in bir pozitif ortak bölenidir. Ardışık Fibonacci sayıları kendi aralarında asal olduğundan $d = 1$ olmalıdır.

$(a, c) = 1$ ise $(a, bc) = (a, b)$ olduğundan $(u_m, u_n) = (u_{qn-1}u_r, u_n) = (u_r, u_n)$ elde edilir. ■

**Teorem 3.** İki Fibonacci sayısının en büyük ortak böleni de bir Fibonacci sayısıdır. $(u_m, u_n) = u_d$ ve $d = (m, n)$.

**İspat.** Varsayalım ki $m \geq n$ . $m$ ve $n$ sayılarına öklid algoritmasını uygulayalım. Aşağıdaki denklem sistemini elde ederiz.

$$\begin{aligned} m &= q_1 n + r_1 && 0 < r_1 < n \\ n &= q_2 r_1 + r_2 && 0 < r_2 < r_1 \\ r_1 &= q_3 r_2 + r_3 && 0 < r_3 < r_2 \\ &\vdots \\ r_{n-2} &= q_n r_{n-1} + r_n && 0 < r_n < r_{n-1} \\ r_{n-1} &= q_{n+1} r_n + 0 \end{aligned}$$

Lemma 1'den $(u_m, u_n) = (u_{r_1}, u_n) = (u_{r_1}, u_{r_2}) = ... = (u_{r_{n-1}}, u_{r_n})$ ve $r_n \mid r_{n-1}$ olduğundan $u_{r_n} \mid u_{r_{n-1}}$elde ederiz. O halde $u_{r_n} \mid u_{r_{n-1}}$ olduğundan $(u_{r_{n-1}}, u_{r_n}) = u_{r_n}$. Ayrıca $m$ ve $n$ için uygulanan öklid algoritmasındaki son sıfır olmayan kalan olan $r_n$ $m$ ve $n$ sayılarının en büyük ortak bölenine eşittir. Bütün bunları birleştirirsek $(u_m, u_n) = u_{(m,n)}$ elde ederiz. ■

Teorem 2'nin tersi yukarıda ispatlanan teoremden elde edilir. Başka bir deyişle; eğer $u_n$ $u_m$ tarafından bölünebiliyorsa o halde $m$ tarafından bölünür. Gerçekten eğer $u_m \mid u_n$ ise $(u_m, u_n) = u_m$ olur. Fakat Teorem 3'e göre $(u_m, u_n)$ değeri $u_{(m,n)}$ değerine eşittir. Bütün bunları aşağıdaki sonuçla özetleyelim.

**Sonuç 1.** Fibonacci dizisinde $n \geq m \geq 3$ için $u_m \mid u_n$ olması için gerek ve yeter koşul $m \mid n$ olmasıdır.

Teorem 3 için bir örnek verelim.

**Örnek.**

$$(u_{16}, u_{12}) = (987, 144) = 3 = u_4 = u_{(16,12)}$$

Teorem 3'ten alt indis $n > 4$ kompositse(bölünebilirse,asal değilse) o halde $u_n$ de kompozittir. $n = rs$ için Sonuç $u_r \mid u_n$ ve $u_s \mid u_n$ olmasını gerektirir. Örnek vermek gerekirse $u_4 \mid u_{20}$ ve $u_5 \mid u_{20}$, farklı şekilde ifade etmek gerekirse 3 ve 5 6765 sayısını böler. Böylece asallar Fibonacci dizisinde sadece asal alt indisler için yer alır( $u_2 = 1$ ve $u_4 = 3$ dışında). Fakat $p$ asalsa $u_p$ asal olmayabilir. Örneğin; $u_{19} = 37.113$ sayısında olduğu gibi.

Şimdi Fibonacci sayılarını kullanarak asalların sonsuzluğunu ispatlayalım.

**Teorem 4.** Asal sayılar sonsuz tanedir.

**İspat.** Sonlu tane asal sayı olduğunu varsayalım. Bunların sayısına $r$ diyelim. Asallar $2, 3, 5, ...p_r$ şeklinde gittikçe artan şekilde sıralansın. Bunlara karşılık gelen Fibonacci sayıları olan $u_2, u_3, u_5, ...u_{p_r}$ sayılarını göz önüne alalım.Teorem 3'e göre $u_2 = 1$ dışında bunlar ikişer ikişer aralarında asaldır. Çünkü bu asal sayılardan alınan $a$ ve $b$ gibi iki asal sayı için $a \nmid b$ olduğundan $u_a \nmid u_b$. Kalan $r-1$ sayı, iki asal çarpana sahip bir sayı dışında hepsi tek bir asal çarpana sahip. Çelişki elde edildi. Çünkü $u_{37} = 73 \cdot 149 \cdot 2221$ üç tane asal çarpana sahip. ■

## 2.1 FİBONACCİ SAYILARIYLA ELDE EDİLEN ÖZDEŞLİKLER

Burada basit özdeşlikleri Fibonacci sayılarını dahil ederek geliştireceğiz. En basitlerinden biri olan, ilk $n$ Fibonacci sayısının toplamının $u_{n+2} - 1$ sayısına olan eşitliğini vererek başlayalım.

Örnek ile gösterecek olursak ;$1+1+2+3+5+8+13+21 = 54 = 55-1 = u_{10}-1$

Bu eşitliği aşağıdaki bağıntıları toplayarak elde etmek mümkündür

$$\begin{aligned} u_1 &= u_3 - u_2 \\ u_2 &= u_4 - u_3 \\ u_3 &= u_5 - u_4 \\ &\vdots \\ u_{n-1} &= u_{n+1} - u_n \\ u_n &= u_{n+2} - u_{n+1} \end{aligned}$$

Sol tarafın toplamı ilk $n$ Fibonacci sayısının toplamıdır. Sağ tarafta ise sadece $u_{n+2} - u_2$ kalır. $u_2 = 1$ olduğundan $u_{n+2} - 1$ kalır diyebiliriz. Sonuç olarak $u_1 + u_2 + u_3 + ... + u_n = u_{n+2} - 1$ elde edilir.

Diğer önemli bir Fibonacci özdeşliği $u_n^2 = u_{n+1}u_{n-1} + (-1)^n$ özdeşliğidir. $n = 6$ ve $n = 7$ için örneklendirecek olursak;

$$u_6^2 = 8^2 = 13.5 - 1 = u_7u_5 - 1$$

$$u_7^2 = 13^2 = 21.8 + 1 = u_8u_6 + 1$$

Şimdi özdeşliği oluşturalım.

$$\begin{aligned} u_n^2 - u_{n+1}u_{n-1} &= u_n(u_{n-1} + u_{n-2}) - u_{n+1}u_{n-1} \\ &= (u_n - u_{n+1})u_{n-1} + u_n u_{n-2} \end{aligned}$$

Fibonacci dizisinin tanımından $-u_{n-1} = u_n - u_{n+1}$ olduğundan

$$u_n^2 - u_{n+1}u_{n-1} = (-1)(u_{n-1}^2 - u_n u_{n-2})$$

elde edilir. Şimdi aynı meteodu $u_{n-1}^2 - u_n u_{n-2}$ için uygulayalım.

$$u_n^2 - u_{n+1}u_{n-1} = (-1)^2(u_{n-2}^2 - u_{n-1}u_{n-3})$$

Bu şekilde devam edersek $n - 2$ adım sonra

$$\begin{aligned} u_n^2 - u_{n+1}u_{n-1} &= (-1)^{n-2}(u_2^2 - u_3 u_1) \\ &= (-1)^{n-2}(1^2 - 2.1) = (-1)^{n-1} \end{aligned}$$

$n = 2k$ için özdeşlik aşağıdaki şekle dönüşür.

$$u_{2k}^2 = u_{2k+1}u_{2k-1} - 1$$

Bu iyi bilinen geometrik bir hilenin özünü oluşturur. Şöyle ki; bir kenarı 8 birim olan bir kareyi şekildeki gibi 4 parçaya ayıralım ve alt taraftaki şekilde görüldüğü gibi yeniden düzenleyelim. Figure 0-5'deki kenarları 5 ve 13 birim olan dikdörtgen ile karenin aynı parçalrdan oluştuğu görülmektedir. Fakat karenin alanı 64, dikdörtgenini alanı 65 birim karedir. Yani alanları farklıdır. Bu durum a,b,c ve d noktalarının aynı doğru üzerinde olmamasından kaynaklanır. Fakat bunun gözle görülmesi imkansızdır. Dikdörtgen ile karenin alanları arasındaki fark bizim göremediğimiz, diktörtgeninin tam ortasında yer alan çok çok küçük bir paralelkenardan kaynaklanır.

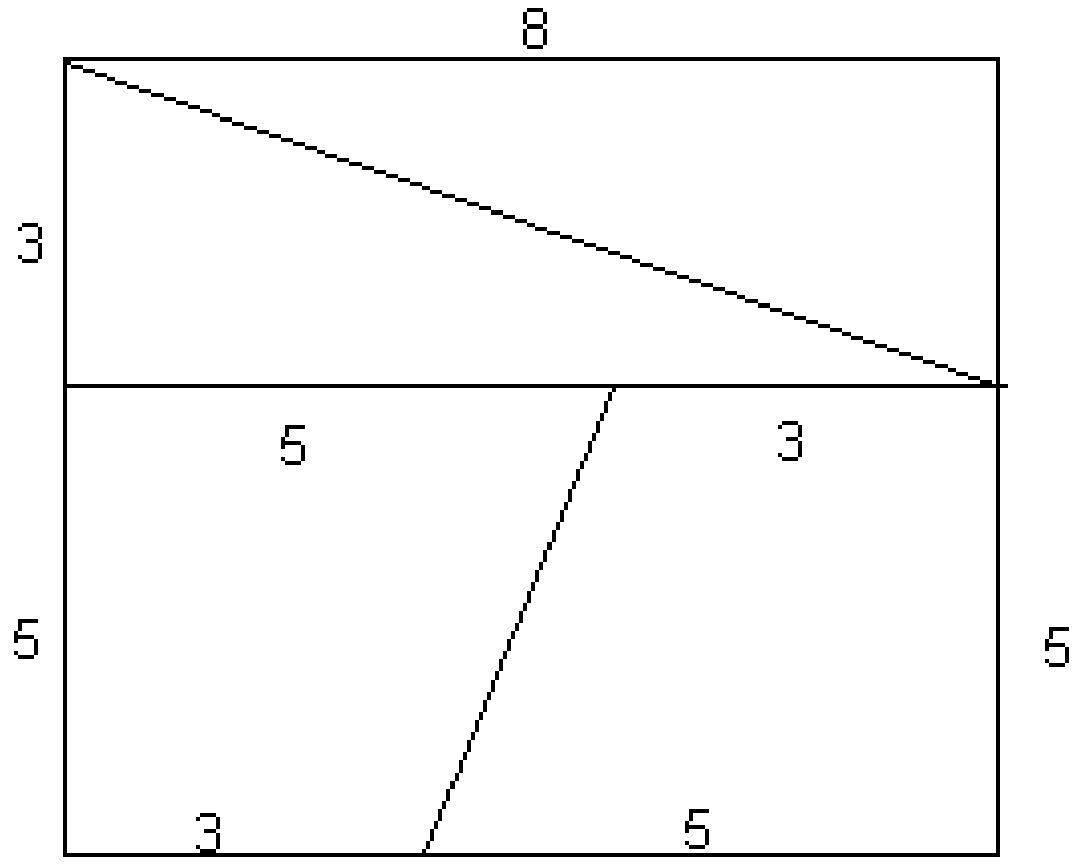

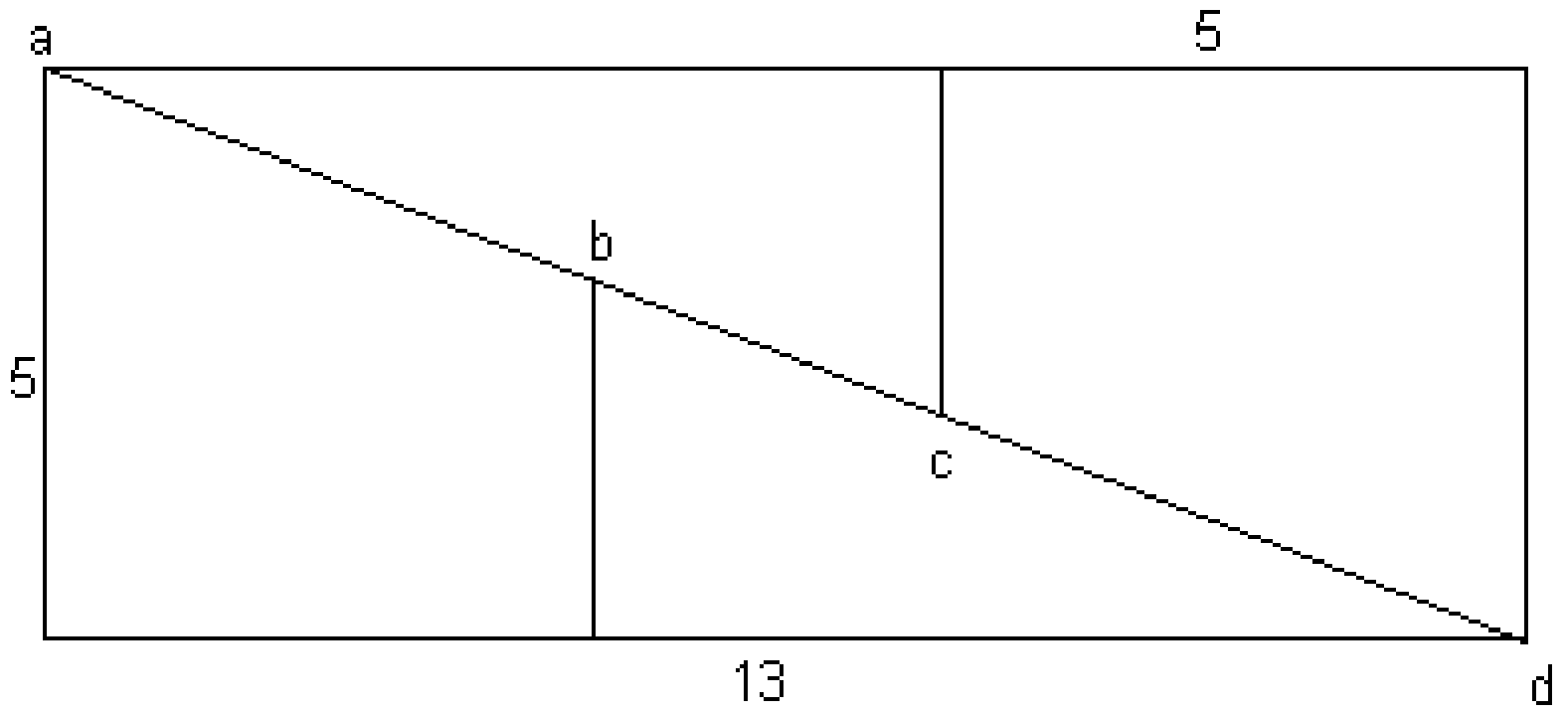

Bu yapı bir kenarı $u_{2k}$ Fibonacci sayısı olan karelere genişletilebilir. Bu tip bir kareyi şekildeki gibi parçalara ayıralım.

Bu parçaları bir dikdörtgen oluşturmak için şekildeki gibi birleştirelim. Bu dikdörtgenin ortasında oluşturduğumuz parçaların dışında, dikdörtgenin içinde çok küçük bir paralelkenar oluşur(şekildeki paralelkenar gerçek halinden oldukça büyüktür).$u_{2k-1}u_{2k+1} - 1 = u_{2k}^2$ ifadesi karenin alanına eşittir

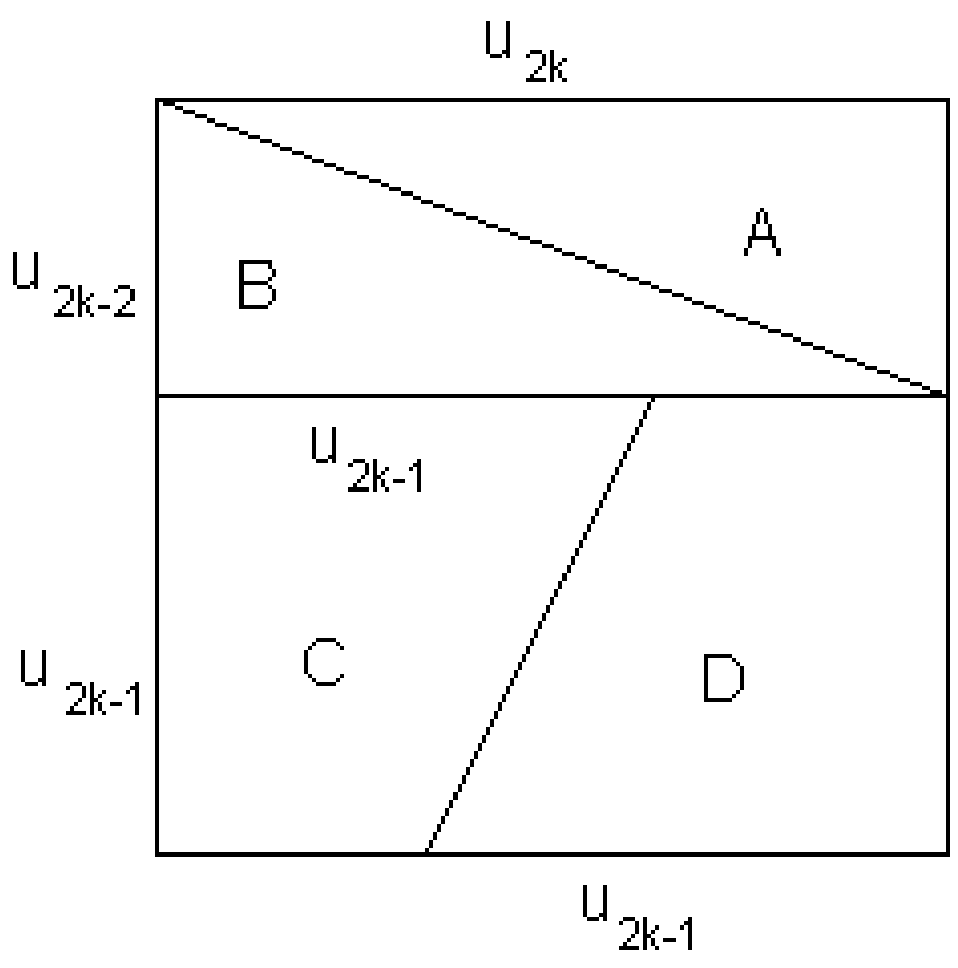

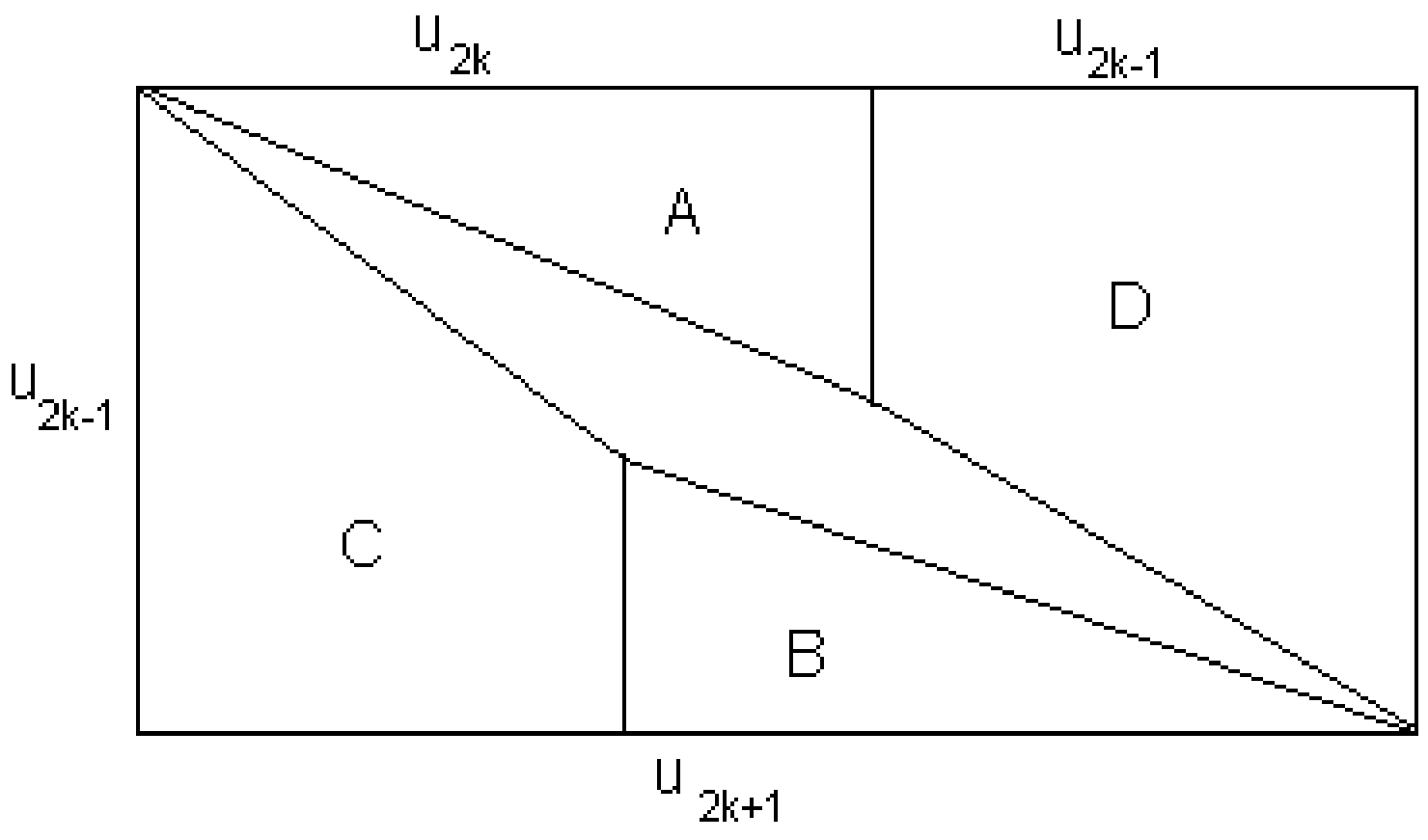

Paralelkenarın Alanı taban kenarı ve o kenara ait yüksekliğin çarpımıdır. A parçasını oluşturan dik üçgende Pisagor teoremini uygularsak paralelkenarın tabanı $\sqrt{u_{2k}^2 + u_{2k-2}^2}$ olarak elde edilir. Paralelkenarın yüksekliğine $h$ diyelim. Paralelkenarın alanı 1'e eşit olduğundan $1 = h \cdot \sqrt{u_{2k}^2 + u_{2k-2}^2}$ ve buradan dikdörtgenin içindeki parelelkenarın yüksekliği $h = \frac{1}{\sqrt{u_{2k}^2+u_{2k-2}^2}}$ olur.

Bilinen sadece üç tane karesel Fibonacci sayısı vardır($u_1 = u_2 = 1, u_{12} = 12^2$). Ayrıca sadece üç tane kübik sayı vardır($u_1 = u_2 = 1, u_6 = 2^3$). Üçgensel Fibonacci sayılarının sayısı da beştir($u_1 = u_2 = 1, u_4 = 3, u_8 = 21, u_{10} = 55$).

Her pozitif tam sayı birbirinden farklı Fibonacci sayılarının toplamı şeklinde yazılabilir. Buna Zeckendorf gösterimi denir. Örnek olarak;

$$
\begin{aligned}
1 &= u_1 & 5 &= u_5 = u_4 + u_3 \\
2 &= u_3 & 6 &= u_5 + u_1 = u_4 + u_3 + u_1 \\
3 &= u_4 & 7 &= u_5 + u_3 = u_4 + u_3 + u_2 + u_1 \\
4 &= u_4 + u_1 & 8 &= u_6 = u_5 + u_4
\end{aligned}
$$

**Teorem 5.** Herhangi bir $N$ pozitif tamsayısı birbirinden farklı olan Fibonacci sayılarının toplamı olarak yazılabilir(herhangi ikisi ardışık olmayacak şekilde).

**İspat.** Tümevarım yöntemini kullanalım. $n > 2$ için $1, 2, 3, ...u_n - 1$ sayılarının herbirinin $\{u_1, u_2, ...u_{n-2}\}$ kümesinden tekrarsız olarak alınan sayıların toplamı olarak yazılabilir olduğunu göstermek yeterlidir. $n = 4$ için $u_4 = 3 = u_1 + u_3$. Şimdi $n = k$ için sağlandığını varsayalım. $N$ sayısını $u_k - 1 < N < u_{k+1}$ olacak şekilde seçelim. Buradan $N - u_{k-1} < u_{k+1} - u_{k-1} = u_k$ elde edilir. Buradan $N - u_{k-1}$ sayısının $\{1, 2, ...u_{k-2}\}$ kümesinden alınan birbirinden farklı sayıların toplamı olacak şekilde yazılabileceği sonucu çıkar. O halde $N$ ve dolayısıyla $1, 2, 3, ...u_{k+1} - 1$ sayılarından herbiri $\{u_1, u_2, ...u_{k-2}, u_{k-1}\}$ kümesinden tekrarsız olarak alınan sayıların toplamı olarak yazılabilir. Bu da ispatı tamamlar. Bu teoremde ardışık olmayan fibonacci sayılarını göz önüne aldık. Çünkü ardışık iki Fibonacci sayısının toplamı bir sonraki Fibonacci sayısını verir. ■

**Örnek.** $N = 50$ olsun. Burada $u_9 < 50 < u_{10}$ ve Zeckendorf gösterimi şu şekilde olur:

$$50 = u_4 + u_7 + u_9$$

$$u_n = \frac{1}{\sqrt{5}}\left[\left(\frac{1+\sqrt{5}}{2}\right)^n - \left(\frac{1-\sqrt{5}}{2}\right)^n\right]$$

genel Fibonacci terimini göz önüne alalım. Bu denklemin köklerini

$$\alpha = \frac{1+\sqrt{5}}{2} \qquad ve \qquad \beta = \frac{1-\sqrt{5}}{2}$$

ile işaretleyelim. Bunlar $x^2 - x - 1 = 0$ denkleminin kökleri olduğundan şu şekilde yazılabilirler:

$$\alpha^2 = \alpha + 1 \qquad ve \qquad \beta^2 = \beta + 1$$

Şimdi her iki tarafı her iki denklemde sırasıyla $\alpha^n$ ve $\beta^n$ ile çarpalım.

$$\alpha^{n+2} = \alpha^{n+1} + \alpha^n \quad ve \quad \beta^{n+2} = \beta^{n+1} + \beta^n$$

Şimdi ise birinci denlemden ikinci denklemi çıkaralım ve her tarafı $\alpha - \beta$ ile bölelim.

$$\frac{\alpha^{n+2} - \beta^{n+2}}{\alpha - \beta} = \frac{\alpha^{n+1} - \beta^{n+1}}{\alpha - \beta} + \frac{\alpha^n - \beta^n}{\alpha - \beta}$$

$H_n = (\alpha^n - \beta^n)/(\alpha - \beta)$ ile işaretleyelim. Böylece yukarıdaki denklem

$$H_{n+2} = H_{n+1} + H_n \qquad n \geq 1$$

şeklini alır. $\alpha$ ve $\beta$ hakkında aşağıdakiler doğrudur.

$$\alpha + \beta = 1 \qquad \alpha - \beta = \sqrt{5} \qquad \alpha\beta = -1$$

burada,

$$H_1 = \frac{\alpha - \beta}{\alpha - \beta} = 1 \qquad H_2 = \frac{\alpha^2 - \beta^2}{\alpha - \beta} = \alpha + \beta = 1$$

Bütün bunlardan $H_1, H_2, H_3...$ bir Fibonacci dizisi olur.

$$u_n = \frac{\alpha^n - \beta^n}{\alpha - \beta} \qquad n \geq 1$$

şeklinde yazabiliriz. Bu gösterim Binet formülü olarak adlandırılır. ve Fibonacci sayılarıyla bağlantılı birçok özelliği elde etmekte kullanılır. Şimdi bunu kullanarak aşağıdaki eşitliği gösterelim.

$$u_{n+2}^2 - u_n^2 = u_{2n+2}$$

Göstermeden önce $\alpha\beta = -1$ olduğundan $k \geq 1$ için $(\alpha\beta)^{2k} = 1$ olduğunu hatır-latalım.

$$\begin{aligned}
u_{n+2}^2 - u_n^2 &= \left(\frac{\alpha^{n+2}-\beta^{n+2}}{\alpha-\beta}\right)^2 - \left(\frac{\alpha^n-\beta^n}{\alpha-\beta}\right)^2 \\
&= \frac{\alpha^{2(n+2)}-2+\beta^{2(n+2)}}{(\alpha-\beta)^2} - \frac{\alpha^{2n}-2+\beta^{2n}}{(\alpha-\beta)^2} \\
&= \frac{\alpha^{2(n+2)}+\beta^{2(n+2)}-\alpha^{2n}-\beta^{2n}}{(\alpha-\beta)^2}
\end{aligned}$$

Paydaki ifade şu şekilde yazılabilir:

$$\alpha^{2(n+2)} - (\alpha\beta)^2\alpha^{2n} - (\alpha\beta)^2\beta^{2n} + \beta^{2(n+2)} = (\alpha^2-\beta^2)(\alpha^{2n+2}-\beta^{2n+2})$$

Bunu yerine koyarsak

$$\begin{aligned}
u_{n+2}^2 - u_n^2 &= \frac{(\alpha^2-\beta^2)(\alpha^{2n+2}-\beta^{2n+2})}{(\alpha-\beta)^2} \\
&= (\alpha+\beta)\left(\frac{\alpha^{n+2}-\beta^{n+2}}{\alpha-\beta}\right) \\
&= 1.u_{2n+2} = u_{2n+2}
\end{aligned}$$

Şimdi daha önce de ispatladığımız $u_{2n}^2 = u_{2n+1}u_{2n-1} - 1$ bağıntısını Binet formülünü kullanarak ispatlayalım.

$$\begin{aligned}
u_{2n+1}u_{2n-1} - 1 &= \left(\frac{\alpha^{2n+1}-\beta^{2n+1}}{\sqrt{5}}\right)\left(\frac{\alpha^{2n-1}-\beta^{2n-1}}{\sqrt{5}}\right) - 1 \\
&= \frac{1}{5}\left(\alpha^{4n}+\beta^{4n}-(\alpha\beta)^{2n-1}\alpha^2-(\alpha\beta)^{2n-1}\beta^2-5\right) \\
&= \frac{1}{5}\left(\alpha^{4n}+\beta^{4n}+\left(\alpha^2+\beta^2\right)-5\right)
\end{aligned}$$

$\alpha^2+\beta^2=3$ olduğundan son eşitlik,

$$\begin{aligned}
\frac{1}{5}\left(\alpha^{4n}+\beta^{4n}-2\right) &= \frac{1}{5}\left(\alpha^{4n}+\beta^{4n}-2(\alpha\beta)^{2n}\right) \\
&= \left(\frac{\alpha^{2n}-\beta^{2n}}{\sqrt{5}}\right)^2 = u_{2n}^2
\end{aligned}$$

şekline dönüşür ki bu da bizim elde etmek istediğimiz eşitliktir.

Binet formülü aynı zamanda Fibonacci sayılarının değerini elde etmekte de kullanılır. $0<|\beta|<1$ eşitsizliği $n\geq 1$ için $|\beta^n| = |\beta|^n$ olmasını gerektirir. Bu

nedenle

$$\left|u_n - \frac{\alpha^n}{\sqrt{5}}\right| = \left|\frac{\alpha^n - \beta^n}{\sqrt{5}} - \frac{\alpha^n}{\sqrt{5}}\right| = \frac{|\beta^n|}{\sqrt{5}} < \frac{1}{\sqrt{5}} < \frac{1}{2}$$

olur. Burada $u_n$ $\frac{\alpha^n}{\sqrt{5}}$ sayısına en yakın tamsayıdır. Örnek olarak, $\frac{\alpha^{14}}{\sqrt{5}} \approx 377.0005$ ve $u_{14} = 377$. Benzer şekilde, $\frac{\alpha^{15}}{\sqrt{5}} \approx 609,9996$ olduğundan $u_{15} = 610$ olur. Yani $u_n$ sayısı aslında $\frac{\alpha^n}{\sqrt{5}} + \frac{1}{2}$ sayısını aşmayan en büyük tamsayıdır. Tam değer fonksiyonu ile gösterecek olursak,

$$u_n = \left[\frac{\alpha^n}{\sqrt{5}} + \frac{1}{2}\right] \qquad n \geq 1$$

**Teorem 6.** $p > 5$ asal sayısı için ya $p \mid u_{p-1}$ ya da $p \mid u_{p+1}$.

**İspat.** Binet formülünden $u_p = (\alpha^p - \beta^p)/\sqrt{5}$ yazılabilir. $\alpha$ ve $\beta$ nın $p$. mertebeden kuvvetlerini binom teoreminden açarsak şunu elde ederiz:

$$\begin{aligned} u_p &= \frac{1}{2^p\sqrt{5}}\left[1 + \binom{p}{1}\sqrt{5} + \binom{p}{2}5 + \binom{p}{3}5\sqrt{5} + ... + \binom{p}{p}5^{(p-1)/2}\sqrt{5}\right] \\ &\quad -\frac{1}{2^p\sqrt{5}}\left[1 - \binom{p}{1}\sqrt{5} + \binom{p}{2}5 - \binom{p}{3}5\sqrt{5} + ... - \binom{p}{p}5^{(p-1)/2}\sqrt{5}\right] \\ &= \frac{1}{2^{p-1}}\left[\binom{p}{1} + \binom{p}{3}5 + \binom{p}{5}5^2 + ... + \binom{p}{p}5^{(p-1)/2}\right] \end{aligned}$$

$1 \leq k \leq p-1$ için $\binom{p}{k} \equiv 0(\text{mod } p)$ ve $2^{p-1} \equiv 1(\text{mod } p)$ olduğunu hatırlayalım. Bunları göz önüne alırsak $u_p$ yi daha basit bir şekilde şöyle ifade edebiliriz:

$$u_p \equiv 2^{p-1}u_p \equiv \binom{p}{p}5^{(p-1)/2} = 5^{(p-1)/2}(\text{mod}\, p)$$

O halde $u_p \equiv (5/p) \equiv \pm1(\text{mod } p)$ ve buradan $u_p^2 \equiv 1(\text{mod } p)$. Şimdi $u_p^2 = u_{p-1}u_{p+1} + (-1)^{p-1}$ bağıntısını modül $p$ ye göre yeniden düzenleyelim. $u_{p-1}u_{p+1} \equiv 0(\text{mod } p)$ elde edilir. Yani $p$ ya $u_{p-1}$ i ya da $u_{p+1}$ i böler. $p$ her ikisini de bölemez. Çünkü $(p-1, p+1) = 2$ ve teorem 3 ü göz önüne alırsak $(u_{p-1}, u_{p+1}) = u_2 = 1$.Yani $u_{p-1}$ ve $u_{p+1}$ aralarında asaldır. Böylece ispat tamamlandı. ■

**Teorem 7.** $p \geq 7$ asal sayısı için $p \equiv 2(\text{mod } 5)$ veya $p \equiv 4(\text{mod } 5)$ olsun. Eğer $2p-1$ sayısı asalsa o halde $2p-1 \mid u_p$.

**İspat.** Varsayalım ki $p$ bazı $k$ sayısı için $5k+2$ formunda olsun. $u_p = (\alpha^p - \beta^p)/\sqrt{5}$ formülünün karesini alalım ve $\alpha^{2p}$ ve $\beta^{2p}$ yi binom teoremini kul-lanarak açalım.

$$5u_p^2 = \frac{1}{2^{2p-1}}\left[1 + \binom{2p}{2}5 + \binom{2p}{4}5^2 + ... + \binom{2p}{2p}5^p\right] + 2$$

elde edilir. $2 \leq k < 2p-1$ için $\binom{2p}{k} \equiv 0(\text{mod } 2p-1)$ olur. Çünkü $2p-1$ asal sayısı $2^{2p-1} \equiv 2(\text{mod } 2p-1)$ özelliğini sağlar. Böylece $u_p^2$ için şunu elde ederiz:

$$2\left(5u_p\right)^2 \equiv (1 + 5^p) + 4(\text{mod}\, 2p-1)$$

veya daha basit olarak $2u_p^2 \equiv 1 + 5^{p-1}$ (mod $2p-1$) buradan $5^{p-1} = 5^{(2p-2)/2} \equiv (5/2p-1)\,(\text{mod } 2p-1)$. Buradan şunu görmek kolaydır:

$$(5/2p-1) = (2p-1/5) = (10k+3/5) = (3/5) = -1$$

Böylece $2u_p^2 \equiv 1 + (-1) \equiv 0(\text{mod } 2p-1)$ ifadesine ulaşırız. Buradan da $2p-1$ böler $u_p$ elde ederiz. $p \equiv 4(\text{mod } 5)$ için ispat benzer yolla yapılır. ■

# 3 GEOMETRİDEKİ VE DOĞADAKİ UYGULAMALARI

Bilindiği gibi Altın Kesit bir $AB$ doğru parçasının, $\frac{AB}{AC} = \frac{AC}{CB}$ orantısına uygun olarak $C$ noktası tarafından bölünmesidir.($AB$ gibi bir doğru parçasının uzunluğunu karışıklığa yol açmadığı takdirde yine $AB$ ile göstereceğiz.)

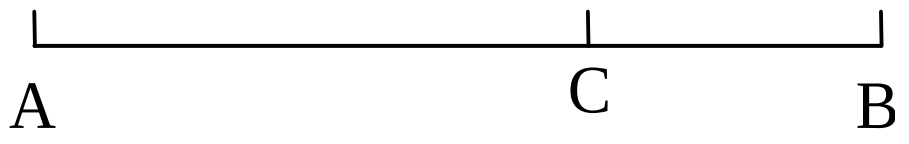

$$AB = AC + CB$$

olduğu göz önüne alınırsa,

$$\frac{AC + CB}{AC} = \frac{AC}{CB}$$

$$1 + \frac{1}{\frac{AC}{CB}} = \frac{AC}{CB}$$

$$\frac{AC}{CB} = \frac{1 + \sqrt{5}}{2}$$

elde edilir.

$$\frac{AB}{AC} = \frac{AC}{CB} = \frac{1 + \sqrt{5}}{2}$$

oranına (sayısına) Altın Oran denir. Altın Oran'ı anlatmanın en iyi yollarından biri, işe bir kare ile başlamaktır. Bir kareyi tam ortasından iki eşit diktörgen oluşturacak şekilde ikiye bölelim. Dikdörtgenlerin ortak kenarının, karenin tabanını kestiği noktaya pergelimizi koyalım. Pergelimizi öyle açalım ki, çizeceğimiz daire, karenin karşı köşesine değsin, yani yarı çapı, bir dikdörtgenin köşegeni olsun.

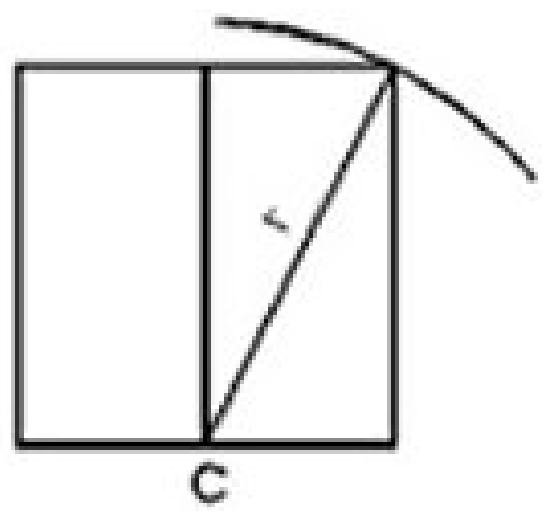

Sonra, karenin tabanını, çizdiğimiz daireyle kesişene kadar uzatalım.

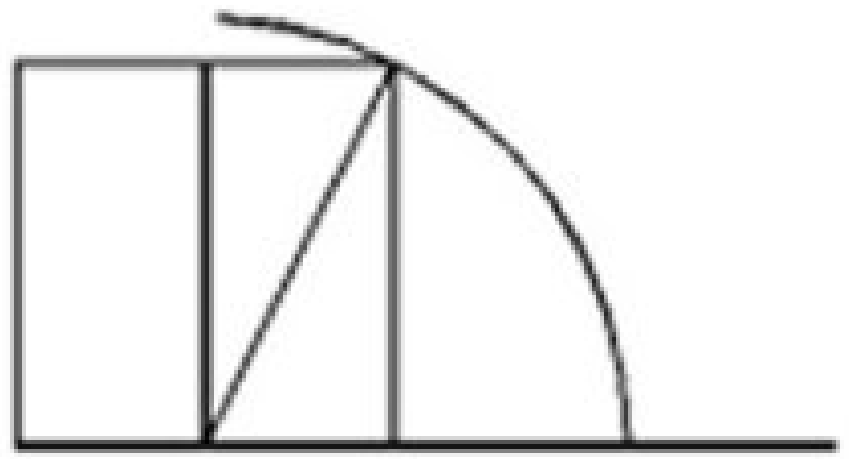

Yeni çıkan şekli bir dikdörtgene tamamladığımızda, karenin yanında yeni bir dikdörtgen elde etmiş olacağız.

İşte bu yeni dikdörtgenin taban uzunluğunun ($B$) karenin taban uzunluğuna ($A$) oranı Altın Oran'dır. Karenin taban uzunluğunun ($A$) büyük dikdörtgenin taban uzunluğuna ($C$) oranı da Altın Oran'dır. $A/B = 1.6180339 =$ Altın Oran $C/A = 1.6180339 =$ Altın Oran

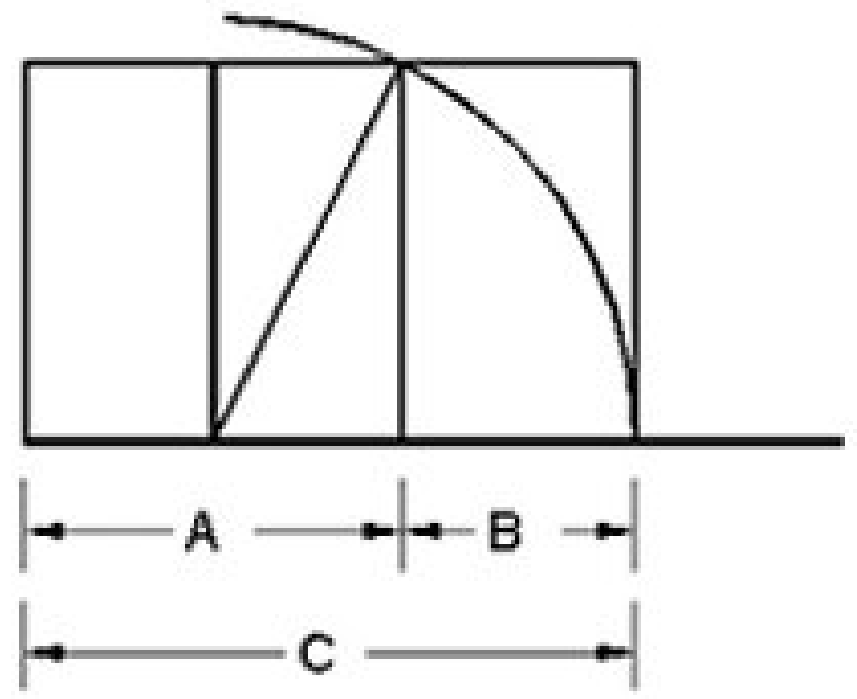

Elde ettiğimiz bu dikdörtgen ise, bir Altın Dikdörtgen'dir. Çünkü kısa kenarının, uzun kenarına oranı 1.618 dir, yani Altın Oran'dır. Artık bu dikdörtgenden her bir kare çıkardığımızda elimizde kalan, bir Altın Dikdörtgen olacaktır. Bunları şu şekilde de ifade edebiliriz. Yukarıdaki düşünceden istifade ederek bir $ABCD$ Altın Dikdörtgeninde sırasıyla noktalarının çizilmesiyle gitgide küçülen (alanları $\Phi$ oranında küçülen ) Altın Dikdörtgenler elde edilir.

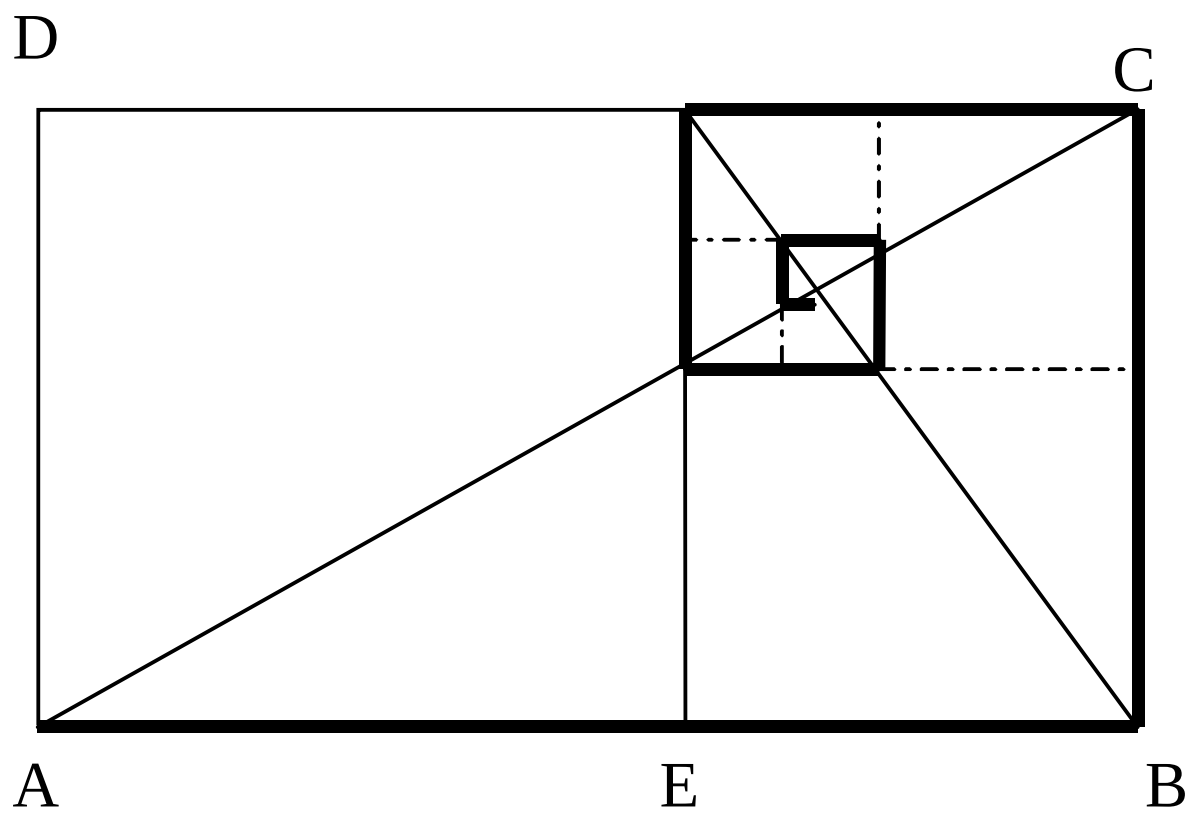

Figure 0-1:

$AB$ kenarını kısa kenar kabul edip gittikçe büyüyen Altın Dikdörtgenler de çizilebilir. $AB, BC, CE_1, E_1E_2, E_2E_3...$ doğru parçalarının oluşturduğu şekle dik çizgili sarmal (spiral) denir.$ABCD$ Altın Dikdörtgeninde,

$$\frac{AB}{BC} = \frac{BC}{CE_1} = \frac{CE_1}{E_1E_2} = \frac{E_1E_2}{E_2E_3} = ... = \Phi$$

İçinden defalarca kareler çıkardığımız bu Altın Dikdörtgen'in karelerinin kenar uzunluklarını yarıçap alan bir çember parçasını her karenin içine çizersek, bir Altın Spiral elde ederiz. Altın Spiral, birçok canlı ve cansız varlığın biçimini ve yapı taşını oluşturur.Buna örnek olarak Ayçiçeği bitkisini gösterebiliriz. Ayçiçeğinin çekirdekleri altın oranı takip eden bir spiral oluşturacak şekilde dizilirler.

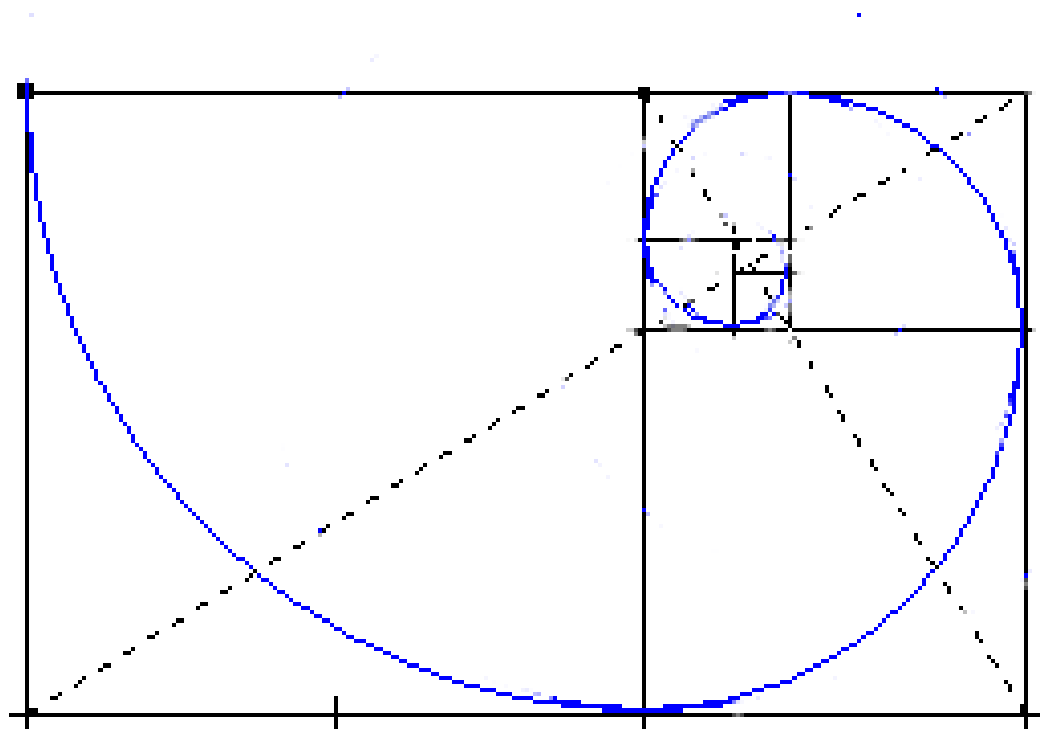

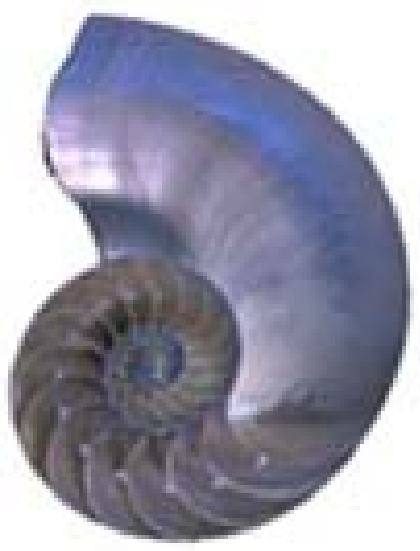

Bu karelerin kenar uzunlukları sırasıyla Fibonacci sayılarını verir.

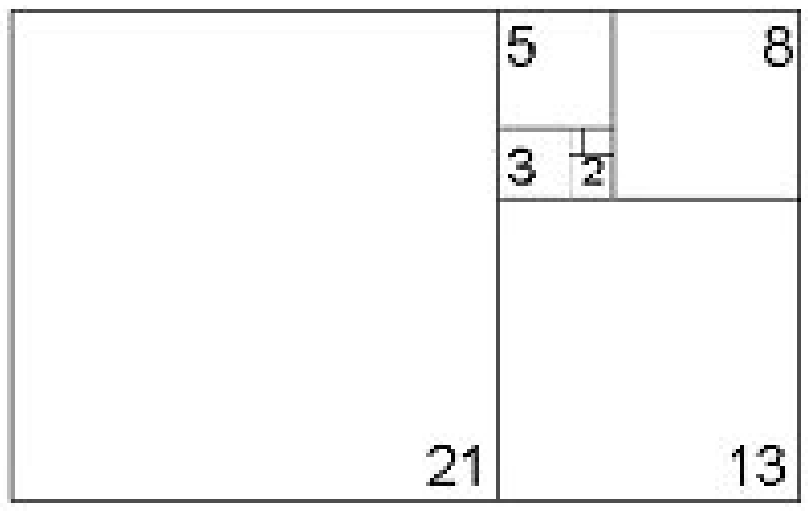

**Beş Kenarlı Simetri**

Φ'yi göstermenin bir yolu da, basit bir beşgen kullanmaktır. Yani, birbiriyle beş eşit açı oluşturarak birleşen beş kenar. Basitçe $\Phi$, herhangi bir köşegenin herhangi bir kenara oranıdır. $AC/AB = 1,618 = \Phi$

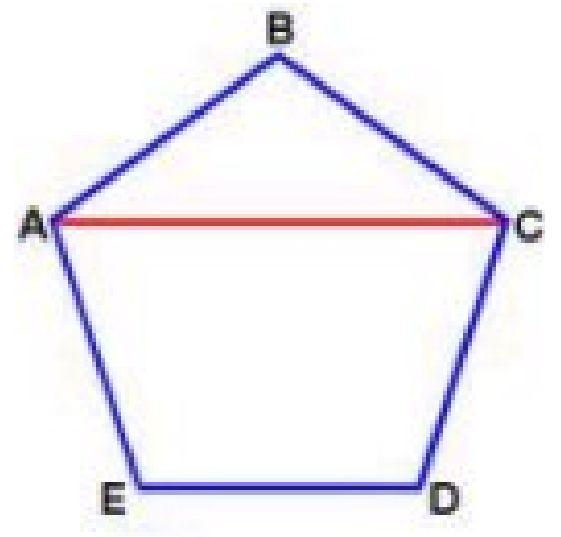

Beşgenin içine ikinci bir köşegen ($[BD]$) çizelim. $AC$ ve $BD$ birbirlerini $O$ noktasında keseceklerdir.

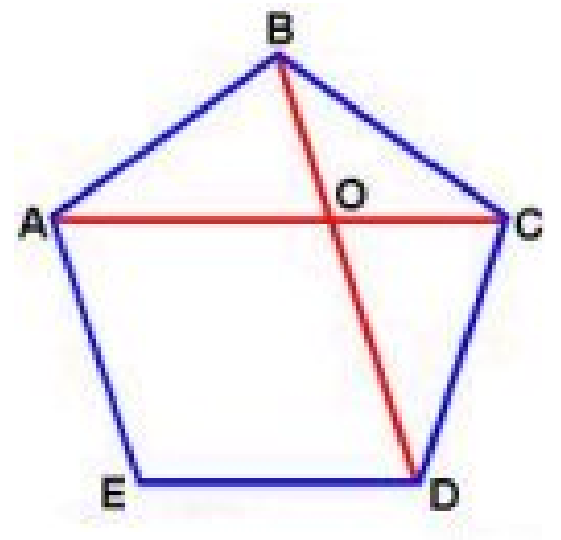

Böylece her iki çizgi de, bir noktadan ikiye bölünmüş olacaktır ve her parça diğeriyle $\Phi$ oranı ilişkisi içindedir. Yani $AO/OC = \Phi$, $AC/AO = \Phi$, $DO/OB = \Phi$, $BD/DO = \Phi$. Bir diğeri ile bölünen her köşegende, aynı oran tekrarlanacaktır. Bütün köşegenleri çizdiğimiz zaman ise, beş köşeli bir yıldız elde ederiz.

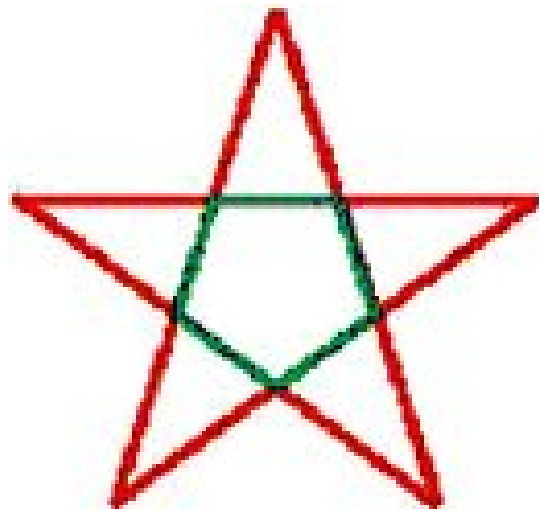

Bu yıldızın içinde, ters duran diğer bir beşgen meydana gelir. Her köşegen, başka iki köşegen tarafından kesilmiştir ve her bölüm, daha büyük bölümlerle ve bütünle, Φ oranını korur. Böylece, içteki ters beşgen, dıştaki beşgenle de Φ oranındadır.

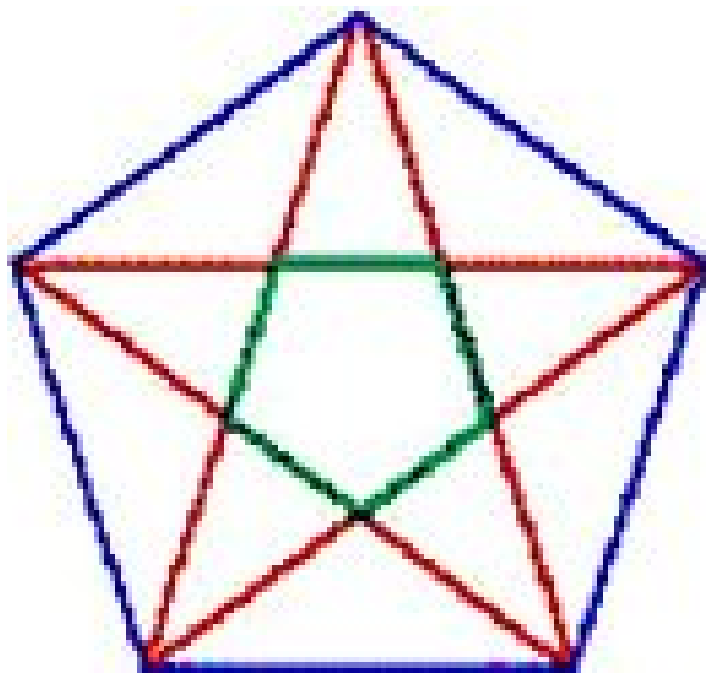

Bir beşgenin içindeki beş köşeli yıldız, Pentagram diye adlandırılır ve Pythagoras'ın kurduğu antik Yunan Matematik Okulu'nun sembolüdür. Eski gizemciler Φ'yi bilirlerdi ve Altın Oran'ın fiziksel ve biyolojik dünyamızın kurulmasındaki önemli yerini anlamışlardı. Bir beşgenin köşegenlerini birleştirdiğimizde, iki değişik Altın Üçgen elde ederiz. Figure 0-20'deki koyu renkli üçgenlerin tabanları kenarları ile Altın Oran ilişkisi içerisindedir.

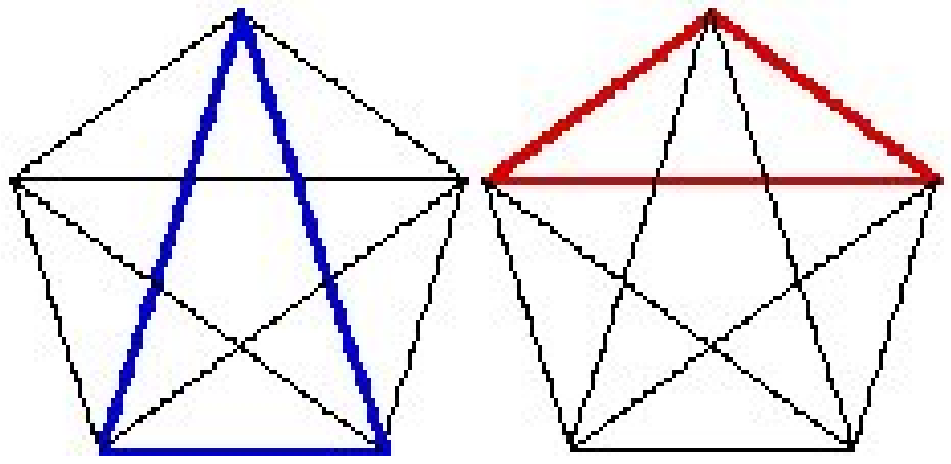

Φ, kendini tekrarlayan bir özelliğe de sahiptir. Altın Orana sahip her şekil, Altın Oranı kendi içinde sonsuz sayıda tekrarlayabilir. Aşağıdaki şekilde, her beşgenin içinde meydana gelen pentagramı ve her pentagramın oluşturduğu beşgeni ve bunun makro kozmik ve mikro kozmik sonsuza kadar Altın Oranı tekrarlayarak devam ettiğini görebiliriz.

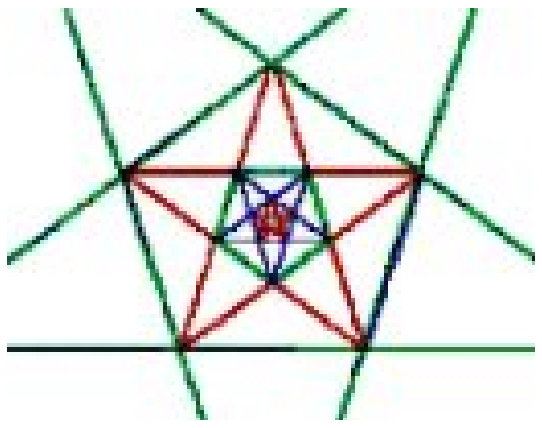

Beşgen, Altın Oranı açıklamak için oldukça basit ve iyi bir yöntem olmakla birlikte, bu oranın belirtilmesi gereken çok daha karmaşık ve anlaşılması zor bir takım özellikleri de vardır. Altın Oran daha iyi anlaşıldıkça, biyolojik ve kozmolojik birçok büyük uygulama örnekleri daha iyi görülebilecektir.

Bir $AB$ doğru parçasının Altın Oranını bulmak için $AB$ ye dik olan ve uzunluğu $BD = AB/2$ olan $BD$ doğru parçası çizilir. Kenar uzunlukları $1 : 2 : \sqrt{5}$ oranında olan $DBA$ dik üçgeninde $D$ merkezli $DB$ yarıçaplı çemberin $AD$ hipotenüsü ile arakesiti ($E$ noktası) bulunur. $A$ merkezli $AE$ yarıçaplı çemberin $AB$ doğru parçası ile arakesiti ($C$ noktası) Altın Oranı vermektedir(Figure 0-22).

Eşit açılı logaritmik spiralin kutupsal koordinatlarda denklemi, $k$ ve $c$ pozitif reel sabitler olmak üzere,

$$p = ke^{c.\theta} \quad ,\theta \in [0, \infty)$$

biçimindedir. $p$ ışın vektörü ile ışın vektörünün spirali kestiği noktadaki teğet arasındaki açının tanjantı,

$$\tan \varphi = \frac{p}{\frac{dp}{d\theta}} = \frac{1}{c}$$

olur.

Figure 0-17'deki $A, B, C, E_1, E_2, E_3....$ noktalarından geçen eşit açılı logaritmik spiral için,

$c = \frac{2}{\pi} \ell n\Phi$ ve $\varphi \approx 73^0$ olur.

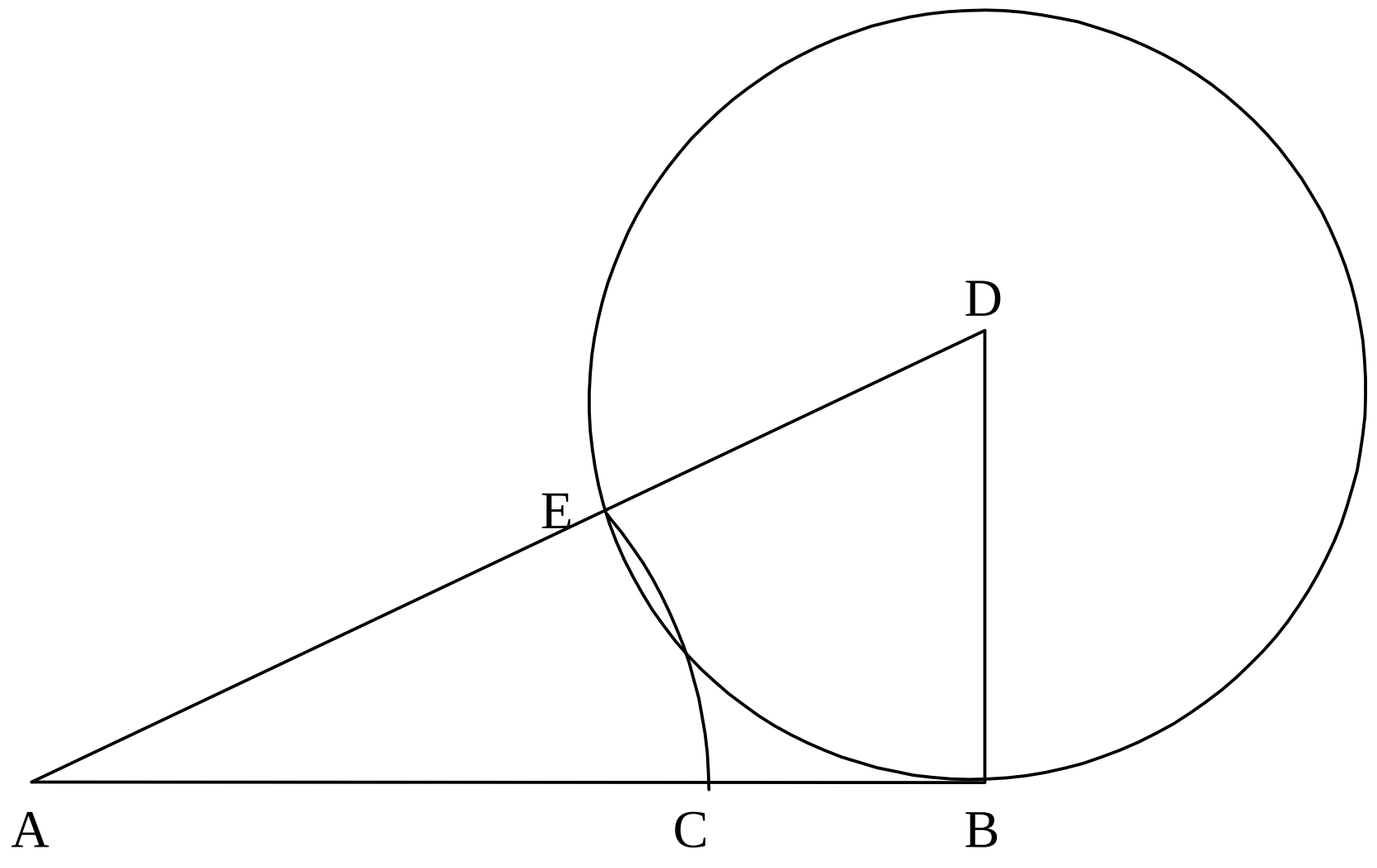

Birçok canlının büyüme sırasında şekilsel olarak logaritmik spirali izlediği gözlenmiştir. Bu konuda D'arcy Thompson: "Bir deniz kabuğunun büyüme sürecinde, aynı ve değişmez orantılara bağlı olarak genişlemesi ve uzamasından daha sade bir sistem düşünemeyiz; nitekim Doğa da son derece basit olan bu yasayı izler. Kabuk giderek büyür, fakat şeklini değiştirmez. İşte, sabit kalan bu büyüme göreceliğinin ya da form özdeşliğinin varlığı, eşitaçılı spiralin özünü ve belki de tanımının esasını oluşturur" diyor.

Altın dikdörtgenler estetik açıdan göze en hoş görünen dikdörtgenlerdir. Kendilerine çeşitli dikdörtgen örnekleri gösterilip, en güzelini ve en çirkinini seçmeleri istenilerek yapılan araştırmalarda, örneğin Fechner (1876) ve Lalo (1908) 'nun anket sonuçları aşağıdaki gibi olmuştur. (Bilim ve Teknik, Cilt 25, Sayı 297,Sayfa 7)

Kenar uzunlukları oranı Altın Orana yakın olan diktörtgenlerin beğenilme yüzdesi büyük olarak gözlenmiştir. Karenin de beğenilmiş olmasına dikkat ediniz. Altın Dikdörtgene çirkin diyen bir tek kişi bile çıkmamıştır.

Mimaride, inşa edilecek yapının cephe görünüşünün daima bir Altın Dikdörtgen içine yerleştirilebilmesi dikkat edilecek ilk husus olmaktadır. Mimaride olduğu gibi, resimde de temel ögelerden biri Altın Oran'dır.

**Altın Üçgen:** Tepe açısı olan ikizkenar üçgenlere Altın Üçgenler denir.

| ORAN Genişlik/Uzunluk | EN GÜZEL DİKDÖRTGEN | | EN ÇİRKİN DİKDÖRTGEN | |
|---|---|---|---|---|
| | Fechner (%) | Lalo (%) | Fechner (%) | Lalo (%) |
| 1.00 | 3.0 | 11.7 | 27.8 | 22.5 |
| 0.83 | 0.2 | 1.0 | 19.7 | 16.6 |
| 0.80 | 2.0 | 1.3 | 9.4 | 9.1 |
| 0.75 | 2.5 | 9.5 | 2.5 | 9.1 |
| 0.69 | 7.7 | 5.6 | 1.2 | 2.5 |
| 0.67 | 20.6 | 11.0 | 0.4 | 0.6 |
| 0.62 | 35.0 | 30.3 | 0.0 | 0.0 |
| 0.57 | 20.0 | 6.3 | 0.8 | 0.6 |
| 0.50 | 7.5 | 8.0 | 2.5 | 12.5 |
| 0.40 | 1.5 | 15.3 | 35.7 | 26.6 |

Altın Üçgen içine çizilen gitgide küçülen Altın Üçgenlerin köşelerinden de eşit açılı logaritmik spiral geçmektedir.

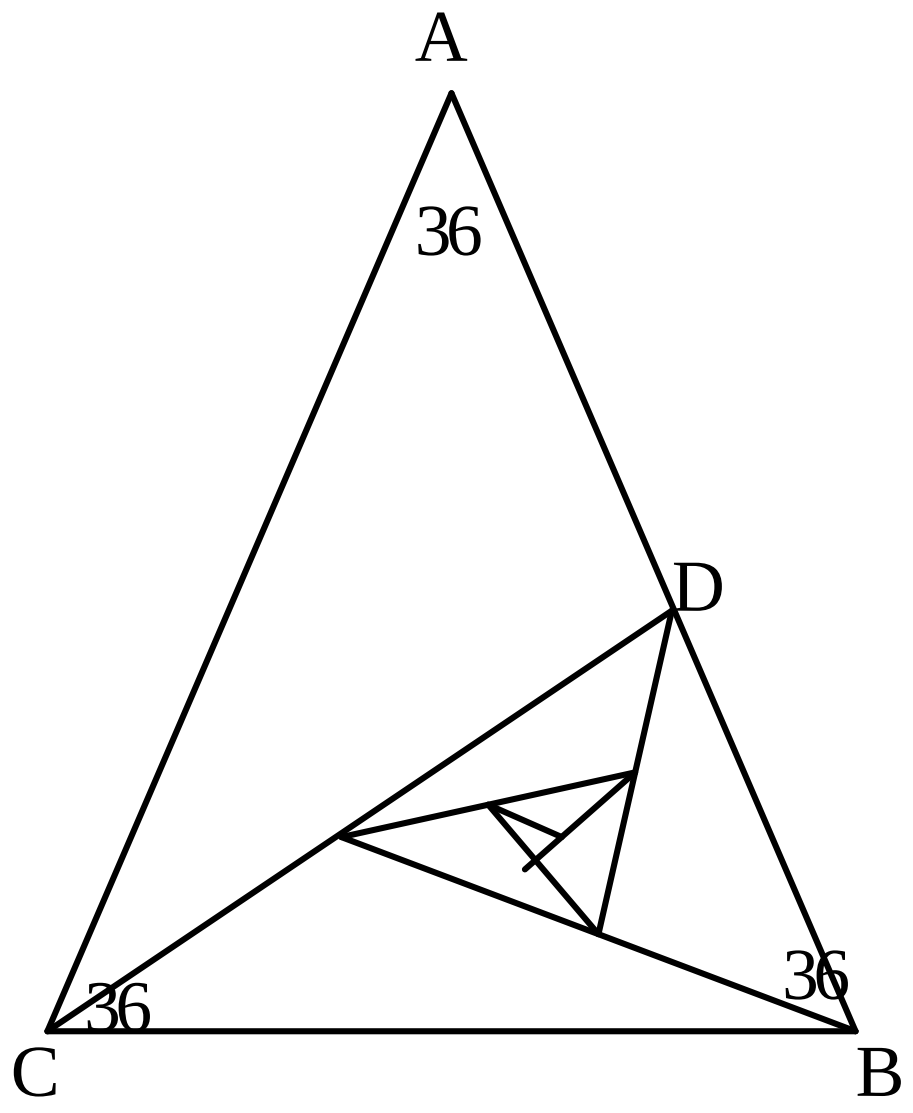

Düzgün bir ongen esasında 10 tane Altın Üçgen diliminden oluşmaktadır. Altın Üçgen düzgün beşgenlerde de karşımıza çıkmaktadır.

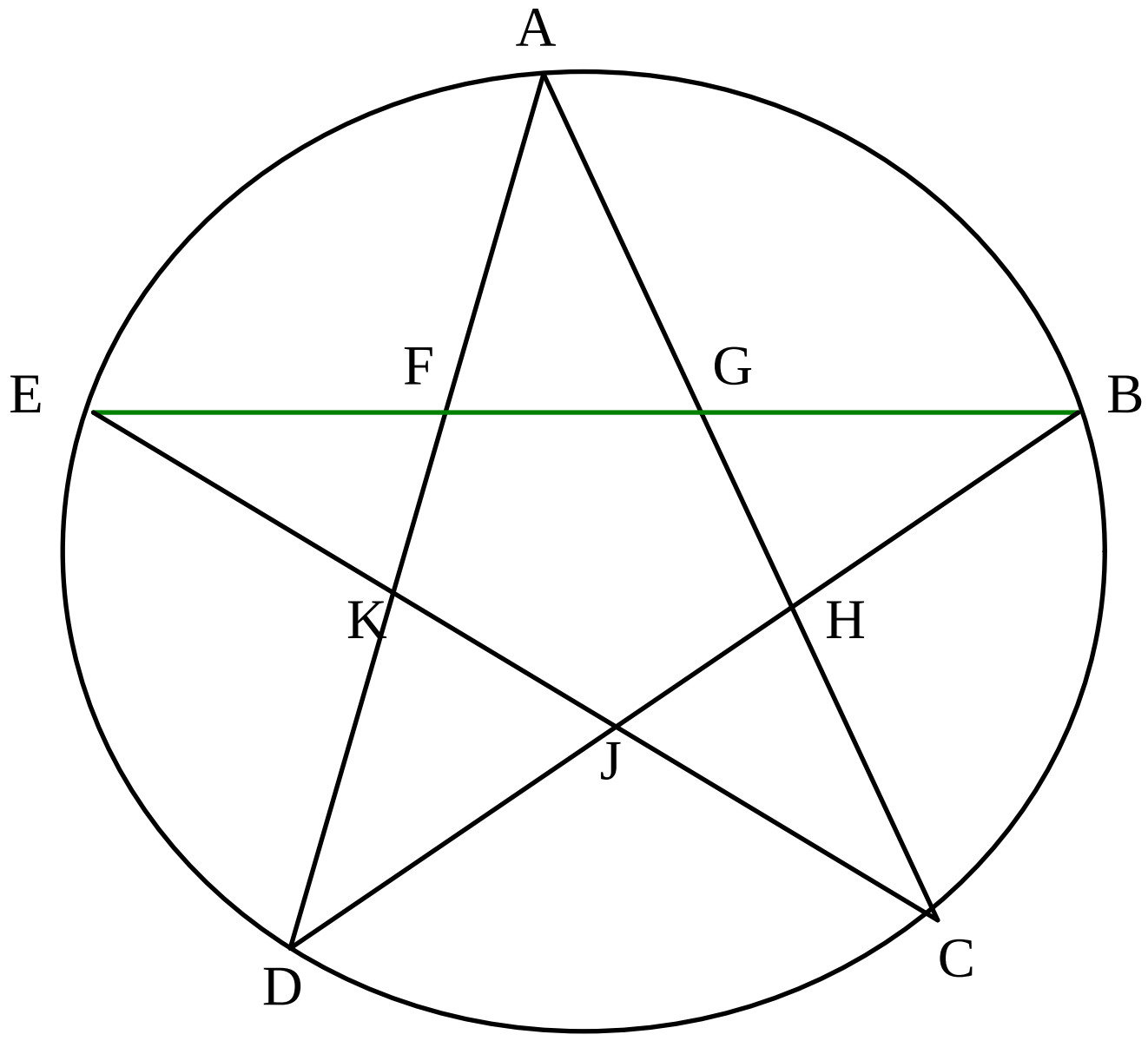

Figure 0-25'deki $FGHJK$ düzgün beşgeninin kenar uzunluğu 1 birim olursa aşağıdaki özellikler yazılabilir.

$$AF = AG = BG = ... = \Phi$$

$$GJ = HK = JK = KG == FH = \Phi$$

$$KP/HP = GP/JP = ... = \Phi$$

$$OB/ON = OC/OS = ... = 2\Phi$$

$$ON/OK = \Phi/2$$

$$OC/OF = \Phi^2$$

$$AB = BC = ... = \Phi^2$$

$$BD/OB = \sqrt{1 + \Phi^2}$$

Yıldızın kolları, çekirdeğini oluşturan $FGHJK$ beşgeninin kenarlarından yukarıya doğru kıvrılıp birleştirilirse oluşan piramidin yüksekliğinin, tabanı çevreleyen çemberin yarıçapına oranı $\Phi$ dir.

Norman Gowar'ın ne dediğine bir bakalım. "Belirli bir sayının, birbirinden bağımsız olarak hem matematik hem de estetik bilimlerini ilgilendiren bir çekiciliği olması, insanı çileden çıkaracak derecede ilginç bir husustur. Üstelik, Altın Oran, insanların tasarımından kaynaklanmaksızın Doğa'da da ortaya çıkmaktadır."

Fibonacci sayılarına binom açılımında terimlerin katsayılarını veren pascal üçgeninde de rastlanılır(Figure 0-26).

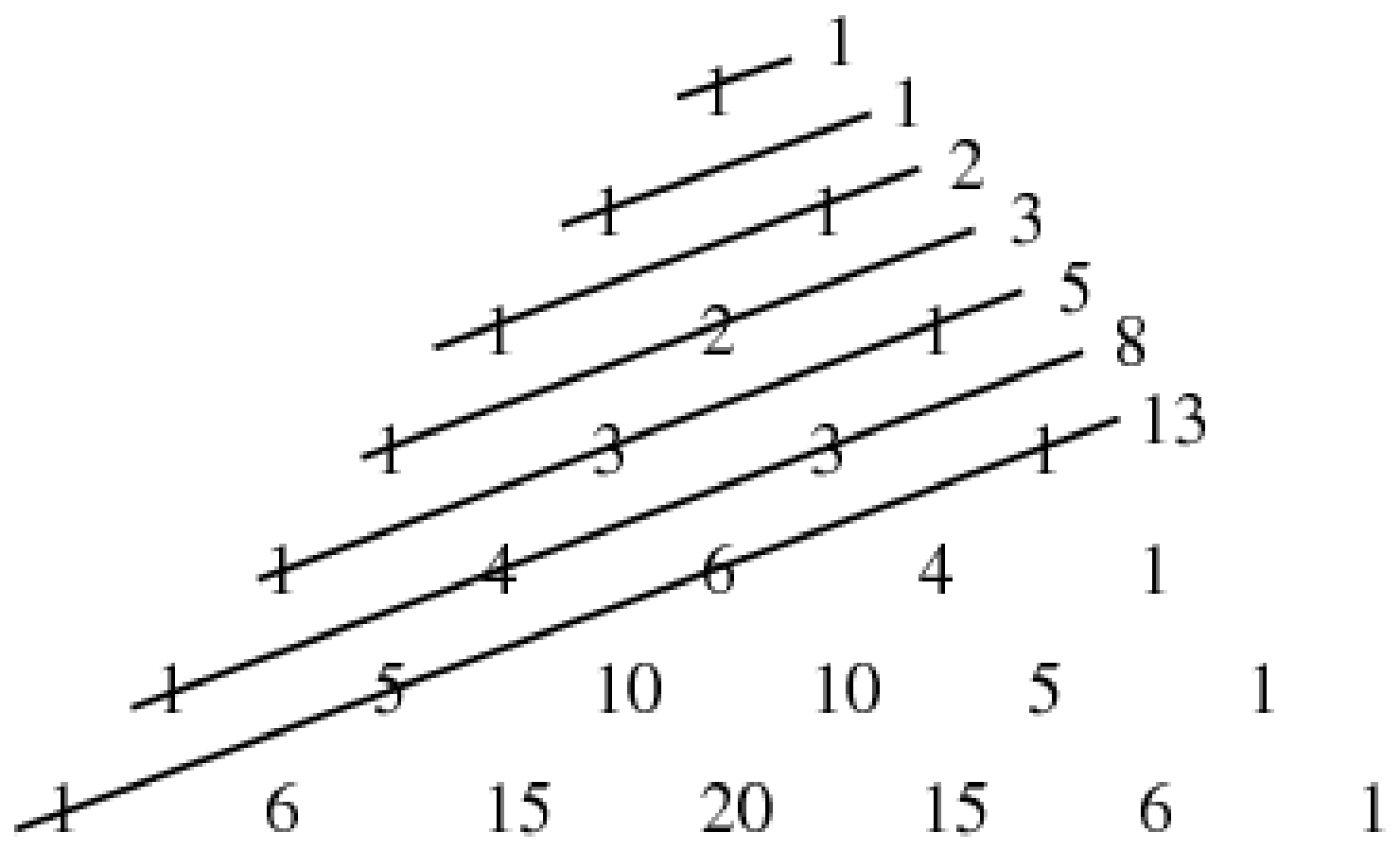

Altın oran Yunan harfi '$\Phi$' ile sembolize edilir. Aslında, Plato zamanının (İ.Ö. 400) matematikçileri onu anlam taşıyan bir değer olarak kabul ediyorlardı ve yunanlı mimarlar $1/\Phi$ oranını tasarımlarının ayrılmaz bir parçası olarak kullanıyorlardı. Bu mimari yapılardan en ünlüsü Atina'daki Parthenon'dur(Figure 0-27).

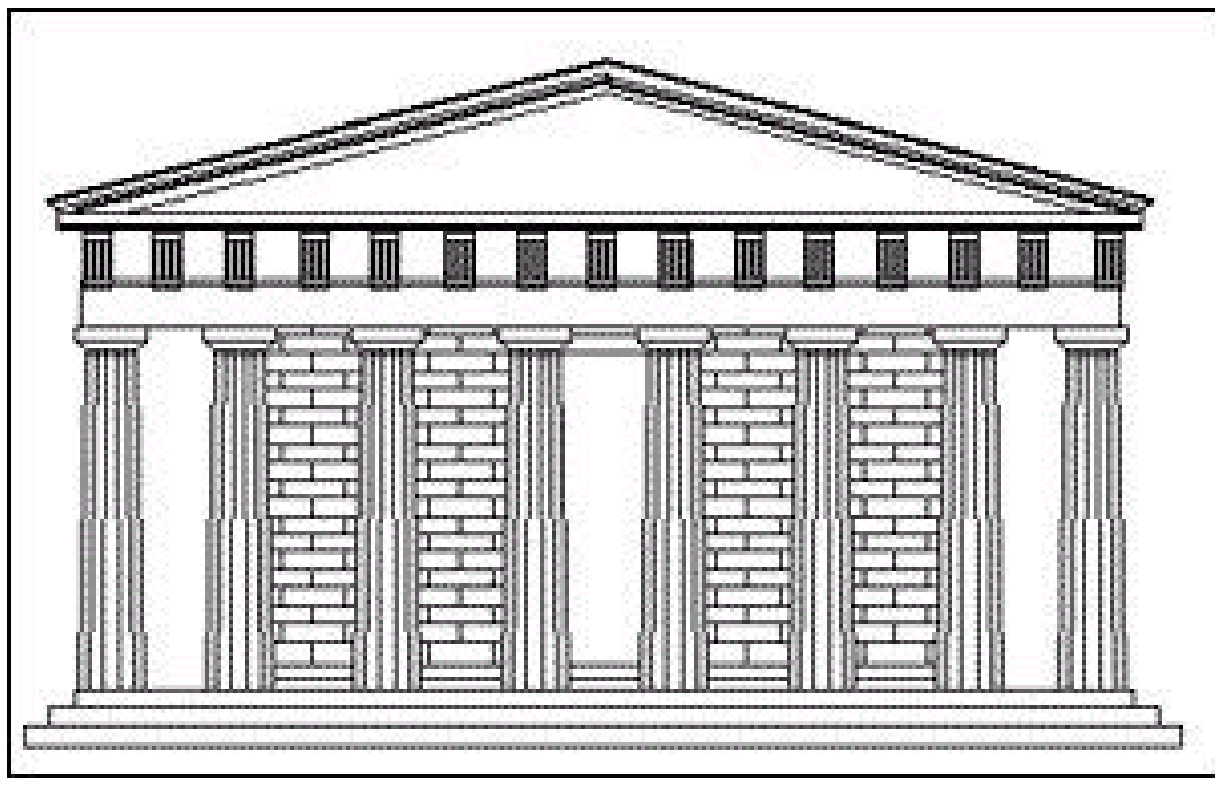

Figure 0-2: Parthenon,Atina

### 3.1 Büyük Piramit ve Altın Oran

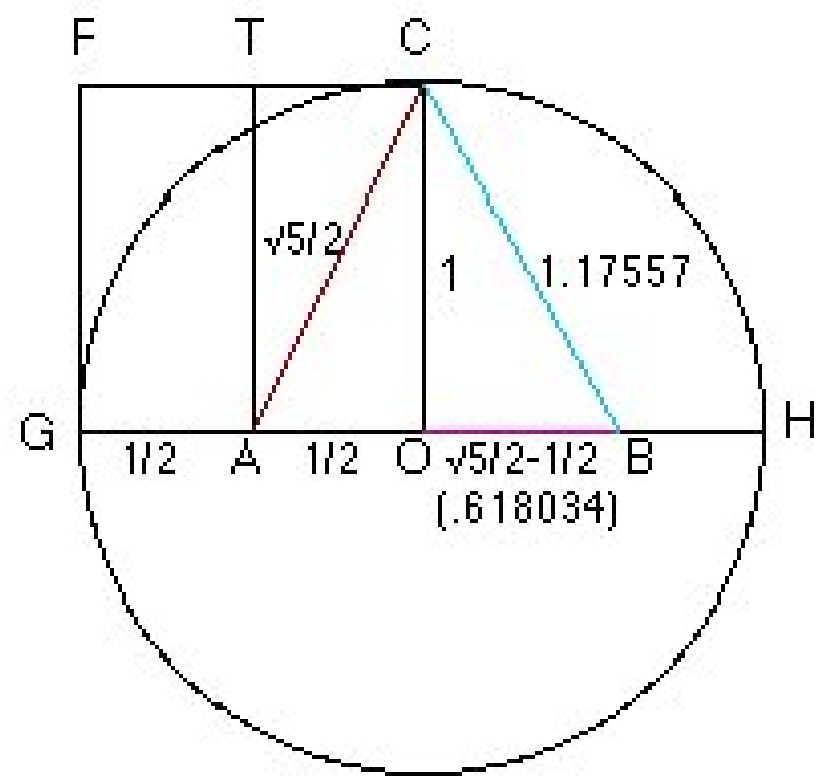

Yukarıdaki diagram, Altın Oran'ın bir çember yarıçapı üzerinde nasıl bulunabileceğini gösterir. Kenar uzunluğu dairenin yarıçapına eşit olan $FCGO$ karesinin $FC$ kenarının orta noktası olan $T$'den $GO$ kenarının orta noktası olan $A$'ya dik çizilen bir çizgi ile ikiye bölünmesinden elde edilen $TCAO$ dikdörtgeninin köşegenini ($AC$) bir ikizkenar üçgenin kenarlarından biri olarak kabul edip $ABC$ üçgenini oluşturursak, üçgenin yüksekliğini 1 kabul ettiğimizde (ki bu dairenin yarıçapıdır) $COB$ üçgeninin $OB$ kenarı, Altın Oran olan 0.618034 olur.

Bir trigonometrik cetvelden baktığımızda, $OCB$ açısının 31"43' ve dolayısıyla $OBC$ açısınında 58"17' olduğunu buluruz. Yukarıdaki diyagram önemini korumak şartıyla bizi başka bir konstrüksiyona götürür ki, bu belki de Mısır'lı rahiplerce çok daha önemli bulunmuş olabilir.

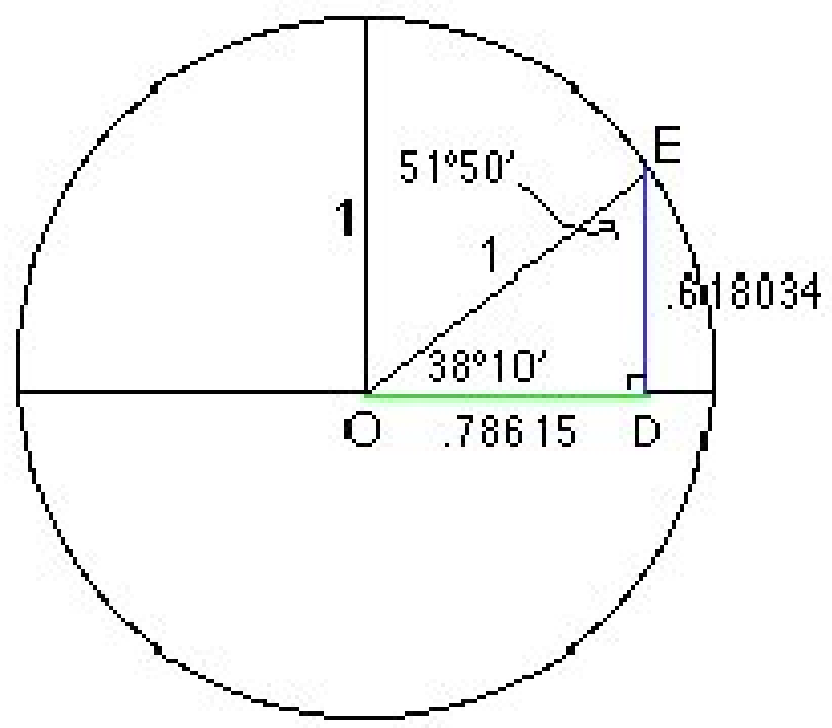

Üstteki diagramda,üçgenin dik açıya ortak kenarlarından biri yine yarıçapın 0.618034'üdür fakat bu defa 1'e eşit olan komşu kenar değil,hipotenüstür.Yine bir trigonometrik tablo yardımıyla, 0.618034'ün karşı açısının 38"10' ve diğer açının da 51"50' olduğunu görürüz. Pisagor Teoremini kullanarak, $OD$ kenarının uzunluğunun da yarıçapın 0.78615'i olduğu görülür.

Bu konstrüksiyonda onu özel yapan iki önemli nokta vardır. Birincisi; $ED$ kenarının uzunluğu (0.618034) $OD$ kenarının uzunluğuna (0.78615) bölünürse sonuç $OD$ kenarının uzunluğuna (0.78615) eşit çıkmaktadır. Trigonometrik ilişkiler açısından bu şu anlama gelmektedir: 38"10' un tanjantı , 38"10' un kosinüsüne eşittir. Tersi, 51"50' nin kotanjantı, 51"50' nin sinüsüne eşittir. İkinci ve belki en önemli husus: $OD$ kenar uzunluğu (0.78615) 4 ile çarpıldığında 3.1446 yı verir ki bu, hemen hemen Pi'ye (3.1416) eşittir. Bu buluş, 38"10' açıya sahip bir dik üçgenin Pi oranı ile Altın Oran fenomeninin çok özel ve ilginç bir kesişimini kapsadığını ortaya koymaktadır.

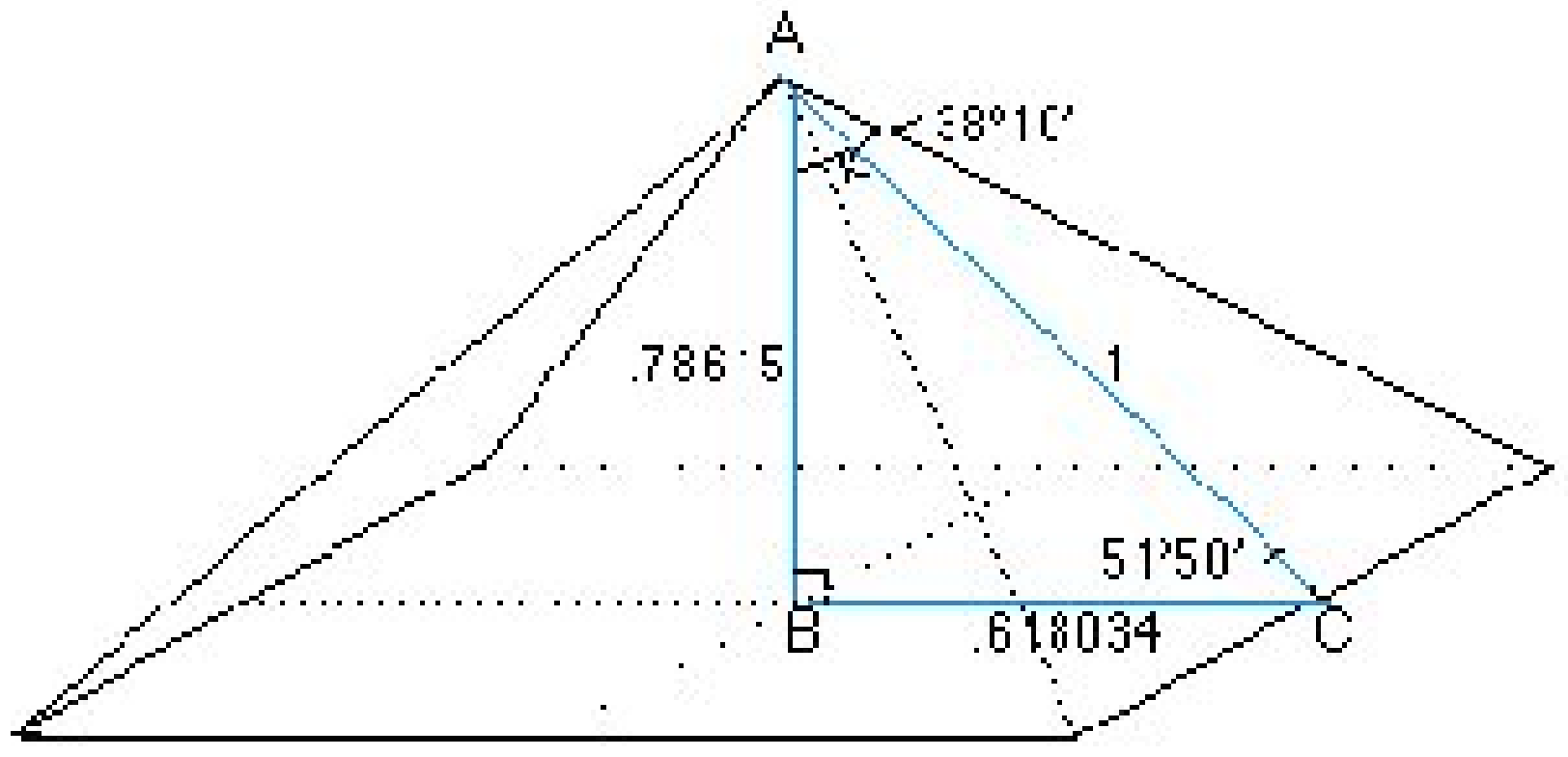

Kadim Mısır Krallığı döneminin rahipleri bu üçgenin özelliklerinden haberdar mıydılar? Bu diagram Büyük Piramit'in dış hatlarını göstermektedir. Bilinçli olarak ya da değil, bu piramit 38"10' lık bir üçgeni ihtiva edecek biçimde inşa edilmiştir. Yüzeyinin eğimi, çok kesin bir şekilde yerle 51"50' lık açı yapmaktadır. Bu piramit kesitini bir önceki ile kıyaslarsak, $BC$ uzunluğunun yarıçapın 0.618034'ü olduğunu, $AB$ uzunluğunun 0.78615 olduğunu ve $AC$ uzunluğunun 1 yani yarıçap olduğunu görebiliriz.Keops Piramidi'nin gerçek ölçüleri şöyledir (feet ölçüsü metreye çevrilmiştir):$AB = 146.6088m$ $BC = 115.1839m$ $AC = 186.3852m$.

Bu XXX noktadan itibaren işler biraz karmaşık ama çok çok ilginç bir hale gelmektedir.

Görüleceği gibi, $BC$ uzunluğu, piramitin kenar uzunluğunun yarısıdır. Bu nedenle piramitin çevresinin uzunluğu $BCx8$ dir. Yani piramitin relatif çevresi $0.618034x8 = 4.9443$ dür. Yine piramitin relatif yüksekliği 0.78615 in bir çemberin yarıçapı olduğu farzedilirse bu çemberin uzunluğu (çevresi) yine 4.9443 olacaktır.

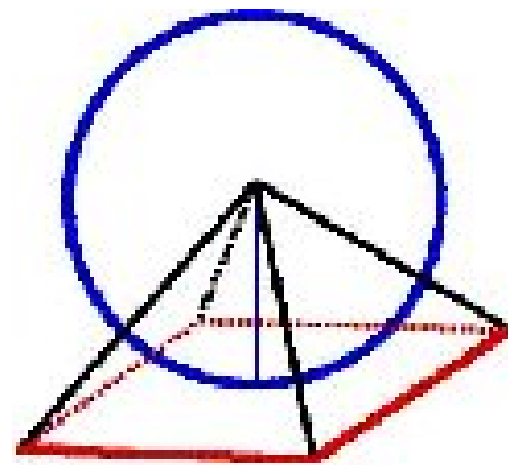

Bundan şu sonuç çıkmaktadır: Büyük Piramit, yatay bir düzlem üzerinden ölçüm yapıldığında sahip olduğu kare şeklindeki çevre uzunluğunun aynına, düşey bir düzlem üzerinde yapılan ölçümde de bu defa daire şeklinde olmak üzere sahiptir.

Birkaç ilginç bilgi olmak kaydıyla şu gerçeklere de kısaca bir göz atalım: Keops Piramidi'nin gerçek taban kenar uzunluğunun ($230.3465m$) 8 katı ya da çevre uzunluğunun iki katı, boylamlar arasındaki 1 dakikalık açının ekvatordaki uzunluğunu vermektedir. Piramitin kenar uzunluğunun, ekvatordaki 1 dakikalık mesafenin 1/8 ine eşit olması ve piramit yüksekliğinin 2 nin 1/8 ine eşit olması korelasyonunu irdelememiz, örneklemeyi evrensel boyutlara taşıdığımızda, dünya ile evrenin Pi ve Altın Oran sabitlerinin ilişkilerini algılamada küçük bir girişim, samimi bir başlangıç sayılabilir.

Şunu akılda tutmak gerekir ki; piramitin kenar uzunluğunun 230.3465m olması tamamen tesadüf de olabilir. Fakat karşılıklı ilişkiler yenilerini doğuruyor ve bunlara yenileri ekleniyorsa, bu korelasyonların kasti düzenlenmiş olduğu ihtimali de ciddi olarak dikkate alınmalıdır.

İnsan vücudunda da fibonacci sayılarına rastlanılır.

2 eliniz var, iki elinizdeki parmaklar 3 bölümden oluşur. Her elinizde 5 parmak vardır ve bunlardan sadece 8'i altın orana göre boğumlanmıştır. 2, 3, 5 ve 8 fibonocci sayılarına uyar.

Ayrıca aşağıdaki oranlar da altın orana eşittir.

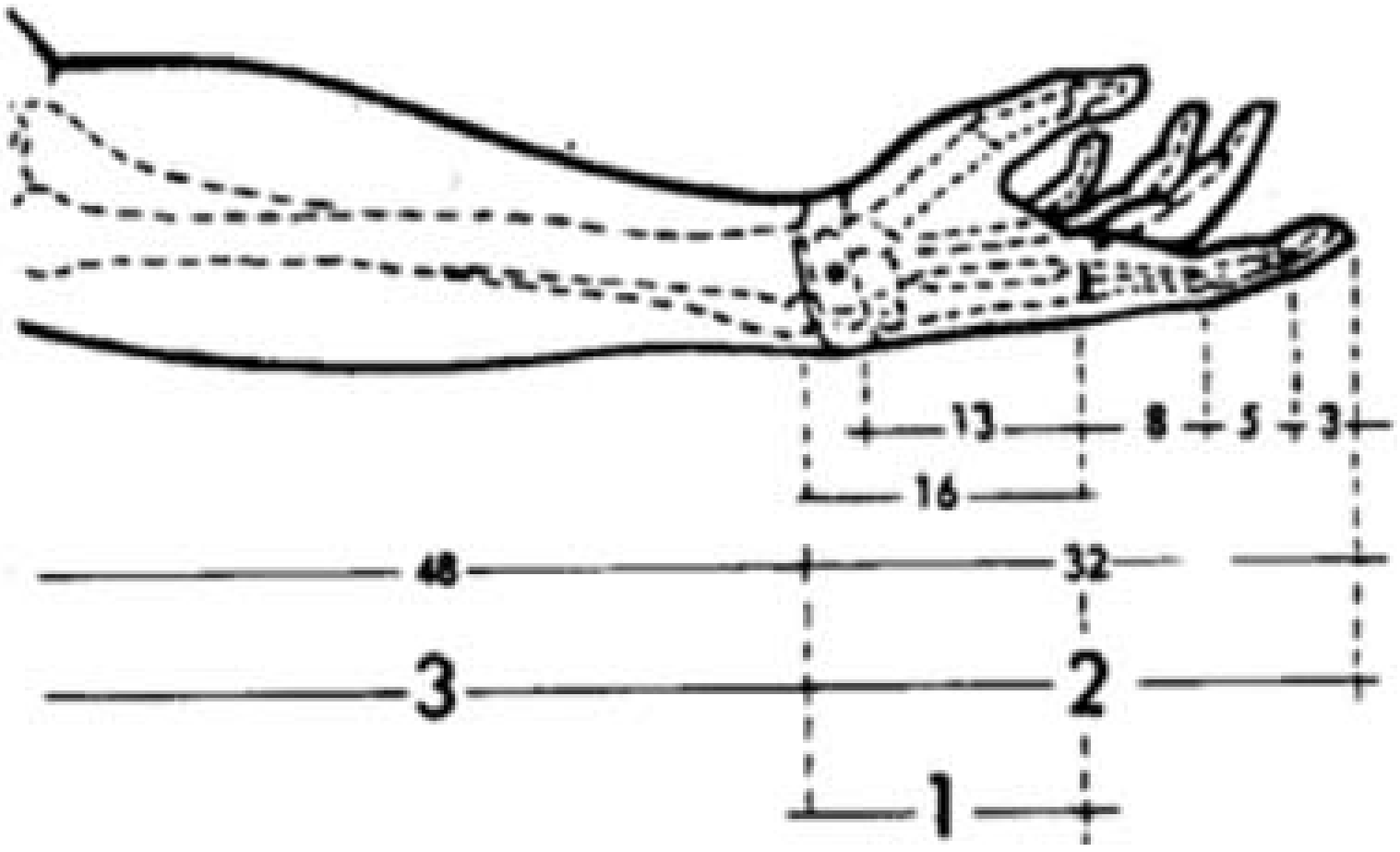

boy/ (bölü)bacak boyu

beden boyu/kolaltı beden boyu

tam kol boyu(boyun-parmak ucu)/dirsek - boğaz

göbek - omuz/göbek – bel

Fibonacci sayılarıyla bitki aleminde karşılaşmanın en çarpıcı örneklerinden biri ayçiçeği tohumlarında mevcut, saat ibresinin hareket yönünde ve buna karşı yönde uzayan iki tür spirallerin (sarmalların) sayısının ardışık iki Fibonacci sayısı olmasıdır. Orta büyüklükte ayçiçekleri için spirallerin sayısı 34 'e karşılık 55 veya 55 'e karşılık 89, daha büyükleri için 89 'a karşılık 144, ve küçükler için de 13 'e karşılık 21 veya 21 'e karşılık 34 olarak gözlenmiştir. Buna benzer bir durum papatya çiçeklerinde 21 'e karşılık 34, ananaslarda 8 'e karşılık 13, çam kozalaklarında 5 'e karşılık 8 veya 8 'e karşılık 13 olarak gözlenmiştir.

Bitki aleminde yaprakların saplar üzerindeki dizilişi (phyllotaxy) ile Fibonacci sayıları arasındaki ilişkiye dair çok sayıda örnek vardır. Örneğin 2/5 kesiriyle ifade edilen bir phyllotaxy, iki yaprağın sap boyunca aynı sıraya gelinceye kadar sap etrafında 2 tur yaptığını ve sap boyunca 5 tane yaprak olduğunu anlatmaktadır. Sap boyunca belli bir yapraktan sonra 6. yaprak aynı hizada olup, ardışık iki yaprak sap etrafında 720/5=144 derecelik açı yapmaktadır. Bazı bitkiler için bu oranlar:

karaağaç, çim 1/2

kayın, ayak otu 1/3

meşe, elma, armut 2/5

kavak, muz 3/8

badem, pırasa 5/13

olarak gözlenmiştir.

# 4 EKONOMİDEKİ UYGULAMALARI

Fibonacci sayıları ile ilgili daha birçok enteresan hesaplamalar vardır.Şöyle ki: Her bir oran, Fibonacci dizisindeki iki sayının binde birler seviyesinde toplamıdır. Biraz karışık görünse de örnekle daha rahat anlayabiliriz. Artan dizide 1.00 oranı 0.987 ile .013 sayılarının toplamıdır. Bir sonraki oran olan 1.618 sayısı 1.597 ile .021 sayılarının, 2.618 oranı 2.464 ile .034 sayılarının toplamıdır ve bu seri böylece gider. Azalan seride .618 oranı .610 ile .008'in farkı, .382 oranı .377 ile .005'in farkı, 0.236 oranı 0.233 ile .003'ün farkı, .146 oranı .144 ile .002'nin farkı, .90 oranı .89 ile .001'in, 0.56 oranı ise 0.55 ile 0.001'in farkıdır. Bundan sonraki her oran artık zaten bir Fibonacci sayısı olmuştur. Böylece on üçüncü sayıdan sonra başladığımız yere yani 0.001'e geri döndük. Nerden sayarsanız sayın her zaman aynı fenomen karşınıza çıkar. İşte bu sürekli kendini tekrar etme ve kendi kendinden yeniden yaratılma Elliott ile Fibonacci arasındaki vazgeçilmez ilişkidir.

**4.1 Fibonacci Oranları:**

Dalga kuramı üç bölümden oluşur; dalga biçimi, oran ve zaman. Şimdi, Fibonacci oranları ve geri-çekilmeleri üzerinde duralım. Bu ilişkiler; fiyata uygulanışı daha güvenilir olduğu kabul edilse de, fiyat ve zamanın ikisine de uygulanabilir.Temel dalga biçiminin Fibonacci sayılarına bölündüğünü görebiliriz. Tam bir çevrim beş yukarı ve üç aşağı olmak üzere sekiz dalgayı kapsar ve bunların tümü Fibonacci sayılarıdır. Bununla birlikte Fibonacci sayı dizinin dalga kuramının matematik temelleri dalga saymanın da ötesine gider. Değişik dalgalar arasındaki oransal ilişkiler sorunu vardır. Aşağıdakiler en genel kullanılan Fibonacci oranları arasındadır.

1) Bir çevrimde üç itki (impulse) dalgadan yalnızca biri uzadığı için, diğer iki dalga zaman ve büyüklük olarak aynıdır. Eğer uzayan dalga 5. dalga ise, 1. ve 3. dalgalar yaklaşık olarak eşit olurlar. Eğer 3. dalga uzarsa, 1. ve 5. dalgalar eşit olma eğilimi taşırlar.

2) 3. dalganın tepesinin minimum hedefini bulabilmek için 1. dalganın uzunluğu 1.618 ile çarpılır ve bu değer 2. dalganın tabanına eklenir.

3) 5. dalganın tepesi, 1. dalganın 3.236 (2X1.618) ile çarpılması ve bu değerin 1. dalganın tepesine ya da tabanına maksimum ve minimum hedefleri bulabilmek için eklenmesiyle belirlenebilir.

4) 1. ve 3. dalgalar yaklaşık olarak eşit olurlarsa ve 5. dalganın uzaması bekleniyor ise fiyat hedefi 1. dalganın tabanından 3. dalganın tepesine kadar olan uzaklığın ölçülüp, bu uzaklığın 1.618 ile çarpılıp elde edilen değerin 4. dalganın tabanına eklenmesiyle bulunur.

5) Düzeltme dalgaları için, normal bir 5-3-5 zig zag düzeltmesinde "C" dalgası "A" dalgasının uzunluğuna çoğunlukla yaklaşık olarak eşit olur.

6) "C" dalgasının muhtemel uzunluğunu ölçmenin diğer bir yolu; "A" dalgasının uzunluğunu önce 0.618'le çarpmak daha sonra bu sonucu "A" dalgasının tabanından çıkarmaktır.

7) 3-3-5 yatay düzeltmede, "B" dalgasının "A" dalgasının tepesine yetiştiği ya da geçtiği durumlarda "C" dalgası "A" dalgasının uzunluğunun yaklaşık 1.618 katı kadar olur.

8) Bir simetrik üçgende, birbirinin arkasından gelen her bir dalga, kendinden bir önceki dalgaya yaklaşık 0.18 sayısı ile orantılı bir ilişki içindedir.

## 4.2 Fibonacci Yüzde Geri Çekilmeleri

Daha başka çeşitli oranlar da vardır fakat yukarıda sıralananlar en çok kullanılanlardır. Oranlar hem itki (impulse) dalgalarda hem de düzeltme (corrective) dalgalarında fiyat hedeflerini bulmakta yardımcı olurlar. Fiyat hedeflerini bulmanın bir diğer yolu, yüzde geri-çekilmelerin kullanımıdır. Geri-çekilme analizlerinde en çok kullanılan sayılar, % 61,8 (% 62'ye yuvarlanır), % 38.2 (%38'e yuvarlanır) ve % 50'dir. Güçlü bir trendde, minimum geri-çekilme genellikle % 38 civarındadır. Daha zayıf bir trendde, maksimum yüzde geri-çekilme çoğunlukla % 62'dir. Daha önce de işaret edildiği gibi Fibonacci oranları, ilk dört sayıdan sonra 0.618 'e yaklaşırlar. (1/3'lük geri-çekilme, Fibonacci oranının bir alternatifi olarak aynı zamanda Elliott'un kuramının da bir parçasıdır.) Bir önceki boğa ya da ayı piyasasının tümüyle geri-çekilmesi de (% 100) önemli bir destek ya da direnç bölgesine işaret eder.

Senedin fiyat hareketlerinde Fibonacci geri alma seviyeleri destek ve direnç olarak çalışmaktadır.

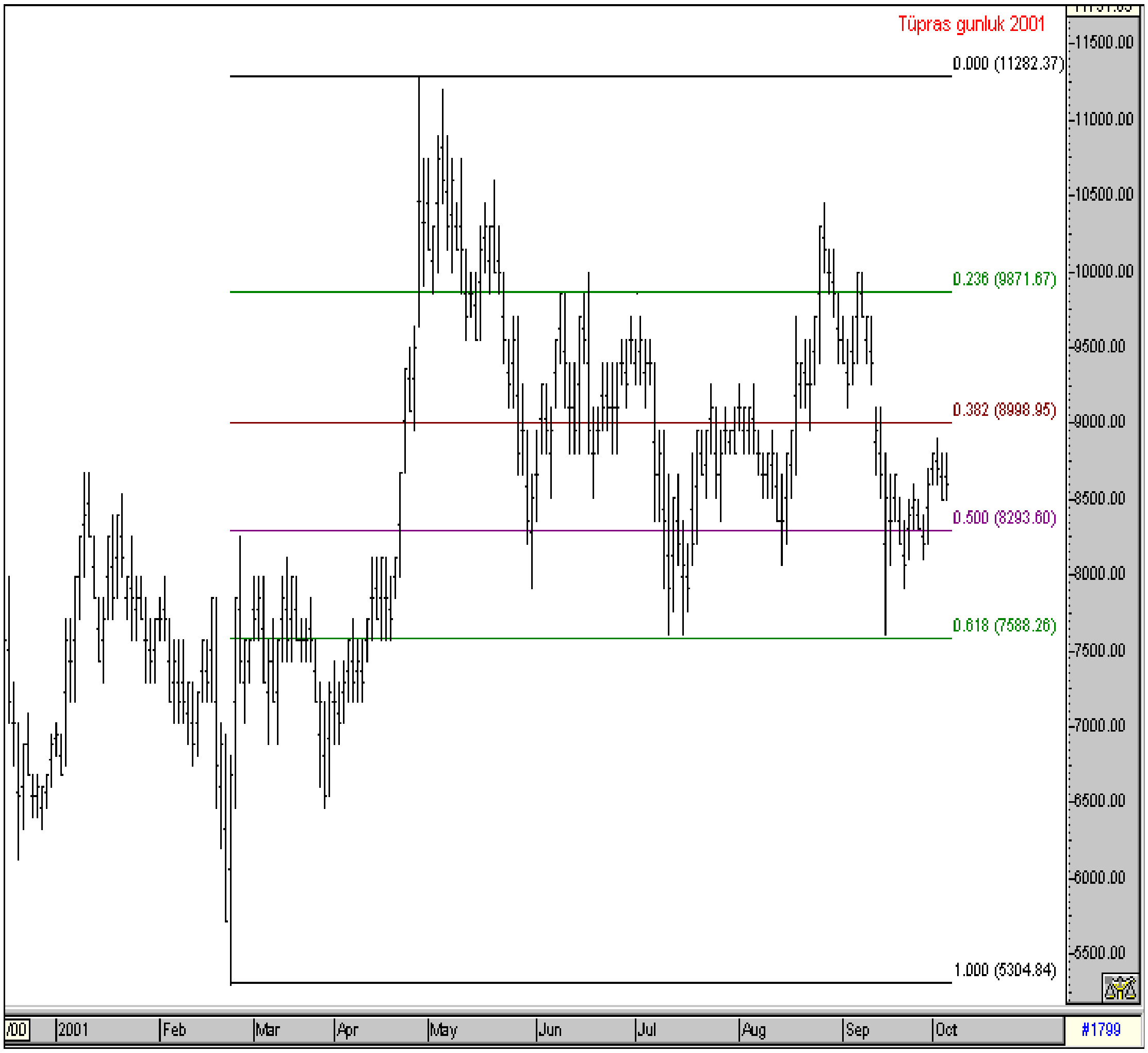

Figure 0-3: TÜPRAŞ

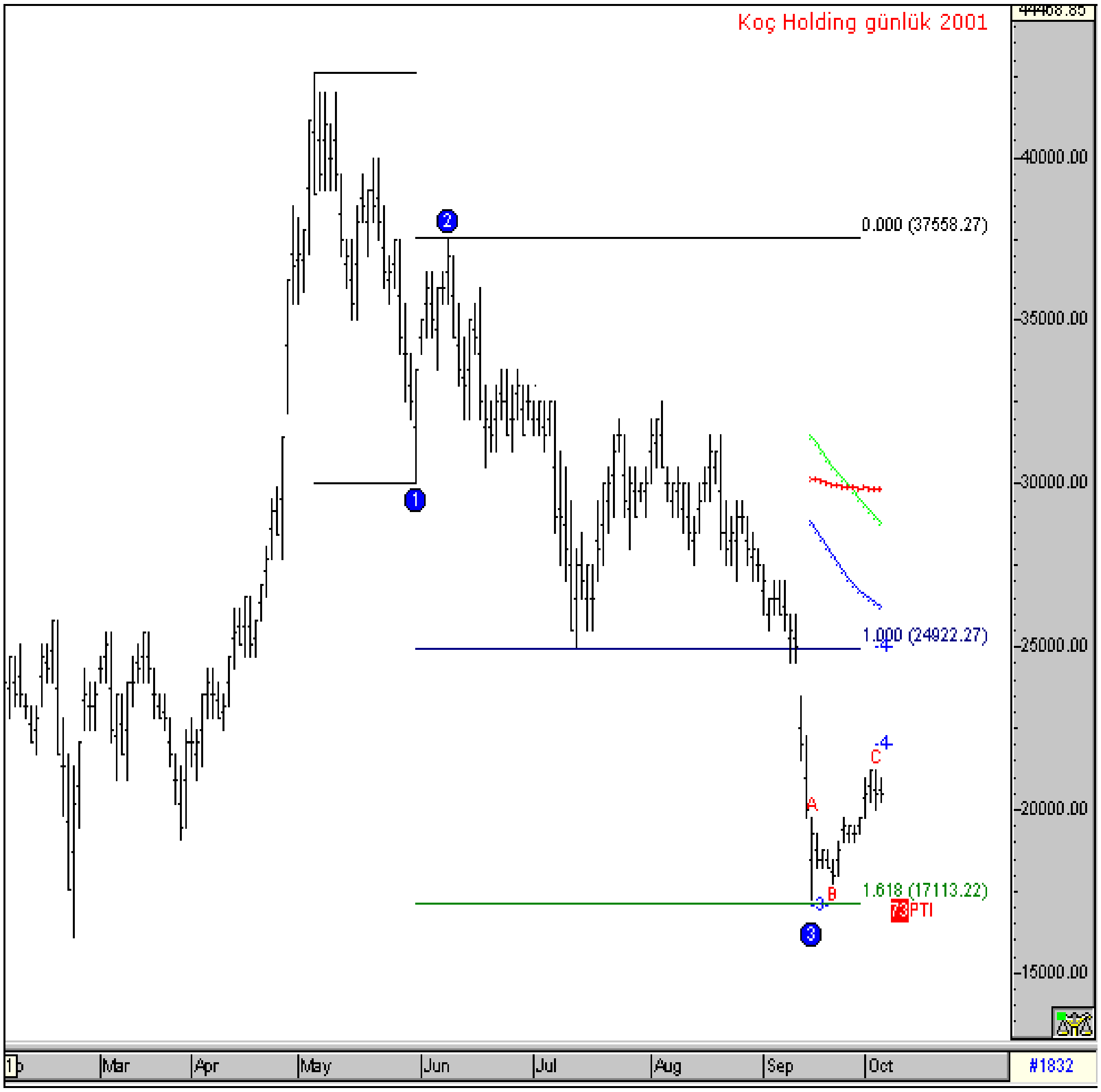

Figure 0-4: KOÇ HOLDİNG

Düşüş trendinde olan senedin 2. dalgası oluştuktan sonra 3. dalga hedefleri Fibonacci oranları ile belirlenir. İlk hedef, 1. dalga 3. dalgaya eşit olması. Bu seviyeden tepki verir. Düşüş 1. dalganın tam katında koparak 1.618 katında durmuştur. 1.618 katından başlayan hareket yine Fibonacci oranları kadar geri alacaktır.

### 4.3 Fibonacci Zaman Çevrimleri

Fibonacci zaman dönemlerini belirlemek için grafiklerdeki önemli bir taban ya da tepeden (Pivot alınarak) ileriye doğru Fibonacci sayılarını sayarak gelecekteki Fibonacci zaman dönemlerini belirlenebilir. Bir başlangıç noktasından sağa doğru 1, 2, 3, 5, 13, 21, 34, 55, 89, 144, 233 vb. günlere düşey çizgiler çizilir. Amaç, Fibonacci zaman hedefleri üzerindeki ya da yakınlarındaki trend değişikliklerini belirleyebilmektir.

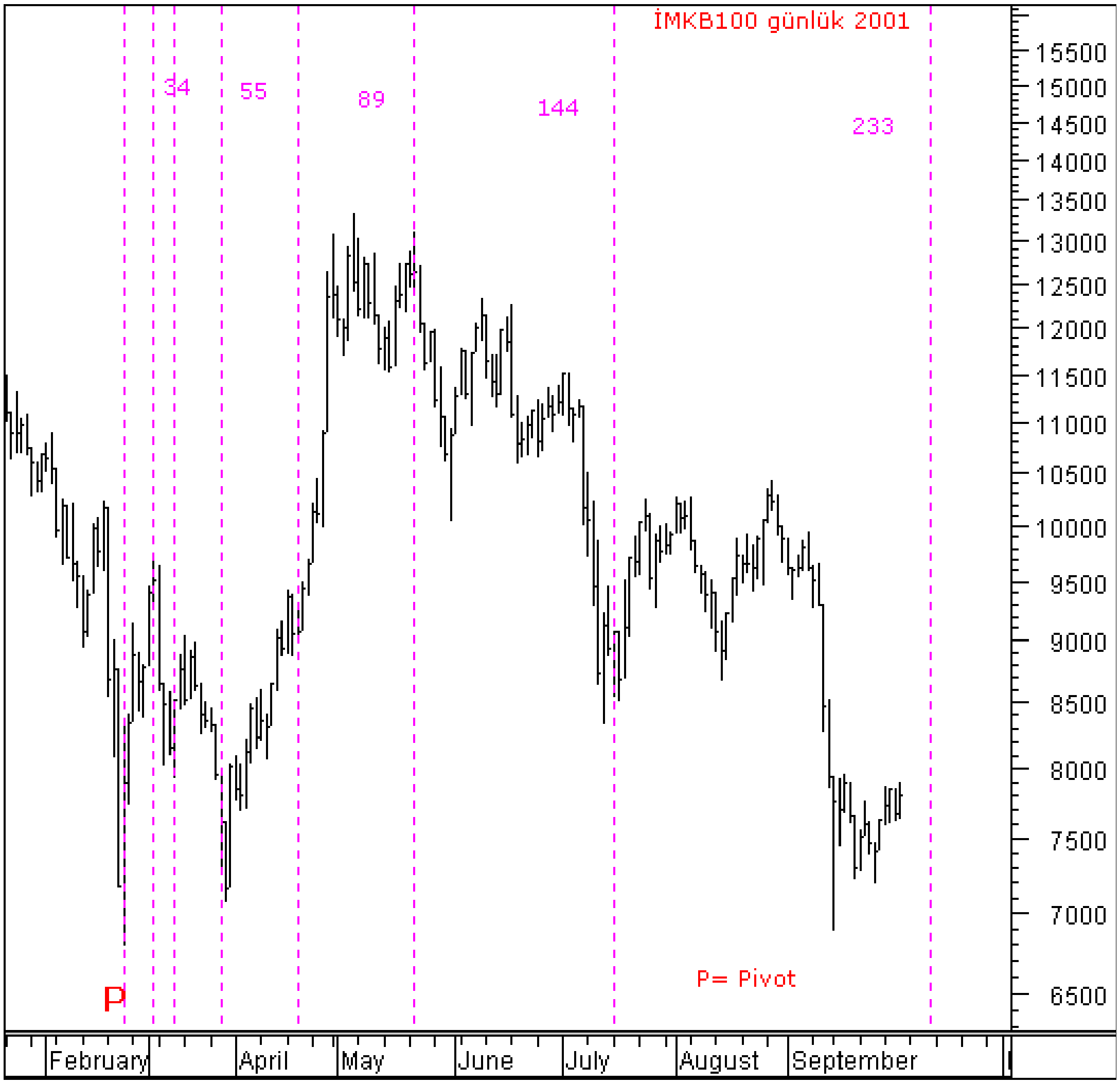

Figure 0-5: İMKB 100

İMKB 100 endeksinin 22.02.2001 tarihinden itibaren Fibonacci sayılarına karşılık gelen günler belirlenmiştir. Bu tarihlere yakın günlerde endekste fiyat değişimleri yaşanmaktadır. Bir başka örnek olarak, Fibonacci sayıları hareketli ortalama analizlerinde çokça kullanılır. Bu şaşırtıcı olmamalı çünkü en başarılı ortalamalar çeşitli piyasaların baskın zaman çevrimlerine bağlı olanlardır.

**4.4 Almaşıklık Kuralı - Rule of Alternation**

Ralph Nelson Elliott'ın almaşıklık kuralı olarak tanımladığı ilke çok kabaca şu varsayıma dayanır :

" Kuvvetli tepkiyi zayıf tepki, kuvvetli ralliyi, zayıf ralli izler." Glenn Neeley

Mastering Elliott Wave kitabında bu yaklaşımı daha ileri götürüyor ve almaşıklığın 5 biçimde görülebileceğini söylüyor:

1. Fiyat (Price) - fiyat mesafesinde katedilen mesafe
2. Zaman (Time) - zaman mesafesinde katedilen mesafe
3. Şiddet (Severity) - Bir önceki dalganın geri alış mesafesi
4. Karışıklık (Intricacy) - Bir modeldeki alt bölüntüler
5. Karmaşıklık (Complexity) - Bir dalganın yatay ya da keskin olarak ilerlemesi

Yukarıdaki yaklaşımdan hareketle endeksin 2000 yılında başlattığı düşüşü bir kez daha analiz edelim :

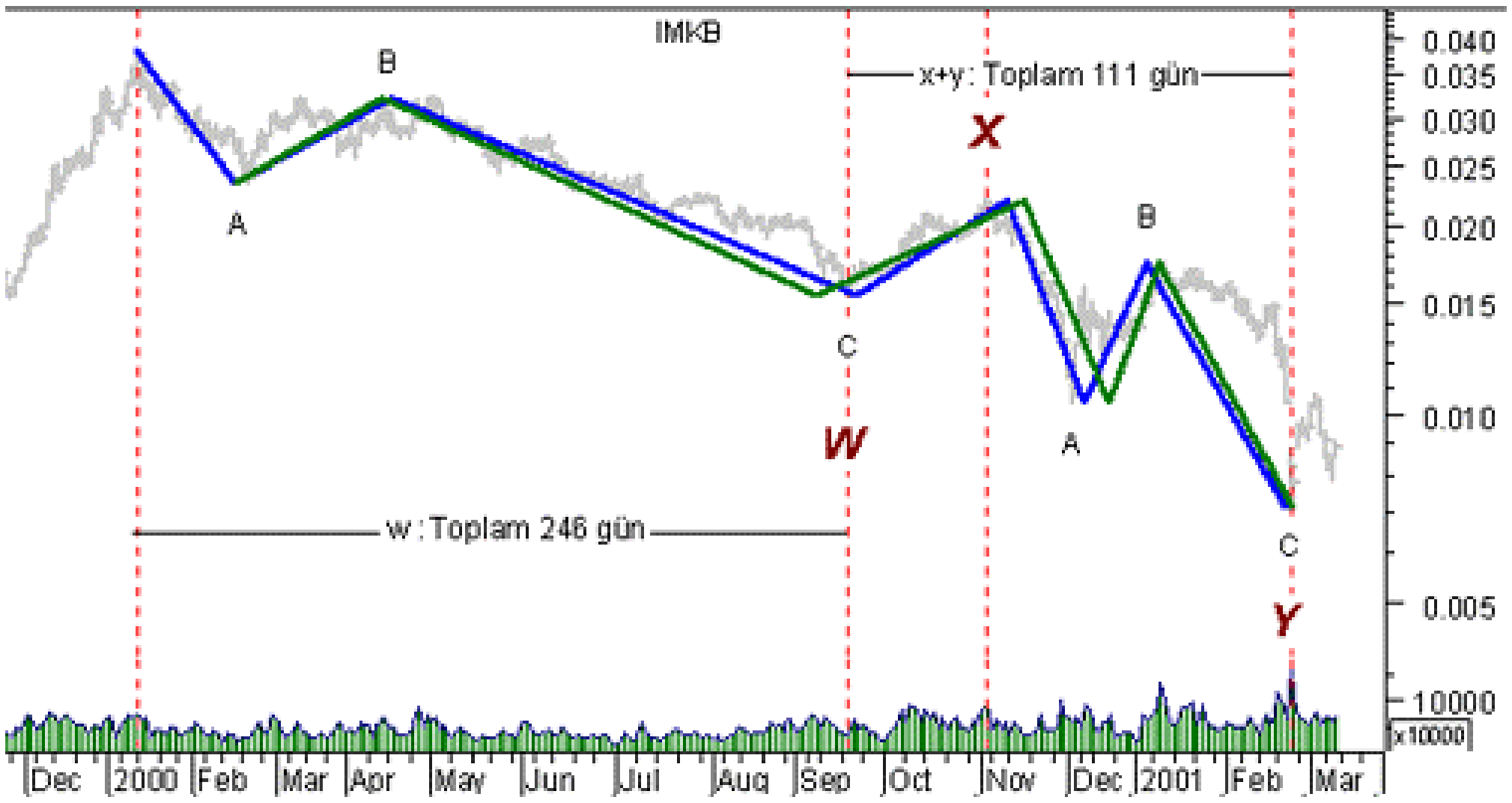

Yukarıdaki analizdeki modellemede, endeks tam Fibonacci dönüş günlerinde yön değiştirseydi hareketin nasıl ilerleyeceği görülüyor. Grafik üzerindeki çizgiler, hep, trendin döndüğü tarihler korunarak ve trend dönüşleri Fibonacci sürelerinde gerçekleşseydi oluşacak modeli göstermektedir. Bu modellerde dikkat edilmesi gereken, endeksin tam dönmesi gereken günlerin değil, dönüşün gerçekleştiği tarihin Fibonacci modelinde öngörülen gün olması varsayımına göre çizilmiş olmalarıdır. Yani, ilk A dalgası 23 Şubat 2000'de tamamlandıktan sonra, ilk B dalgası 69 takvim günü sürmüştü. Bunun yerine model, Fibonacci sayısı olan 55 günde tamamlansaydı ne olurdu yaklaşımı ile çizilmiştir. Buradaki yaklaşımın endeksin muhakkak dönmesi gereken gün gibi bir yaklaşıma dayanmamaktadır.

Zaman analizinde ilginç olan $W$ dalgasının toplam 246 takvim süresi sürmüş olmasıdır. O halde daha sonra oluşan $X$ ve $Y$ dalgalarının toplam $246 \cdot 0.618 = 152$ takvim gününde oluşması gerekirdi. 22 Şubat 2001'de 0.71 cent görülürken, ilk düşüş dibinden sonra tam 155 takvim günü geçmişti. Buradan çıkartılacak sonuç : Endeksin 19 Şubat 2001 tarihinde çok önemli bir Fibonacci dönüş zamanına ulaşmış olmasıdır. Ancak bu çıkarsama hiçbir türlü, endeksin bu seviyeden dönmesi gerektiği anlamına gelmez. Bu tarihe ulaşıldığında endeksin bu seviyede dip yapması sadece beklenebilir. Endeks dip yapamaz ve düşüşe devam ederse, bir sonraki zaman hedefine yöneldiği varsayılabilir. Burada almaşıklık kuralının öngördüğü, uzun zaman süren düşüşü daha kısa sürecek bir düşüşün izlemesi gerektiğidir. Buradaki ilişkinin tam altın oran olan 0.618 olması ise beklenendir. (Glenn Neeley'in Mastering Elliott Wave kitabında bir zigzag düzeltmede zaman ilişkisinin $W = X + Y$ veya $W = 0.618 \cdot (X + Y)$, ya da $W + X = Y$ veya $0.618 \cdot (W + X) = Y$ olması gerektiği öngörülmektedir.)

Not : Yukarıdaki analiz, endeksin bir double zigzag düzeltme yaptığı ve bu düzeltmenin tamamlandığı varsayımına dayanmaktadır. Doğaldır ki, endeks yeni fiyat dibine dönerse dalga sayımı ile beraber tüm zaman ve fiyat modelleri de değişecektir.

Şimdi de analizi fiyat değişiminde yapalım:

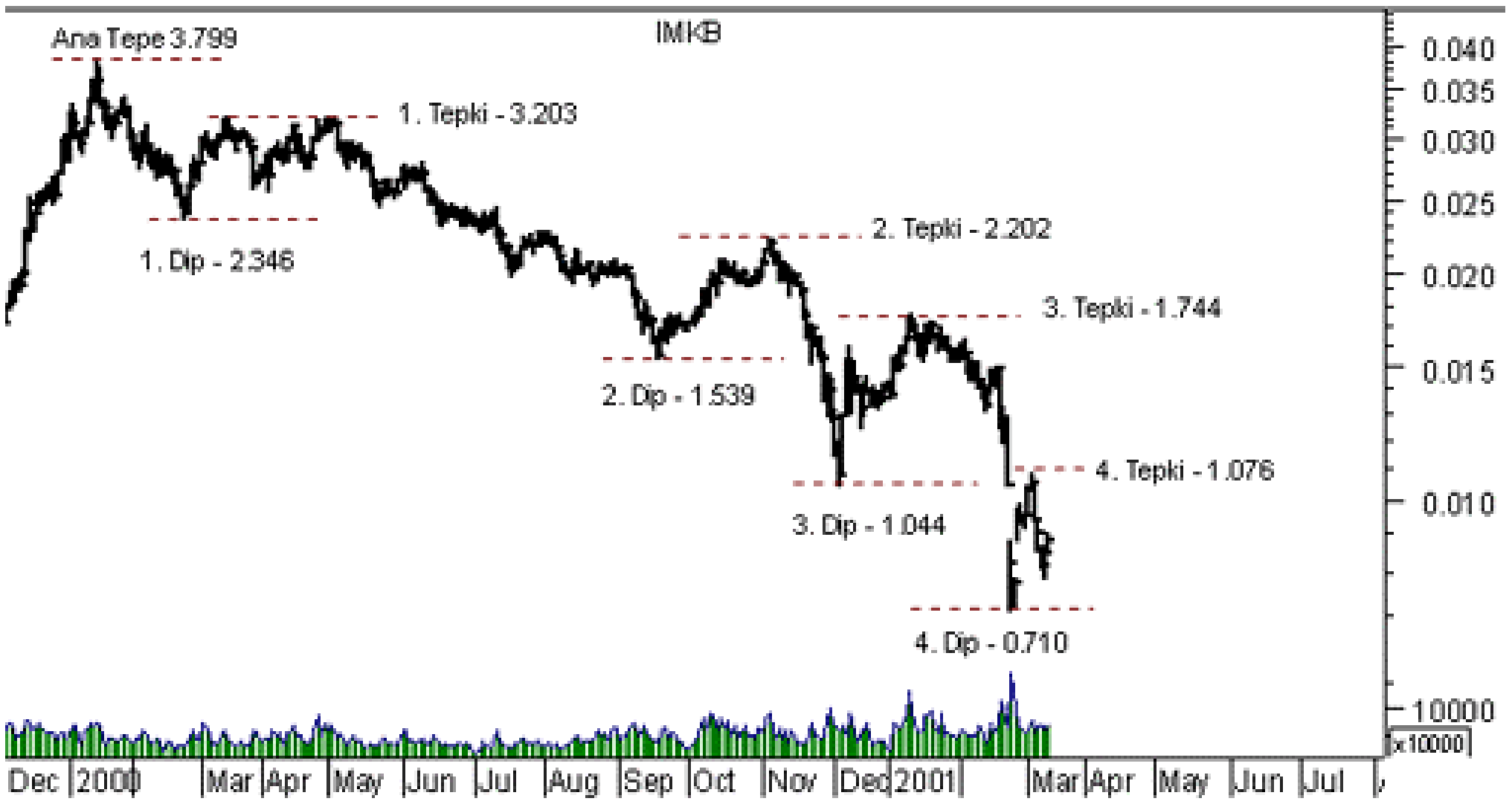

Önceki grafikte endekste şu ana kadar oluşan 4 dip ve bu 4 dipten gelen tepkilerin katettikleri fiyat mesafeleri görülmektedir. Endekste 4 tepki dönemindeki fiyat değişimi şu sekilde oluşmuştur :

1. Dip düşüş mesafesi = 2.346 - 3.799 = -1.453

1. Tepki fiyat mesafesi = 3.203 - 2.346 = 0.857

Tepki mesafesi / Düşüş mesafesi = 0.857 / 1.453 = 0.590 ( ~0.618)

2. Dip düşüş mesafesi = 1.539 - 3.203 = -1.664

2. Tepki fiyat mesafesi = 2.202 - 1.539 = 0.663

Tepki mesafesi / Düşüş mesafesi = 0.663 / 1.664 = 0.398 ( ~0.382)

3. Dip düşüş mesafesi = 1.044 - 2.202 = -1.158

3. Tepki fiyat mesafesi = 1.744 - 1.044 = 0.700

Tepki mesafesi / Düşüş mesafesi = 0.700 / 1.158 = 0.604 ( ~0.618)

4. Dip düşüş mesafesi = 0.710 - 1.744 = -1.034

4. Tepki fiyat mesafesi = 1.076 - 0.710 = 0.366

Tepki mesafesi / Düşüş mesafesi = 0.366 / 1.034 = 0.354 ( ~0.382)

Görüldüğü gibi, son derece net bir şekilde güçlü tepkiyi, zayıf tepki izleyerek Almaşıklık kuralı çok küçük yüzde sapmalarla hemen hemen mükemmel Fibonacci oranlarında gerçekleşmiş.

Buraya bir küçük not düşelim ve bir yükseliş başlarsa bu yükselişin ya yeni bir trend başlangıcı, ya da almaşıklık kuralına göre güçlü bir tepki olması gerektiğini belirtelim.

Üçüncü olarak tepki yükselişlerini karışıklık ve karmaşıklık ilkeleri doğrultusunda inceleyelim. Yani tepki yükselişlerini oluşturan dalgaların bir alt dereceli dalgalarına ve hangi modellerde ilerlediklerine bakalım:

1. Tepki kendi içinde bariz bir bölüntü ile ilerlemiş ve yeni bir fiyat tepesine ulaşamamış, yani yükseliş bir yassı (flat) düzeltme olarak kalmış. Bir sonraki tepkide bu kez alt dereceli dalgalar kendi içlerinde zaman ve fiyat mesafesi olarak ciddi oranlarda bölünmeden daha basit bir düzeltme olarak ilerlemiş. 2. tepkide ilk tepenin görülmesinin ardından ikinci bir tepe daha görülmüş. 3. tepki yükselişi kendi içinde bariz bir bölüntü ile yeni bir tepeye ulaşan bir zigzag olarak ilerlerken bir sonraki tepki, yani 4 No'lu tepki eğer yeni bir trende dönüşebilirse bir tek dalgada oluşmuş olacak, eğer 1.076 cent tepesi geçilemezse bu kez de yeni tepe yapamayan bir yassı (flat) düzeltme olarak kalması en muhtemel gelişim olacaktır.

Yukarıdaki analiz zaman ve fiyat analizlerinde Fibonacci sayıları ve oran-larının dalga oluşumunda nasıl gerçekleştiğini ortaya koymaktadır. Fibonacci sayıları ile çok daha sofistike analizler yapmak mümkün. Dalga yapısının ilerlemesi, Fibonacci fanları, Fibonacci zaman bölgeleri, Fibonacci Spiralleri ve Fibonacci box ile incelenebilir.

Böylece neden tüm fiyat grafiklerinin birbirlerine benzediği de anlaşılmış oluyor. Irkı, dili, kültürü, içinde yaşadığı sosyal çevresi, aldığı eğitimi, refah düzeyi ne olursa olsun tüm insan kitleleri aynı finansal davranış kalıpları içinde hareket ederler. Kaotik görünen sistemlerin özünde düzenli yapılar vardır ve bu yapılar daha geniş bir zaman boyutunda incelendiklerinde hep aynı sayısal ilişkiyi vermektedir : 0.618 veya (1-0.618) = 0382. Tüm büyüme, gelişme ve ilerleme modellerinin özünde aynı ilişki vardır. Bu ilişki sayısal olarak incelenebilir ve bu sayısal ilişki de Fibonacci sayı serileridir.

**4.5 Fibonacci Box**

Bir piyasa, birbirinden farklı düşünen ve farklı vadelerde işlem yapan binlerce yatırımcıdan oluşur. Her bir yatırımcının, fiyatların gelecekteki yönü konusundaki farklı görüşü, fiyatların bir trend çizgisinden diğerine taşınmasına neden olur.

Fiyat değişimi, daha önceki fiyat hareketlerinin bir sonucudur ve bir yükseliş ya da düşüşün katettiği fiyat mesafesi ve oluşma zamanı, kendisinden sonra gelecek olan fiyat hareketlerinin de belirleyicisidir.

Geçmiş fiyat hareketinden hareketle, gelecekte etkili olması mümkün trend çizgilerini ve destek ve dirençleri tespit etmek mümkündür.

Aşağıda detayları tarif edilen teknik zaman-fiyat analizlerinin bir parçası olarak ortaya çıktı. Bu tekniğin varsayımlarını ve neden-sonuç ilişkilerini kavrayabilmek için önce diğer analiz tekniklerine bir göz atalım :

**4.5.1 Trend çizgileri**

Bir trend çizgisini çizebilmek için, iki belirgin dip ya da tepenin fiyat grafiğinde belirmesi gerekir. Bu iki noktadan çizilen trend çizgisi, ilerleyen bir hareketi izlemek ve pozisyon değişimlerini gerçekleştirmek için gereklidir. Ancak aşağıdaki örnekte de görüldüğü gibi, iki noktadan çizilen trend çizgisi, özellikle trend özelliği taşımayan dönemlerde, üçüncü denemede kırılıyor.

IMKB endeksinin 2000 yılı içindeki hareketinden alınan aşağıdaki örneklerde,

iki noktadan çizilen yükselen ve alçalan trend çizgilerinin sonuçları görülüyor.

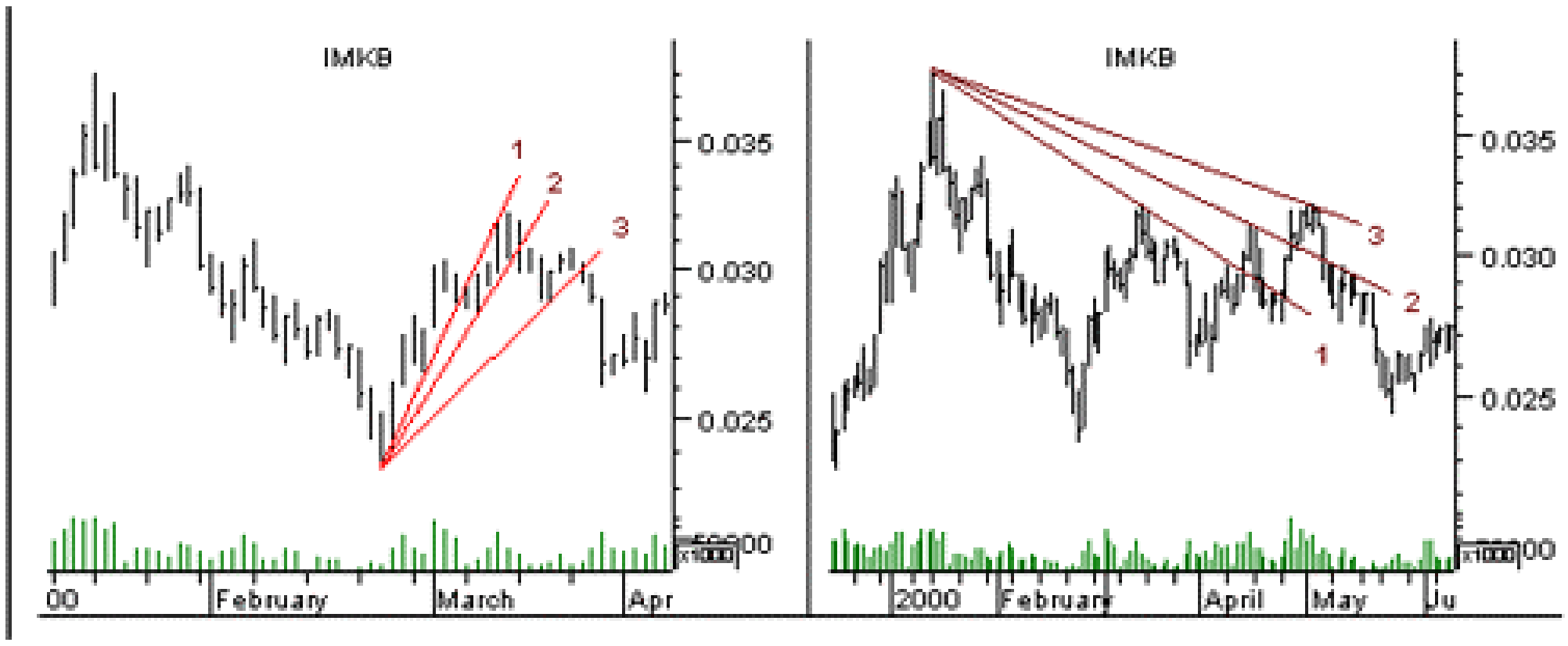

### 4.5.2 Hız direnç çizgileri

Bir tepe ile bir dip arasında çizilen hız çizgileri, iki noktadan çizilen trend çizgilerine göre daha doğru sonuçlar verme eğilimindedir. Ancak oynaklığın arttığı ve fiyat hareketlerinin düzenli trendler izlemediği dönemlerde hız çizgileri de yanıltıcı sonuçlar verebilir. Aşağıda aynı dönemden seçilen örneklerde de görüldüğü gibi, hız çizgileri ile yapılan analizlerde de İMKB'nin aynı döneminde yapılan değerlendirmeler hatalı sonuçlar verebilmektedir.

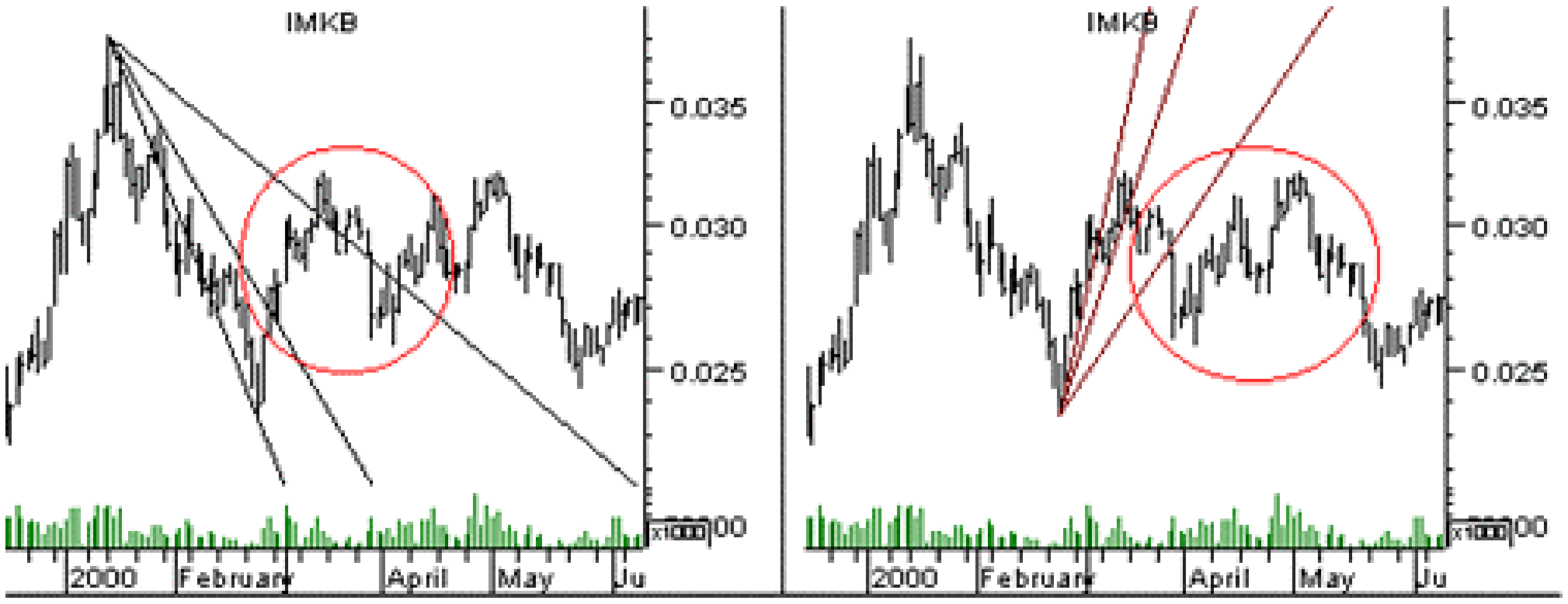

### 4.5.3. Fibonacci Fan

Fiyat hareketini izlemek için kullanılabilecek yöntemlerden biri de fibonacci fan çizmektir. Ancak benzer şekilde fibonacci fan da, hatalı ve yanıltıcı sinyaller verebilmektedir.

### 4.5.4 Gann Fan

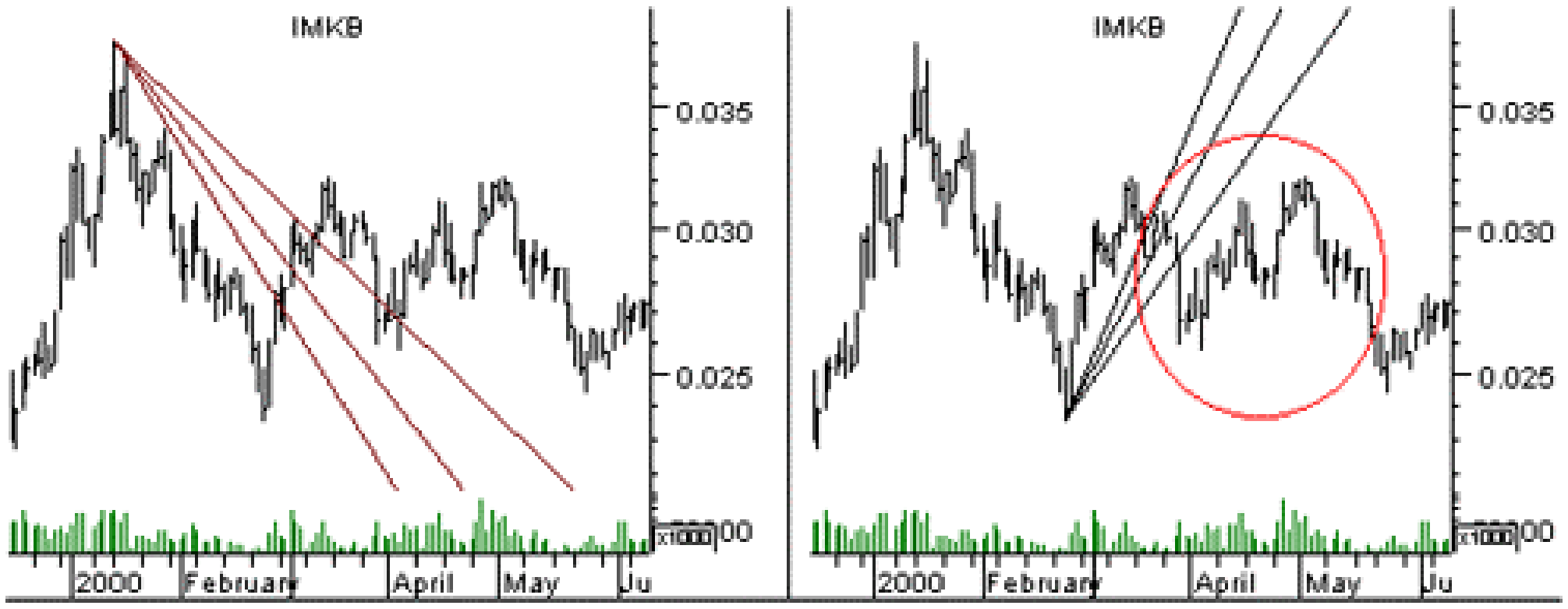

Diğer yöntemlere göre daha doğru sonuçlar verme eğiliminde olan Gann fan kullanmak iki açıdan problem doğurabilir. Birincisi, Gann fan, incelenen grafiğin ekran üzerinde kapsadığı alana göre değişkenlik göstermektedir. Bu nedenle, grafiğe zoom in, ya da zoom out yapıldığında fanın konumu değişkenlik göstermekte, her çizilişte fan ekrandaki görüntüye bağlı olarak değişmektedir. Bu nedenle önerilen, gann fanın mükemmel bir kare üzerinde çizilmesi ve 1x1 gann açısının bu karenin diyagonalinden geçirilmesidir. Özetle doğru çizilmeyen Gann fan, hatalı sonuçlar verebilmektedir. İkincisi, zaman ilerledikçe, fanın içindeki fiyat hareketinin mesafesi büyümekte, yakın dirençler ve destekler görülmez olmaktadır.

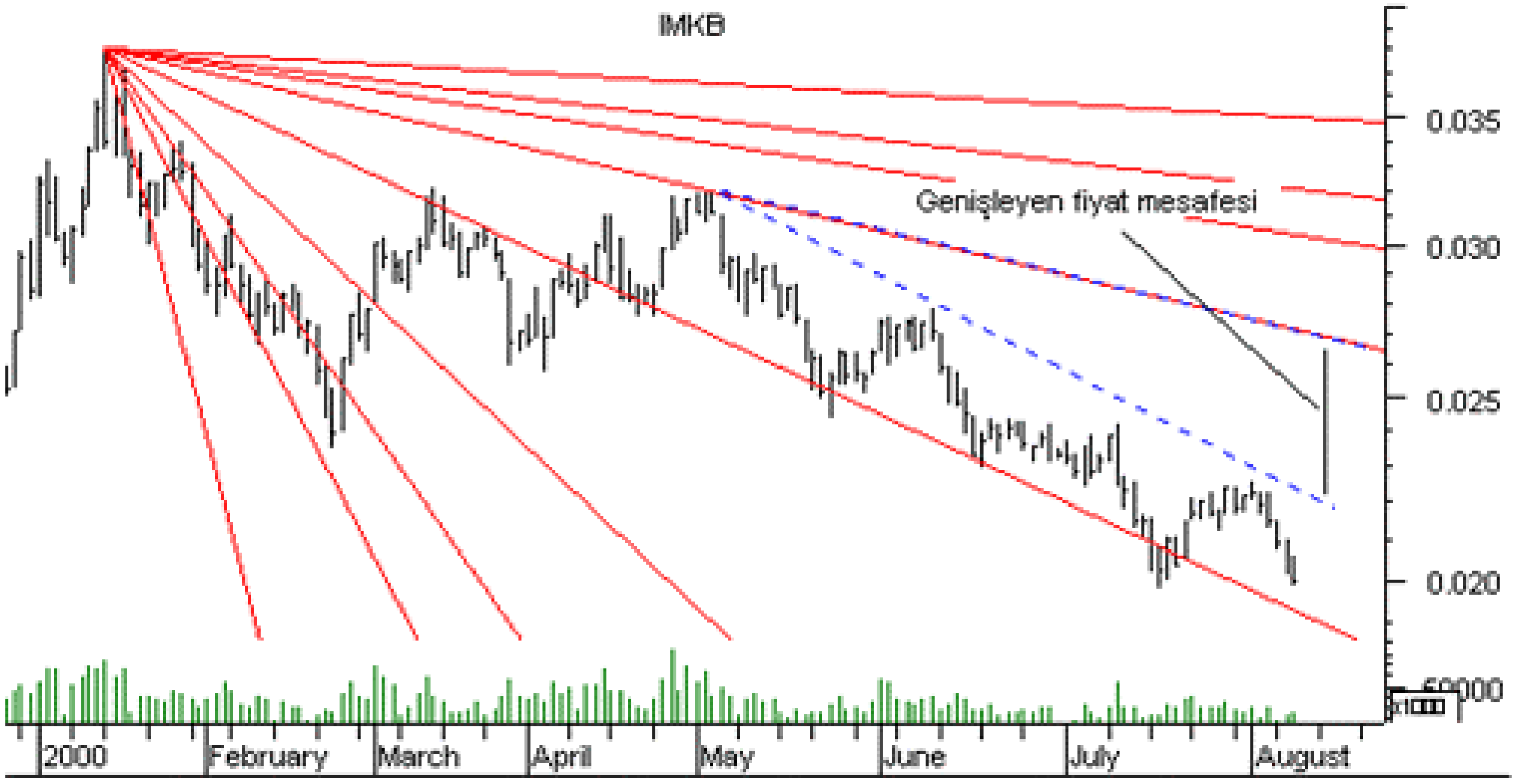

Yukarıda sıralanan yöntemlerin her birinin temel varsayımı, bir ana (ya da belirgin) dipten itibaren fiyatı izlemek esasına dayanır. Dönüş noktalarının, yatırımcı psikolojisini değiştirdiği ve daha önceki döneme ait duygu ve düşüncelerin yerlerini yeni döneme ait beklentilere bıraktığı varsayımı ile yeni dönemde fiyatların nasıl bir trend üzerinde ilerleyebileceği ve hangi seviyelerde destek ve dirençle karşılaşabileceğini tespit etmeye yönelik olan bu varsayımlarda sadece belirgin dibe (ya da tepeye) odaklanılmış olması, yön değiştikten sonra gelişen fiyat hareketinin ne kadarlık bir zaman içinde ilerleyeceği sorusunu havada bırakmaktadır.

**4.6 Zaman ve Fiyat hedefleri**

Bir fiyat hareketinin katettiği fiyat mesafesi ve zaman süresi, kendisinden sonra gelen fiyat hareketinin fiyat ve zaman mesafesini de belirler. Birbirini takip eden iki fiyat hareketi arasında fibonacci oranları kadar bir ilişki olduğu varsayılır.

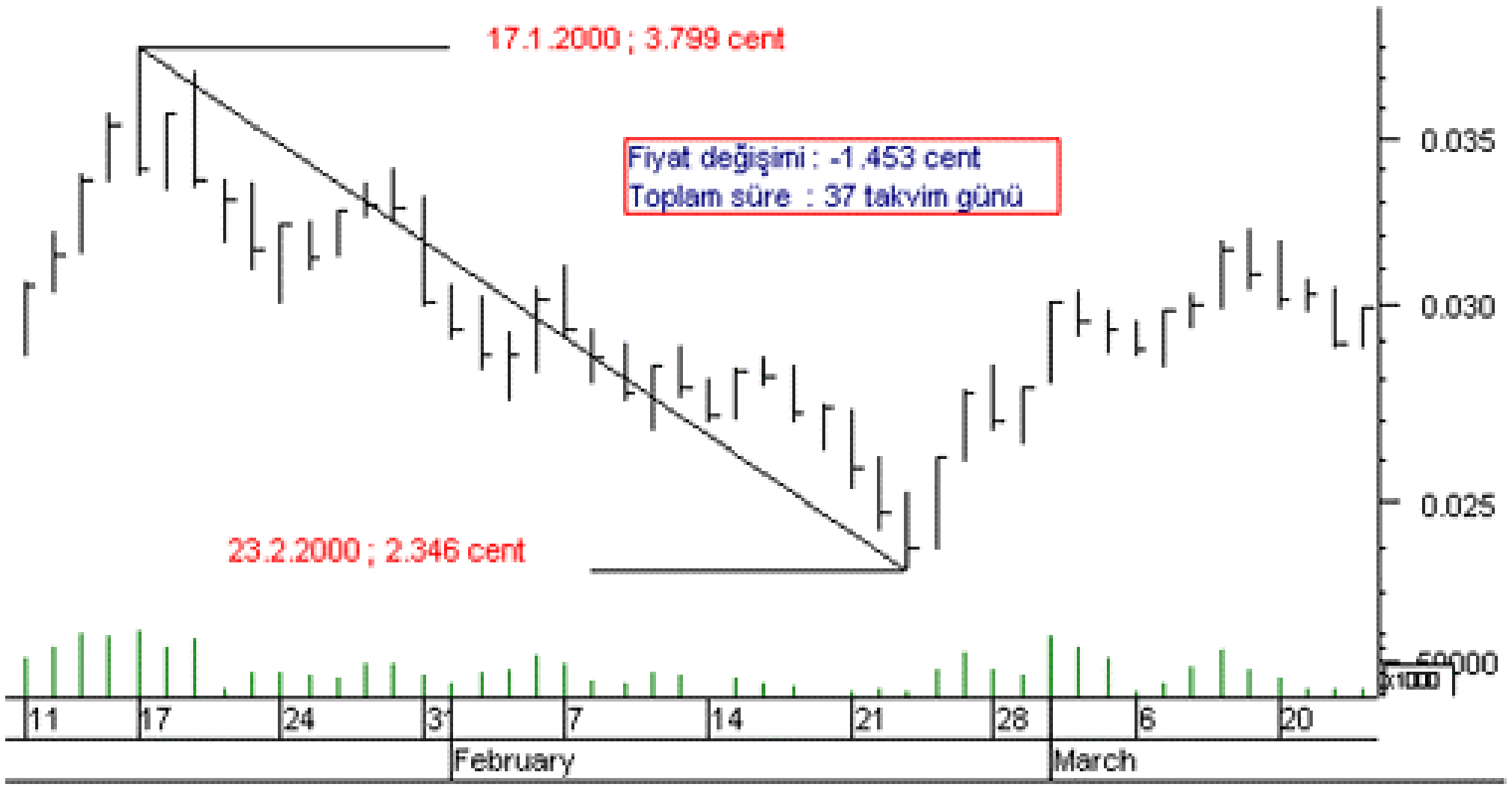

İMKB endeksi, 17 Ocak'tan, 23 Şubat'a kadar devam eden düşüşte, 3.799 cent'ten, 2.346 cent'e gerilemiş ve bu düşüş toplam 37 takvim günü sürmüştü. Bu düşüşü izleyecek olan fiyat hareketinin fiyat ve zaman hedefleri, şu şekilde hesaplanabilir

Yukarıdaki tabloda hesap yöntemi gösterilen zaman ve fiyat hedefleri, grafiğe yerleştirilip, ana tepe ve ana dip ile uç zaman ve fiyat hedeflerinin çakışma noktaları birleştirildiğinde, aşağıdaki grafik elde edilir :

Bu grafikte çizilen her bir trend çizgisi, bir fiyat-zaman doğrultusunu göstermektedir ve birbirinden farklı düşünen ve piyasaya farklı büyüklüklerde etkide

| Başlangıç Günü (İG) | Bitiş Günü (SG) | Başlangıç Fiyatı (IF) | Bitiş Fiyatı (SF) | Toplam takvim günü (DT) | Fiyat değişimi (DP) | | |
|---|---|---|---|---|---|---|---|
| 17.1.2000 | 23.2.2000 | 3.799 | 2.346 | 37 | -1.453 | | |
| Geri Alış | | | | Eşitlik | Uzatma | | |
| 0.236 | 0.382 | 0.500 | 0.618 | 1.000 | 1.618 | 2.618 | 4.236 |
| **Zaman Hedefleri** | | | | | | | |
| | | | | T1 | T2 | T3 | T4 |
| | | | | DT * 1.0 +SG | DT * 1.618 +SG | DT * 2.618 +SG | DT * 4.236 +SG |
| | | | | 31.3.2000 | 22.4.2000 | 29.5.2000 | 28.7.2000 |
| **Fiyat Hedefleri** | | | | | | | |
| P1 | P2 | P3 | P4 | P5 | | | |
| SF – DP * 0.236 | SF– DP * 0.382 | SF – DP * 0.500 | SF – DP * 0.618 | SF – DP * 1.00 | | | |
| 2.689 | 2.901 | 3.073 | 3.244 | 3.799 | | | |

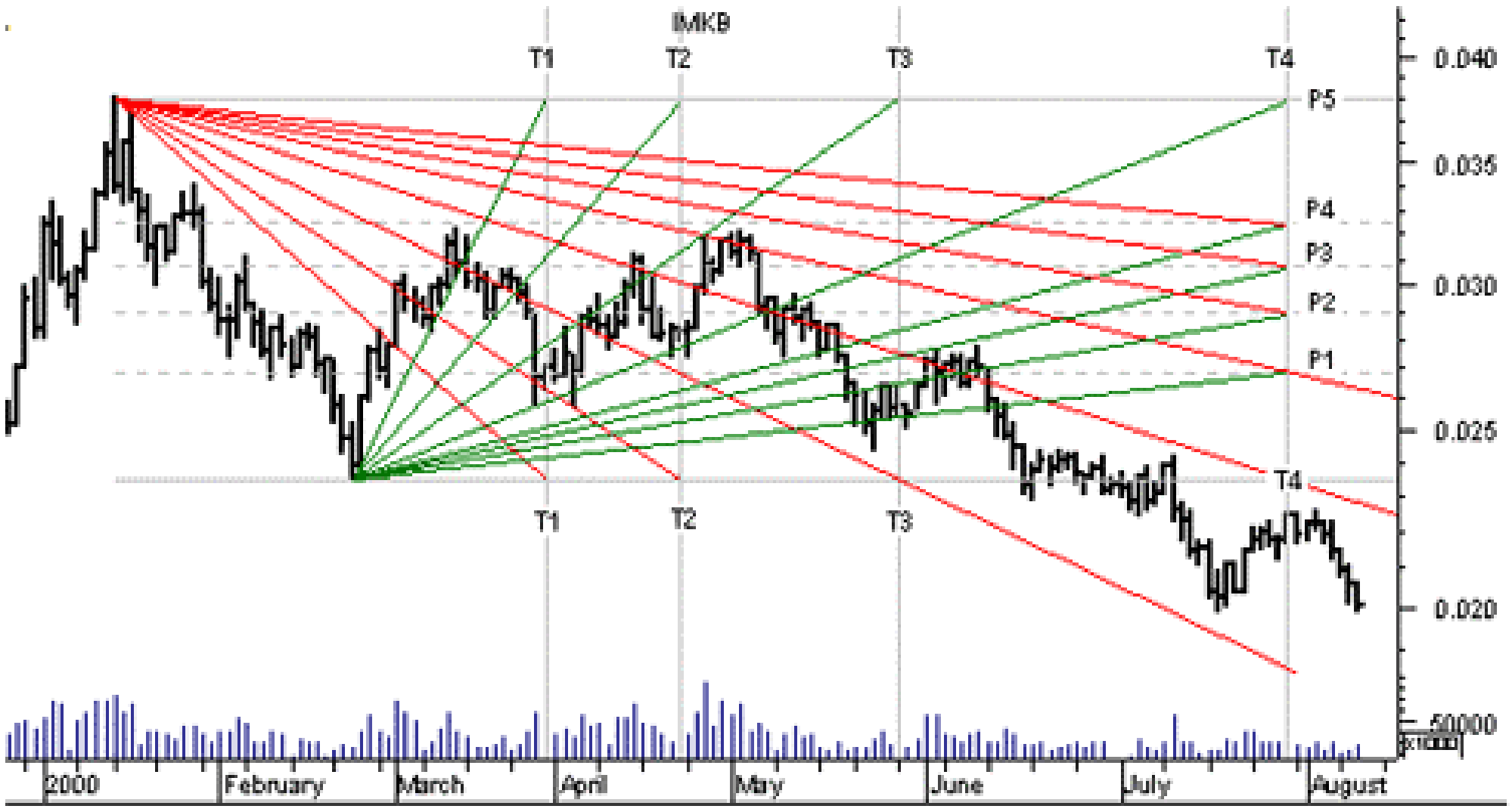

bulunan binlerce yatırımcının farklı eğilimlerini temsil etmektedir. Düşüşe ilk tepki, fiyatı hızla T1-P5 zaman-fiyat hedefine çekme eğilimini temsil eden trend çizgisi üzerinde ilerledi. Bu trend çizgisi kırılmasa idi, düşüşte verilen miktar, en geç düşüş zamanı içinde (yani 37 takvim gününde) geri alınacaktı. Fibonacci box'ın yorumlanışı ile temel ilkeleri koyduktan sonra, alım satımda nasıl kullanılacağını da tespit edelim :

1. Fibonacci box'ı oluşturan fiyat-zaman trend çizgilerinin birbirlerini keserek oluşturduğu alanlar, fiyat bölgeleridir ve fiyat bu bölgeler içinde hareket etme eğilimindedir. Bir fiyat bölgesi aşılıp diğerine geçildiğinde, yeni bölge içinde hareketin devam edeceği varsayılabilir.

2. Fiyat-zaman trend çizgileri, fibonacci box'ın dışına taştıktan sonra da etkilerini sürdürürler.

3. Fiyat, alçalan ve yükselen fiyat-zaman çizgilerinin kesişme noktalarında bulunma eğilimindedir. Kesişme noktaları, konsensusu temsil eder ve fiyat, içinde bulunduğu fiyat-zaman bölgesinde kesişme noktalarından birine gitme eğilimindedir.

4. Fibonacci box'ın diyagonalini oluşturan trend çizgileri, (Ana tepeden ve ana dipten T4'e yönelmiş olan trend çizgileri), tüm trend çizgileri içinde en önemli olanlarıdır ve en uzun süre etkili olmaları beklenir.

5. Kitle psikolojisi, fiyatı, fiyat-zaman hedeflerine çekme eğilimindedir.

Yukarıdaki ilkelerden hareketle, 23 Şubat'taki dibin bir pivot noktası olduğu belirginleştiğinde fibonacci box çizilebilir. (Bir fiyat dibinin pivot noktası olduğu, bazı dönüş formasyonları ile tespit edilebilir.) Fibonacci box'ın T1'e yönelmiş olan ilk trend çizgisi kırılana kadar pozisyonlar sürdürülmeli, bu trend çizgisi kırıldığında satış yapılmalıdır. Gerileyen fiyatın, diyagonal trend açısında destek bulması beklenir ve bu trend çizgisine ulaştığında yeniden alım yapılır. Bu trend çizgisi üzerinde kalındığı sürece pozisyonlar sürdürülür ve bu trend çizgisinin kırılması önemli bir dönüş işareti olduğundan, trend çizgisi kırıldığında satış yapılır. Fiyat fibonacci box'ın son yükselen trend çizgisi üzerinde de tutunamadığında, yani 2.60 civarında son yükselen trend çizgisini de kırdığında, artık fibonacci box'ın dışına çıklmakta olduğu kesinleştiğinden kapatılmamış pozisyonlar da kapatılır ve fiyat hareketi daha farklı bir fibonacci box üzerinde incelenir.

Aşağıdaki grafik, İMKB endeksinin 26 Ağustos - 17 Ocak yükselişi baz alınarak çizilmiş olan fibonacci box'ı göstermektedir. Tüm düşüş döneminde bir ana çakışma noktasında bulunma eğiliminde olan fiyat, bu çalışmanın yapıldığı gün de hızla trend çakışma hedefine yönelmişti. En önemli trend desteğine ulaşmış görünen endekste bu trend çizgisinin bundan sonraki dönemde etkili olması beklenebilir. Çok kısa zaman içinde, aşağı ve yukarı eğimli iki diyagonal trend çizgisinden birini kıracak olan endekste, hangi diyagonalin kırıldığı, orta vadeli yönü de belirleyecektir.

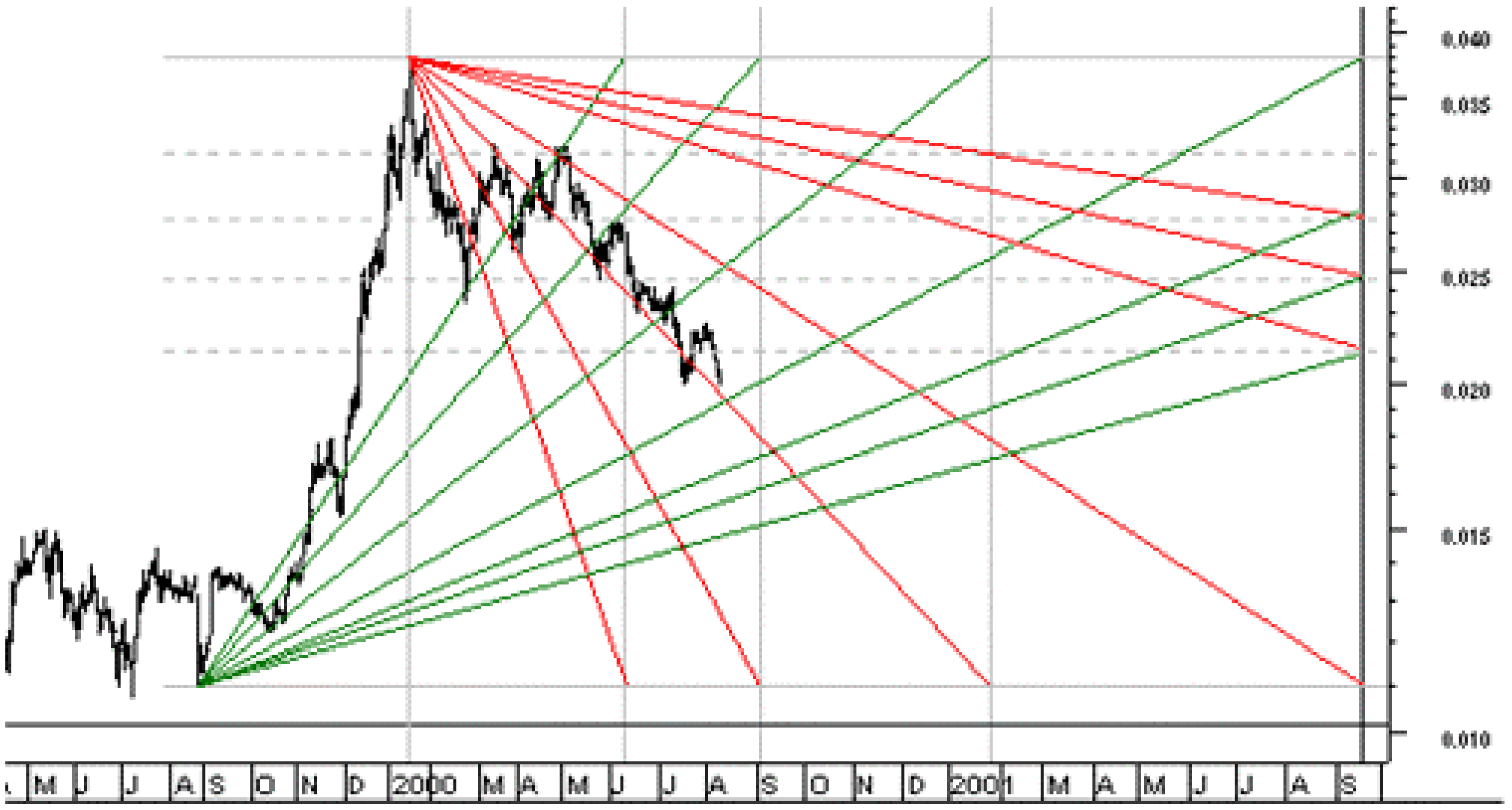

## EK-1: İLK 50 FİBONACCİ SAYISI VE ÇARPANLARINA AYRILMIŞ HALİ

0 : 0

1 : 1

2 : 1

3 : 2

4 : 3

5 : 5

6 : 8 = 23

7 : 13

8 : 21 = 3 x 7

9 : 34 = 2 x 17

10 : 55 = 5 x 11

11 : 89

12 : 144 = 24 x 32

13 : 233

14 : 377 = 13 x 29

15 : 610 = 2 x 5 x 61

16 : 987 = 3 x 7 x 47

17 : 1597

18 : 2584 = 23 x 17 x 19

19 : 4181 = 37 x 113

20 : 6765 = 3 x 5 x 11 x 41

21 : 10946 = 2 x 13 x 421

22 : 17711 = 89 x 199

23 : 28657

24 : 46368 = 25 x 32 x 7 x 23

25 : 75025 = 52 x 3001

26 : 121393 = 233 x 521

27 : 196418 = 2 x 17 x 53 x 109

28 : 317811 = 3 x 13 x 29 x 281

29 : 514229

30 : 832040 = 23 x 5 x 11 x 31 x 61

31 : 1346269 = 557 x 2417

32 : 2178309 = 3 x 7 x 47 x 2207

33 : 3524578 = 2 x 89 x 19801

34 : 5702887 = 1597 x 3571

35 : 9227465 = 5 x 13 x 141961

36 : 14930352 = 24 x 33 x 17 x 19 x 107

37 : 24157817 = 73 x 149 x 2221

38 : 39088169 = 37 x 113 x 9349

39 : 63245986 = 2 x 233 x 135721

40 : 102334155 = 3 x 5 x 7 x 11 x 41 x 2161

41 : 165580141 = 2789 x 59369

42 : 267914296 = 23 x 13 x 29 x 211 x 421

43 : 433494437

44 : 701408733 = 3 x 43 x 89 x 199 x 307

45 : 1134903170 = 2 x 5 x 17 x 61 x 109441

46 : 1836311903 = 139 x 461 x 28657

47 : 2971215073

48 : 4807526976 = 26 x 32 x 7 x 23 x 47 x 1103

49 : 7778742049 = 13 x 97 x 6168709

50 : 12586269025 = 52 x 11 x 101 x 151 x 3001

## EK-2: FİBONACCİ SAYILARINI VEREN ALGORİTMA VE PROGRAM

Algorithm Fibonacci (n)

(* Bu algoritma n.inci Fibonacci sayısını indirgeme bağıntısını kullanarak hesaplar.*)

Begin(*algorithm*)

if n=1 or n=2 then (* base cases *)

Fibonacci → 1

else (* general case *)

Fibonacci → Fibonacci (n-1) + Fibonacci (n-2)

End (* algorithm *)

```
//Fibonacci Sayilari
#include<stdio.h>
void main(void)
{
int a, b, c, dongu_sayisi, index;
printf("Kac tane Fibonacci sayisi gormek istiyorsunuz:");
scanf("%d", &dongu_sayisi);
while(dongu_sayisi<=0)
{
printf("Lutfen pozitif bir sayi giriniz:");
scanf("%d", &dongu_sayisi);
}
a=0;
b=1;
index=0;
do
{
index++;
c=a+b;
a=c-a;
b=c;
```

```
printf("%d ", c);
} while(index<dongu_sayisi);
}
```

# KAYNAKLAR :

# ÖZGEÇMİŞ

Erdoğan ŞEN, 28.12.1985 tarihinde Bulgaristan'ın Kırcaali kentinde doğdu. Orta öğrenimini Pertevniyal Anadolu Lisesi 'nde tamamladıktan sonra 2004 yılında Gebze Yüksek Teknoloji Enstitüsü Fen Fakültesi Matematik Bölümü 'nü kazanıp halen matematik bölümü öğrencisi olarak çalışmalarına devam etmektedir.